\documentclass[a4paper,12pt]{article}
\usepackage{amsmath}
\usepackage[mathscr]{eucal}
\usepackage{amssymb}
\usepackage{latexsym}
\usepackage{amsthm}
\usepackage{amscd}
\theoremstyle{plain}
\newtheorem{theorem}{Theorem}%[section]

\newtheorem{lemma}{Lemma}[section]
\newtheorem{proposition}{Proposition}

\theoremstyle{remark}
\newtheorem{remark}{Remark}
\theoremstyle{definition}

\usepackage[dvips]{graphicx}
\usepackage{psfrag}
\begin{document}
\title{ Remarks on surfaces with $c_1^2 =2\chi -1$ \\ 
having non-trivial $2$-torsion}
\author{Masaaki MURAKAMI%\footnote{} 
%\medskip \\
%University of Bayreuth\\
%Lehrstuhl Mathematik VIII \\
%Universitaetsstrasse 30 \\
%D-95447 Bayreuth, Germany\\
%\texttt{Masaaki.Murakami@uni-bayreuth.de} 
}
\date{}
\maketitle
\begin{center}
On the occasion of $60$-th birthday of Prof.\,Fabrizio Catanese
\end{center}
\begin{abstract}
We shall show that any complex minimal surface of general type 
with $c_1^2 = 2\chi -1$ having non-trivial $2$-torsion divisors, 
where $c_1^2$ and $\chi$ are the first Chern number of a surface 
and the Euler characteristic of the structure sheaf respectively,  
has the Euler characteristic $\chi$ not exceeding $4$. 
Moreover, we shall give a complete description for the surfaces of the 
case $\chi =4$, 
and prove that the coarse moduli space for surfaces of this case 
is a unirational variety of dimension $29$. 
Using the description, we shall also prove that our surfaces of  
the case $\chi = 4$ have non-birational bicanonical maps and no pencil 
of curves of genus $2$, hence being of so called non-standard case for 
the non-birationality of the bicanonical maps. 
\end{abstract}

\section{Introduction}  \label{sectn:introduction}

In classification of regular surfaces of general type, 
the torsion parts of the Picard groups (the torsion groups for short) 
sometimes play an important role. 
One of the reasons for this lies in variety of topological types
under single values of numerical invariants, 
which is common especially in cases of small geometric genus;  
the torsion group of a regular surface, isomorphic to 
the first homology group with integral coefficients, 
carries information that the numerical invariants $c_1^2$ and $\chi$ do not.

Studies on surfaces of general type done using   
the torsion groups are well-known   
for cases of  vanishing geometric 
genus (see, e.g., 
Barth-Peters-Van de Ven \cite[p.\ 237]{complexsurf}). 
In those studies, they tried to determine 
the structures of surfaces with given isomorphism classes of 
the torsion groups. 
There are, however, some other cases of numerical invariants for which 
similar studies have been successfully developed. Consider the case 
$c_1^2 = 2\chi -2$. In this case, by Ciliberto-Mendes Lopes \cite{on2chi-2}, 
the orders of the torsion groups 
do not exceed $2$, and the Euler characteristics $\chi$'s for the 
cases of non-trivial torsion do not exceed $5$.      
Complete descriptions for the surfaces with 
non-trivial torsion with $\chi =2$, $3$, $4$, and $5$ 
are given in  Catanese-Debarre \cite{pg1k2=2'}, 
Ciliberto-Mendes Lopes \cite{on2chi-2}, 
Bartalesi-Catanese \cite{withtorsion}, and 
Ciliberto-Mendes Lopes \cite{on2chi-2} respectively. 
We remark that even in cases of vanishing geometric genus, 
complete descriptions are known only for a small number of 
classes.

%%%
%In the present paper, we study minimal surfaces with  
%$c_1^2= 2 \chi -1$ having non-trivial $2$-torsion divisors:  
%we shall give a bound for the Euler 
%characteristics $\chi$'s (Theorem \ref{thm:maintheorem}), 
%describe the surfaces of 
%the case of maximal $\chi$ (Theorem \ref{thm:completedescription},
%Remark \ref{rem:meaning}), 
%and study the moduli space for surfaces of this case  
%(Theorem \ref{thm:moduli}). 
%By the main theorem of \cite{bound'''}, 
%the order of the torsion group 
%of a minimal surface with $c_1^2 = 2\chi -1$ is 
%at most $3$ if $\chi =2$, and at most $2$ if $\chi \geq 3$.  
%Thus for our surfaces with $\chi \geq 2$, two conditions 
%$\mathbb{Z}/2 \subset \mathrm{Tors}$ and   
%$\mathrm{Tors} \simeq \mathbb{Z}/2$ are equivalent, 
%where $\mathrm{Tors}$ denotes the torsion group. 
%The case $\chi =1$ on this line is that of 
%the numerical Godeaux surfaces (i.e., minimal surfaces 
%of general type with $c_1^2=1$ and $p_g= 0$).  
%%%

%%
In the present paper, we study minimal surfaces with  
$c_1^2= 2 \chi -1$ having non-trivial $2$-torsion divisors.
Note that if $X$ is a minimal surface with 
$c_1^2 = 2 \chi -1$, then $X$ has vanishing irregularity, 
hence geometric genus $p_g = \chi -1$. 
We shall prove the bound $\chi \leq 4$ for the Euler 
characteristics $\chi$'s (Theorem \ref{thm:maintheorem}), 
describe the surfaces of 
the case $\chi = 4$ (Theorem \ref{thm:completedescription},
Remark \ref{rem:meaning}), 
and study the moduli space for surfaces of this case  
(Theorem \ref{thm:moduli}). 
By the main theorem of \cite{bound'''}, 
the order of the torsion group 
of a minimal surface with $c_1^2 = 2\chi -1$ is 
at most $3$ if $\chi =2$, and at most $2$ if $\chi \geq 3$.  
Thus for our surfaces with $\chi \geq 2$, two conditions 
$\mathbb{Z}/2 \subset \mathrm{Tors}$ and   
$\mathrm{Tors} \simeq \mathbb{Z}/2$ are equivalent, 
where $\mathrm{Tors}$ denotes the torsion group. 
The case $\chi =1$ on this line is that of 
the numerical Godeaux surfaces (i.e., minimal surfaces 
of general type with $c_1^2=1$ and $p_g= 0$).  

Surfaces with $c_1^2 = 2 \chi -1$, $\chi = 4$, and 
$\mathrm{Tors} \simeq \mathbb{Z} / 2$ are known to exist and 
can be found in \cite{nonstandardpg3}. 
In \cite{nonstandardpg3}, Ciliberto and Mendes Lopes 
completely classified regular surfaces with $p_g = 3$ having 
non-birational bicanonical maps and without genus $2$ 
pencils, i.e., regular surfaces with $p_g = 3$ and of non-standard
case for the non-birationality of the bicanonical maps.  
Among their results, they showed that 
any regular surface of non-standard 
case with $c_1^2 = 7$ and $p_g =3$ 
is obtained by performing a certain operation on what is known 
as Du Val's ancestor with $c_1^2=8$ and $p_g = 4$. Since  
these surfaces have non-trivial $2$-torsion divisors, 
as has been shown in \cite{nonstandardpg3}, 
these are examples of our surfaces for the case $\chi =4$.  
In fact, our structure theorem for surfaces with $c_1^2 = 2 \chi -1$, 
$\chi = 4$, and $\mathrm{Tors} \simeq \mathbb{Z} / 2$ 
shows that although we start from the different assumption, 
the resulting surfaces are exactly those seen in 
the paper \cite{nonstandardpg3}.

Our complete description for the surfaces with $\chi =4$ 
asserts that any such surface $X$ is obtained roughly as 
a free quotient by $\mathbb{Z} /2$ 
of a double cover of the 
Hirzebruch surface 
$\varSigma_d = \mathbb{P} (\mathcal{O}_{\mathbb{P}^1} \oplus 
\mathcal{O}_{\mathbb{P}^1} (d))$ ($d=0$ or $2$).   
We shall describe the branch divisor of the double cover,  
and determine the free action by $\mathbb{Z} /2$ 
(Theorem \ref{thm:completedescription}, 
Remark \ref{rem:involutiondescription}).   
The branch divisor of the double cover 
turns out to be a member of the 
quadruple anticanonical system having exactly two $[3,3]$-points.     
The action by $\mathbb{Z} /2$ turns out to be a lifting of 
that on the Hirzebruch surface $\varSigma_d$.  
This description induces another description 
of our surfaces of the case $\chi =4$ 
(Proposition \ref{prop:anotherdescription}),  
which is almost the same as a description appearing   
in Ciliberto-Mendes Lopes \cite{nonstandardpg3}. 
Using our descriptions, we shall show that our 
surfaces of the case $\chi = 4$ has non-birational 
bicanonical maps and no pencil of curves of genus $2$
(Proposition \ref{prop:bicanonical}), 
hence completely coinciding   
with those seen in \cite{nonstandardpg3}
(see also Remark \ref{rem:notecmdescr}).     

The coincidence of the resulting surfaces 
certainly implies possibility of another proof of our complete 
description, i.e., of a proof, 
like one for the case $c_1^2 = 2 \chi -2$ in 
Ciliberto-Mendes Lopes \cite{on2chi-2},  
by showing that our surfaces with $\chi =4$ are 
of non-standard case for the non-birationality of 
the bicanonical maps. 
We however do not chose this way. 
We remark that our method has an advantage in the 
sense that we can show the irreducibility of the 
moduli space in a very explicit and elementary way.

The present paper is organized as follows. 
In order to show our main theorem, 
we follow Miyaoka \cite{MiyaokaTri} and Reid \cite{Sfpgk2'},  
and take the unramified double cover 
$Y \to X$ corresponding to a torsion divisor. 
We study its canonical map $\varPhi_{K_Y}$ using 
the action by the Galois group of $Y$ over $X$. 
In Section \ref{scn:statement}, we state our main results and 
show, on the assumption $\chi \geq 4$, that we have  
$\deg \varPhi_{K_Y} =1$ or $2$, 
and that $\deg \varPhi_{K_Y} =1$ implies $\chi =4$.  
Note here that to obtain our main theorem, we only need to study the 
case $\chi \geq 4$. 
In Section \ref{scn:deg=2}, we study the case $\deg \varPhi_{K_Y} =2$. 
We divide this case into three according to the degree of the 
canonical image $Z = \varPhi_{K_Y} (Y) \subset \mathbb{P}^n$: 
the case $\deg Z = n+1$, the case $\deg Z =n$, and the case $\deg Z = n-1$.  
We shall classify non-degenerate surfaces 
in $\mathbb{P}^n$ of degree $n+1$  
of which minimal desingularizations have vanishing irregularities 
(Proposition \ref{prop:deg=n+1}), 
and use this classification to study the case $\deg Z = n+1 $.    
%We employ Horikawa's method to study the canonical map $\varPhi_{K_Y}$.
In Section \ref{scn:deg=1}, 
we study the case $\deg \varPhi_{K_Y} =1$ 
and $\chi =4$, 
and then prove Theorems \ref{thm:maintheorem} 
and \ref{thm:completedescription}.  
In the case $\deg \varPhi_{K_Y} =1$ and $\chi =4$,  
the surface $Y$ has the first Chern number $14$, 
geometric genus $7$, and irregularity $0$.  
Hence the surface $Y$ in this case is a canonical 
surface whose invariant lies on the Castelnuovo line.
We use results given in Ashikaga-Konno \cite{3pg-7}
%, a paper by Ashikaga and Konno, 
to exclude this case. 
Finally in Section \ref{scn:onmodulispace}, 
we study the coarse moduli space for the surfaces of the 
case $\chi =4$, and prove Theorem \ref{thm:moduli}.  
To prove the unirationality of the moduli space and the uniqueness 
of the deformation type, we describe our surfaces of the case 
$\chi =4$ as double planes, which is 
almost the same as the description in 
Ciliberto-Mendes Lopes \cite{nonstandardpg3} 
for the surfaces of the non-standard case 
(see also Ciliberto-Francia-Mendes Lopes \cite{onbicanonicalmaps}).  
Using the two descriptions of our surfaces, we show that 
our surfaces of the case $\chi =4$ in fact are 
of the non-standard case 
for the non-birationality of bicanonical maps. 
\medskip

{\sc Acknowledgment}

The author performed the first half of the computations 
included in this article 
during his stay at Max-Planck-Institut f\"{u}r Mathematik in Bonn   
from October 2004 until March 2005. 
He expresses his deepest gratitude to the institute for the comfortable 
environment and the financial support which he received during his stay. 
During the final preparation of the manuscript, the author was supported 
by the research grant fellowship by the University of Padua, 
and then later by DFG Forschergruppe 790 ``Classification of 
algebraic surfaces and compact complex manifolds''. 
The author is thankful to the referee, who suggested him the 
proof of Proposition \ref{prop:bicanonical}, shorter than the proof 
in the earlier version. 
The author is thankful also to Prof.\,Kazuhiro Konno for letting him 
know the paper \cite{evenI}, to Prof.\,Margarida Mendes Lopes for her 
kind comments on the first manuscript, 
and to Prof.\,Fabrizio Catanese and Prof.\,Ingrid Bauer for giving 
him comfortable environment for the study at Bayreuth, where he 
performed the final revision of the manuscript.

\medskip

{\sc Notation and Terminology}

Let $S$ be a compact complex manifold of dimension $2$. 
We denote by $c_1 (S)$, $p_g(S)$, and $q(S)$ the 
first Chern class, the geometric genus, and the irregularity 
of $S$ respectively. 
The torsion group of $S$, denoted by $\mathrm{Tors} (S)$, is 
the torsion part of the Picard group of $S$. 
If $V$ is a complex manifold, $K_V$ is a canonical divisor of $V$. 
For a coherent sheaf $\mathcal{F}$ on $V$, 
we denote by $H^i (\mathcal{F})$, $h^i (\mathcal{F})$, and 
$\chi (\mathcal{F})$ the $i$-th cohomology group, its dimension 
$\dim_{\mathbb{C}} H^i (\mathcal{F})$, and 
the Euler characteristic $\sum (-1)^i h^i (\mathcal{F})$ respectively.  
Let $f: V \to W$ be a morphism to a complex manifold $W$, and 
$D$, a divisor on $W$. Then $f^* (D)$ and $f^{-1}_* (D)$  
denote the total transform and the strict transform respectively of $D$. 
The symbol $\sim$ means the linear equivalence of divisors. 
We denote by $\varSigma_d \to \mathbb{P}^1$ 
the Hirzebruch surface of degree $d$. 
The divisors $\varDelta_0$ and $\varGamma$ are 
its minimal section and its fiber respectively. 
Let $C$ be a curve on $S$. We denote by 
$\mathrm{mult}_x\, C$ the multiplicity of $C$ at a point $x \in S$. 
Let $x$ be a triple point of a reduced curve $C$ on $S$, 
and $S^{\prime} \to S$, the blowing-up at $x$. Assume that 
the strict transform $C^{\prime}$ of $C$ has 
an infinitely near triple point $x^{\prime}$.     
Then the point $x$ is called a $[3, 3]$-point of $C$, 
if the strict transform $C^{\prime \prime}$ to $S^{\prime \prime}$, 
where $S^{\prime \prime} \to S^{\prime}$ 
is the blowing-up at $x^{\prime}$, 
has at most negligible singularities on the 
exceptional locus of $S^{\prime \prime} \to S$.

\section{Statement of the main theorem}  \label{scn:statement}

In \cite{bound'''}, we obtained a bound for the orders of the torsion 
groups of minimal surfaces with $c_1^2 = 2\chi -1$ and $\chi \geq 2$. 
In the present paper, we study the case of $2$-torsion divisors, and 
sharpen the bound.  
Our goals are a bound for the Euler characteristic $\chi$,   
a complete description for the surfaces of the case of maximal $\chi$, 
and the unirationality of the moduli space for surfaces of this case. 
The following three are the main theorems:   

\begin{theorem} \label{thm:maintheorem}
Let $X$ be a minimal surface of general type with 
$c_1^2 =2\chi -1$ and torsion group 
$\mathrm{Tors} (X) \simeq \mathbb{Z} /2$. 
Then the Euler characteristic $\chi$ 
of the structure sheaf does not exceed $4$. 
\end{theorem}

\begin{theorem} \label{thm:completedescription}
Let $X$ be a minimal surface with $c_1^2 = 2\chi -1$, $\chi =4$, 
and torsion group $\mathrm{Tors} (X) \simeq \mathbb{Z} /2$. 
Then the unramified double cover $Y$ of $X$ admits a 
generically two-to-one morphism $f$ onto the 
Hirzebruch surface $\varSigma_d$ of degree $d=0$ or $2$ 
satisfying the following conditions$:$ %\smallskip

%\noindent
$\mathrm{i}$%
$)$ the action by the Galois group 
    $G = \mathrm{Gal} (Y/X) \simeq \mathbb{Z} /2$ 
    of $Y$ over $X$ induces one on $\varSigma_d$,  
    of which fixed locus is a set of four points on $\varSigma_d$%
$;$  

%\noindent
$\mathrm{ii}$%
$)$ the branch divisor $B$ of $f$ is a member of the 
     linear system $|-4K_{\varSigma_d}|$ passing no fixed 
     point of the action by $G$%
$;$   

%\noindent
$\mathrm{iii}$%
$)$ the branch divisor $B \in |-4K_{\varSigma_d}|$ 
      has exactly two $[3,3]$-points, and all other singularities, 
      if any, are negligible ones.
\end{theorem}

\begin{theorem} \label{thm:moduli}
Any two minimal surfaces with $c_1^2 = 2 \chi -1$, $\chi =4$, and 
$\mathrm{Tors} \simeq \mathbb{Z}/2$ are equivalent under deformation 
of complex structures.  The coarse moduli space for minimal surfaces 
with these invariants is a unirational variety of dimension $29$.  
\end{theorem}

Theorem \ref{thm:maintheorem} sharpens the bound given in \cite{bound'''} 
into the following: 

\begin{theorem} \label{thm:modifiedbound}
Let $X$ be a minimal algebraic surface  
with $c_1^2=2\chi -1$. 
Then the following hold$:$

$\mathrm{i}$%
$)$ if $\chi =2$, then $\sharp \mathrm{Tors} (X) \leq 3$%
$;$ 

$\mathrm{ii}$%
$)$ if $\chi \geq 3$, 
then $\sharp \mathrm{Tors} (X) \leq 2$%
$;$ 

$\mathrm{iii}$%
$)$ if $\chi \geq 5$, 
then $\sharp \mathrm{Tors} (X) =1$.  
\end{theorem}

\begin{remark}  \label{rem:involutiondescription}
In Theorem \ref{thm:completedescription}, 
we can describe the action by $G$ on $\varSigma_d$ more concretely:  
if an involution of the Hirzebruch surface $\varSigma_d$ 
has exactly four fixed points ($d$: even),  
then there exists an open cover $\{ U_i \}_{i=0, 1}$ of 
$\varSigma_d$ satisfying 
$U_i = \{ (u_i, (t_i : 1))\} 
= \mathbb{C} \times \mathbb{P}^1$, 
$u_0 = 1/ u_1$, and 
$t_0 = u_1^d t_1$,  
such that this involution is given by  
\begin{equation} \label{eq:involution}
 (u_0, t_0) \mapsto (-u_0, -t_0). 
\end{equation}

\end{remark}

\begin{remark} \label{rem:meaning}
Theorem \ref{thm:completedescription} 
asserts that any minimal surface $X$ 
with $c_1^2 =2\chi -1$, $\chi =4$, 
and $\mathrm{Tors} (X) \simeq \mathbb{Z} /2$ is obtained 
by the following procedure: 
$1$) set $d=0$ or $2$; the involution (\ref{eq:involution}) 
defines an action by $G= \mathbb{Z} /2$ on the 
Hirzebruch surface $\varSigma_d$;  
$2$) take a reduced member $B \in |-4K_{\varSigma_d}|$ 
stable under this action that satisfies 
the conditions ii) and iii) in Theorem \ref{thm:completedescription};   
$3$) take the double cover of $\varSigma_d$ branched along $B$, and 
denote by $Y$ its minimal desingulraization;  
there exists a unique free lifting to $Y$ of the action by $G$ on 
$\varSigma_d$; 
$4$) take the quotient of $Y$ by this free action. 

It is not difficult to check that this procedure in fact gives 
surfaces of the case $\chi =4$ for sufficiently general $B$.  
\end{remark}

\begin{remark}
Let $\varSigma_d$ be the Hirzebruch surface which appears in 
Theorem \ref{thm:completedescription}.
It is obvious from Rmarks \ref{rem:involutiondescription} 
and \ref{rem:meaning} 
that the fibration $\varSigma_d \to \mathbb{P}^1$
induces a hyperelliptic fibration $Y \to \mathbb{P}^1$ of genus $3$ 
and that the divisor class of a fiber of this fibration is stable under 
the action by the Galois group $G = \mathrm{Gal} (Y / X)$.
So we obtain a hyperelliptic fibration $X = Y/G \to \mathbb{P}^1 / G$ 
of genus $3$ with two multiple fibers $2 A_1$ and $2 A_2$ 
corresponding to the fixed points 
of the action by $G$ on $\mathbb{P}^1$.   
As is explained also in \cite[p.\,$85$]{nonstandardpg3}, 
the difference $A_1 - A_2$ gives a non-trivial $2$-torsion divisor 
of our surface $X$. 
\end{remark}

In what follows, $X$ is a minimal surface with $c_1^2 = 2\chi -1$, 
$\chi = \lambda \geq 4$, and $\mathrm{Tors} (X) \simeq \mathbb{Z} /2$. 
We denote by $\pi : Y \to X$ the unramimfied double cover 
corresponding to the torsion group $\mathrm{Tors} (X)$. 
Note that we have assumed $\lambda \geq 4$. 
The following lemma follows from 
%Deligne's argument \cite[Theorem 14]{canonicalmdl} and 
the unbranched covering trick. 

\begin{lemma}
$K_Y^2 = 2(2\lambda -1)$, $p_g(Y) = 2\lambda -1$, and $q(Y) =0$. 
\end{lemma}

In order to show Theorems \ref{thm:maintheorem} 
and \ref{thm:completedescription}, we study the 
canonical map $\varPhi_{K_Y} : Y \to \mathbb{P}^n$ of $Y$, 
where $n= 2\lambda -2$.
We denote by $Z = \varPhi_{K_Y} (Y)$ 
the canonical image of the surface $Y$. 

\begin{proposition}   \label{prop:degphiKY} 
The canonical image $Z$ is a surface. The equality 
$\deg \varPhi_{K_Y} =1$ or $2$ holds. 
Moreover,  
if $\deg \varPhi_{K_Y} =1$, then $\lambda =4$.  
\end{proposition}

Proof. Since we have assumed $\lambda \geq 4$, we have 
\[
 K_Y^2 -3p_g(Y) = -(2\lambda -1) \leq -7 . 
\]
By this together with $q(Y) =0$ and 
\cite[Theorem 1.1]{smallc1-3}, we see that  
$|K_Y|$ is not composite with a pencil. 
Thus we have 
\[
 \deg \varPhi_{K_Y} \leq \frac{K_Y^2}{\deg Z} 
 \leq \frac{2(n+1)}{n-1} = 2 + \frac{4}{n-1} \leq 2 + \frac{4}{5}, 
\]
hence $\deg \varPhi_{K_Y} \leq 2$. The second assertion follows from 
Castelnuovo's inequality. \qed

If $\lambda =4$, then the Chern invariant of $Y$ is on the 
Castelnuovo line. 
Thus we can use results given in \cite{3pg-7} to study the 
case $\deg \varPhi_{K_Y} =1$. 
   
\section{The case $\deg \varPhi_{K_Y} =2$} \label{scn:deg=2}

In this section, we study the case $\deg \varPhi_{K_Y} =2$. 
We begin with the study of the base locus of the canonical 
system $|K_Y|$. 
Let $|M|$ and $F$ be the variable part and the fixed part of the 
linear system $|K_Y|$. 
We take the shortest composite $p: \tilde{Y} \to Y$  of 
quadric transformations such that the variable part $|L|$ of 
$p^{*} |M|$ is free from base points, and denote by $E$ the 
fixed part of $p^{*} |M|$. 
Then we have $p^*|K_Y| = |L| + E + p^*F$ and 
\begin{equation} \label{eq:KY^2}
  K_Y^2 = L^2 + LE + MF + K_Y F,
\end{equation}
where each term of the right hand is a non-negative integer. 
Note that the eigenvectors of the natural action 
by $G = \mathrm{Gal} (Y/X)$ span the 
the space of global section $H^0 (\mathcal{O}_Y (K_Y))$. 
This implies that the linear systems 
$|K_Y|$, $|M|$, and $F$ are spanned by the pull-backs of divisors on $X$.   
Hence, for example, we have  $MF \equiv 0$ $\mathrm{mod}$ $2$, 
since $\pi : Y \to X$ is of mapping degree $2$. 
In the same way, we obtain  
\begin{equation} \label{eq:mod2}
 L^2 \equiv LE = -E^2 \equiv MF \equiv K_YF \equiv 0 
 \quad \mathrm{mod} \quad 2 
\end{equation}
(for the detail, see \cite[Section 3]{3tors'}).

\begin{proposition} \label{prop:L^2}
Let $M$, $F$, $L$, and $E$ be divisors as above. 
Then one of the following holds$:$

$1$%
$)$ $|K_Y|=|L|$% 
$:$ 
the canonical system $|K_Y|$ is free from base points$;$ 

$2$%
$)$ $L^2 = K_Y^2 -2$, $F=0$, and $LE =2$%
$;$

$3$-$1$%
$)$ $L^2 = K_Y^2 -4$, $F=0$, and $LE =4$%
$;$

$3$-$2$%
$)$ $L^2 = K_Y^2 -4$, $|L|=|M|$, $K_Y F=0$, and $F^2 = -4$.
  
\end{proposition}
Proof. 
First, note that we have  
$L^2 = K_Y^2$, $K_Y^2 -2$, or $K_Y^2 -4$. 
This follows from (\ref{eq:mod2}) and \cite[Lemma 2]{quintic}. 
Second, note that  
\begin{equation} \label{eq:MFmod4}
 MF \equiv 0 \quad \mathrm{mod} \quad 4.
\end{equation}
This follows from the Riemann--Roch theorem, since we have 
$MF = M(M + K_Y) - 2M^2 = M(M + \pi^* K_X) - 2 M^2$, 
$\deg \pi = 2$,    
and $M \sim \pi^*M^{\prime}$ for a certain divisor $M^{\prime}$ on $X$.  
Then the assertion follows from (\ref{eq:KY^2}),  
(\ref{eq:MFmod4}), (\ref{eq:mod2}), and Hodge's index theorem. \qed 

In case $3$-$1$), the number of the base points of $|M|$ cannot be $1$, 
since the action by $G$ on $Y$ has no fixed point. 
Thus in this case, the morphism $p : \tilde{Y} \to Y$ is a composite 
of four quadric transformations. 
In the same way, 
we see that, in case $2$), 
the morphism $p : \tilde{Y} \to Y$  is a blowing-up 
of $Y$ at two distinct points. 
In case $3$-$2$), the divisor $F$ is a sum of two fundamental cycles of 
rational double points. 

We denote by $\varPhi_L : \tilde{Y} \to Z \subset \mathbb{P}^n$ 
the morphism associated with the linear system $|L|$. 
The action by $G$ on $Y$ induces one on $\tilde{Y}$. 
We study the morphism $\varPhi_L$ using this action. 

\subsection{The case $|K_Y|=|L|$}

Let us first exclude case $1$) in Proposition \ref{prop:L^2}. 
In what follows, we assume $|K_Y| = |L|$.
Thus we have $\deg Z = n+1$. 
We shall prove the following proposition in Appendix. 

\begin{proposition} \label{prop:deg=n+1}
Let $n \geq 4$ be an integer, $Z$, a non-degenerate surface 
in $\mathbb{P}^n$ of degree $n+1$, and $Z^{\prime} \to Z$, 
its minimal desingularization.  
Assume that the morphism $Z^{\prime} \to Z$ is given by  
a complete linear system $|D^{\prime}|$      
and that $q(Z^{\prime}) = 0$ holds. 
Then $n$ does not exceed $11$. 
Further, there exist an integer $0 \leq d \leq 3$ and a 
blowing-up $r : Z^{\prime} \to \varSigma_d$ at 
$($possibly infinitely near$)$ $11 - n$ points such that 
the equivalence $D^{\prime} \sim - K_{Z^{\prime}} + r^* \varGamma $ holds. 
Here, the divisor $\varGamma$ is a fiber of the Hirzebruch surface 
$\varSigma_d \to \mathbb{P}^1$.    
\end{proposition}

In our case, we have $n = 2\lambda -2$, $\lambda \geq 4$, and $q (Y) =0$. 
Moreover $Z$ is the canonical image of $Y$. Thus our surface 
$Z = \varPhi_{K_Y} (Y)$ satisfies all the conditions 
in the proposition above.  
It follows that there exist an integer $0 \leq d \leq 3$ and 
a blowing-up $r : Z^{\prime} \to \varSigma_d$ at 
$11-n$ points such that 
the morphism $\varPhi_{D^{\prime}} : Z^{\prime} \to Z$, 
where $\varPhi_{D^{\prime}}$ is a morphism corresponding to 
the complete linear system 
$|D^{\prime}| = | -K_{Z^{\prime}} + r^* \varGamma|$, 
gives the minimal desingularization of $Z$.

\begin{proposition} \label{prop:branchdivisor}
The canonical map $\varPhi_{K_Y} : Y \to Z$ lifts to 
a morphism $f^{\prime} : Y \to  Z^{\prime}$. 
The branch divisor $B^{\prime}$ of $f^{\prime}$ is a member of 
the linear system $|2(2D^{\prime} - r^* \varGamma)|$ 
having at most negligible singularities.   
\end{proposition}

Proof. 
Let us first show the liftability of the canonical map $\varPhi_{K_Y}$. 
Let $p^{\prime} : Y^{\prime} \to Y$ be the shortest 
composite of quadric transformations such that 
the morphism $\varPhi_{K_Y} \circ p^{\prime}$ factors 
through $\varPhi_{D^{\prime}} : Z^{\prime} \to Z$.  
We denote by $f^{\prime} : Y^{\prime} \to Z^{\prime}$ 
the unique morphism satisfying 
$\varPhi_{K_Y} \circ p^{\prime} = \varPhi_{D^{\prime}} \circ f^{\prime}$. 
Then we have $K_{Y^{\prime}} \sim {p^{\prime}}^* K_Y + \eta$ for 
a certain effective divisor $\eta$ on $Y^{\prime}$.   
%The morphism $f^{\prime}$ contracts no $(-1)$-curve, 
%since $Y$ is of general type. 
%To prove the liftability, we only need to show 
%$f^{\prime}_* \eta =0$, since, in this case, 
%$p^{\prime} : Y^{\prime} \to Y$ is an isomorphism. 
If $f^{\prime}_* \eta =0$, 
then $p^{\prime} : Y^{\prime} \to Y$ is an isomorphism
%, since $p^{\prime}$ is the shortest one
. 
Thus we only need to show $f^{\prime}_* \eta =0$. 

So we prove the equality above. 
Let $R^{\prime}$ be the ramification divisor of $f^{\prime}$, 
and $B^{\prime} = f^{\prime}_* R^{\prime}$, 
its direct image. 
Then from 
$ R^{\prime} \sim K_{Y^{\prime}} - {f^{\prime}}^* K_{Z^{\prime}}
%             \sim {p^{\prime}}^* K_Y + \eta - f^{\prime}_* K_{Z_{\prime}}
             \sim {f^{\prime}}^* (2D^{\prime} - r^*\varGamma) + \eta $,
we infer  
\begin{equation} \label{eq:equivBprime}
 B^{\prime} \sim 2(2D^{\prime} -r^*\varGamma + \alpha), 
\end{equation}
where $\alpha$ is a divisor  
satisfying $2\alpha \sim f^{\prime}_* \eta$.  
We denote by $Y^{\prime \prime} \to Z^{\prime}$  
the double cover branched along $B^{\prime}$,  
and by $Y^{\sharp} \to Y^{\prime \prime}$ its canonical resolution. 
To show the equality $f^{\prime}_* \eta =0$,
we compute the Euler characteristic 
$\chi (\mathcal{O}_{Y^{\sharp}})$ in two ways and compare them. 
Note that $\dim (\varPhi_{K_Y} \circ p^{\prime}) (\eta) =0$, 
and that any general member of $|r^*\varGamma|$ is a 
$0$-curve. It follows
%\[
% \alpha D^{\prime} = \frac{D^{\prime} f^{\prime}_* \eta }{2} =0 
% \quad \textrm{and} \quad 
% D^{\prime} (r^* \varGamma) = -K_{Z^{\prime}} (r^* \varGamma) =2.  
%\]
$D^{\prime} \alpha = {D^{\prime} f^{\prime}_* \eta }/2 =0$ and 
$D^{\prime} (r^* \varGamma) = -K_{Z^{\prime}} (r^* \varGamma) =2$. 
Thus by (\ref{eq:equivBprime}) and 
\cite[Lemma 6]{quintic}, we obtain   
\begin{align}
         \chi (\mathcal{O}_{Y_{\sharp}}) 
              &= 2 + \frac{1}{2} (2D^{\prime} -r^*\varGamma + \alpha)
                 ((2D^{\prime} -r^*\varGamma + \alpha) + K_{Z^{\prime}}) 
                  - \beta                              \notag    \\
              &= 2 + \frac{1}{2} (2D^{\prime} -r^*\varGamma + \alpha)   
                  (D^{\prime} + \alpha )  
                  - \beta                             \notag     \\
              &= {D^{\prime}}^2 + 1 
                     - \frac{1}{4} (r^*\varGamma) (f^{\prime}_* \eta ) 
                     + \frac{1}{8}(f^{\prime}_* \eta)^2  
                     - \beta ,
                   \label{eq:chiOYsharp}  
\end{align}
where $\beta$ is a term coming from  
essential singularities of the branch divisor $B^{\prime}$. 
Here, we have three inequalities 
\begin{equation}
     - \frac{1}{4} (r^*\varGamma) (f^{\prime}_* \eta ) \leq 0 , \quad     
     \frac{1}{8}(f^{\prime}_* \eta)^2  \leq 0 ,  \quad 
     \text{and}   \quad 
     - \beta \leq 0.        \label{eq:twoinequalities}
\end{equation} 
The first one follows from the absence of base points of $|r^*\varGamma|$, 
the second one from 
${D^{\prime}}^2 > 0$ and $D^{\prime} f^{\prime}_* \eta = 0$,
and the last one from the definition of $\beta$.  
Meanwhile we have 
$\chi (\mathcal{O}_{Y^{\sharp}}) =
 \chi (\mathcal{O}_Y) = n +2 = {D^{\prime}}^2 +1$. 
Thus by (\ref{eq:chiOYsharp}) and (\ref{eq:twoinequalities}),  
we obtain $(f^{\prime}_* \eta)^2 =0$, from which together with 
Hodge's index theorem, we infer $f^{\prime}_* \eta =0$. 
Hence the canonical map $\varPhi_{K_Y}$ lifts. 

The remaining assertion easily follows from the proof above.  \qed 

Note that the  action by $G = \mathrm{Gal} (Y/X)$ 
on $Y$ induces one on $Z^{\prime}$. 
We can verify it as follows.   
Since $Z$ is the canonical image of our surface $Y$, 
the action on $Y$ induces one on $Z$. 
Meanwhile the surface $Z^{\prime}$ is the minimal 
desingularization of our surface $Z$. 
Thus this action on $Z$ induces one on $Z^{\prime}$. 
   
\begin{lemma} \label{lm:onedimcomp}
The induced action by $G$ on $Z^{\prime}$ is non-trivial. 
The fixed locus of this action has a one-dimensional 
irreducible component $C_0^{\prime}$ 
satisfying ${C_0^{\prime}}^2 \equiv 1 $ $\mathrm{mod}$ $2$.  
\end{lemma}

Proof. 
The first assertion is trivial, since the action on 
$Y$ has no fixed point. 
Let us show the second assertion. 
Let $\{ z_1, \ldots , z_b \}$ be the set of 
isolated fixed points of the action on $Z^{\prime}$, 
and $r^{\prime \prime} : Z^{\prime \prime} \to Z^{\prime}$ 
the blowing-up at these $b$ points. 
We denote by $C_i^{\prime \prime}$ the 
$(-1)$-curve lying over $z_i$.
Let $\{ C_1^{\prime}, \ldots , C_a^{\prime} \}$ 
be the set of $1$-dimensional irreducible components  
of the fixed locus of the 
action on $Z^{\prime}$.
We use the same symbol $C_i^{\prime}$ for the total transform 
to $Z^{\prime \prime}$ of the divisor  $C_i^{\prime}$.   
Note that the divisor $\sum_{i=1}^{a} C_i^{\prime} 
+ \sum_{i=1}^{b} C_i^{\prime \prime}$ has no singularity,  
since we have $G \simeq \mathbb{Z}/2$.  
It follows that the quotient $Z^{\prime \prime} /G$ is smooth, 
where the action by $G$ is the lifting of that on $Z^{\prime}$. 
We denote by $\bar{C_i^{\prime}}$ and $\bar{C_i^{\prime \prime}}$  
the image to $Z^{\prime \prime} /G$ of 
the divisor $C_i^{\prime}$ 
and that of the divisor $C_i^{\prime \prime}$, respectively. 
Then since the branch divisor 
$\sum_{i=1}^{a} \bar{C_i^{\prime}} 
+ \sum_{i=1}^{b} \bar{C_i^{\prime \prime}}$ 
is linearly equivalent to twice a divisor on $Z^{\prime \prime} /G$,   
we have  
\[ 
 \sum_{i=1}^{a}{C_i^{\prime}}^2 - b 
= 
  (\sum_{i=1}^{a} C_i^{\prime} 
+ \sum_{i=1}^{b} C_i^{\prime \prime})^2  
=
  (\sum_{i=1}^{a} \bar{C_i^{\prime}} 
+ \sum_{i=1}^{b} \bar{C_i^{\prime \prime}})^2 /2 
   \equiv 0 \quad \mathrm{mod} \quad 2.
\]
Meanwhile, since 
$K_{Z^{\prime \prime}}$ is linearly equivalent to a pull-back of 
a divisor on $Z^{\prime \prime} /G$, we have 
$K_{Z^{\prime \prime}}^2 = K_{Z^{\prime}}^2 -b =n-3-b \equiv 0 $ 
$\mathrm{mod}$ $2$, hence $b \equiv 1$ $\mathrm{mod}$ $2$. 
Thus we infer 
$\sum_{i=1}^{a}{C_i^{\prime}}^2 \equiv 1$ $\mathrm{mod}$ $2$, 
which implies the second assertion. \qed 

\begin{lemma} \label{lm:bcneqzero}
Let $C_0^{\prime}$ be an irreducible curve as in Lemma 
$\ref{lm:onedimcomp}$. Then $B^{\prime}C_0^{\prime} \neq 0$ holds. 
\end{lemma}

Proof. 
We derive a contradiction by assuming $B^{\prime}C_0^{\prime} =0$. 
Assume that $B^{\prime}C_0^{\prime} =0$ holds. 
Then by Proposition \ref{prop:branchdivisor},  
we have 
\[ 
 (2D^{\prime} - r^* \varGamma ) C_0^{\prime} =
 (-2K_{Z^{\prime}} + r^* \varGamma ) C_0^{\prime} =0.
\]
If $(r^* \varGamma) C_0^{\prime} =0$, 
then by the equality above, we obtain  
$K_{Z^{\prime}} C_0^{\prime} =0$, 
which contradicts  
${C_0^{\prime}}^2 \equiv 1$ $\mathrm{mod}$ $2$. 
Thus we have $(r^* \varGamma) C_0^{\prime} > 0$, 
hence $-2K_{Z^{\prime}} C_0^{\prime} 
= -(r^* \varGamma) C_0^{\prime} < 0$. 
It follows $C_0^{\prime}$ is a fixed component of 
the anti-canonical system $|-K_{Z^{\prime}}|$. 
Then since $-K_{\varSigma_d} \sim 2\varDelta_0 + (2+ d) \varGamma$, 
we obtain $(r^* \varGamma) C_0^{\prime} \leq 2$, 
hence $(r^* \varGamma) C_0^{\prime} 
= 2K_{Z^{\prime}} C_0^{\prime} =2$. 
Thus $r_* C_0^{\prime} \sim 2\varDelta_0 + c \varGamma$ 
holds for a certain integer $c \geq 1$. 
Meanwhile since $0 \leq d \leq 3$, we have 
$h^0(\mathcal{O}_{Z^{\prime}}(-K_{Z^{\prime}}))
\geq h^0(\mathcal{O}_{\varSigma_d}(-K_{\varSigma_d})) -(11 -n)
= n-2$. 
Thus we obtain 
\begin{multline}
n-2 \leq h^0(\mathcal{O}_{Z^{\prime}}(-K_{Z^{\prime}})) 
  = h^0(\mathcal{O}_{Z^{\prime}}(-K_{Z^{\prime}} - C_0^{\prime})) \\  
  \leq h^0(\mathcal{O}_{\varSigma_d}
     (-K_{\varSigma_d} - r_* C_0^{\prime})) 
  = 3 + d -c, \notag 
\end{multline}
hence $ c -d \leq 5 - n <0$. 
It follows $ (r_* C_0^{\prime}) \varDelta_0 = (c -d) -d <0$, 
which contradicts the irreducibility of $C_0^{\prime}$. 
Hence $B^{\prime} C_0^{\prime} \neq 0$ holds.  \qed

Now let us exclude case $1$) in Proposition \ref{prop:L^2}. 

\begin{proposition} \label{prop:case1exclusion}
Case $1$%
$)$ in Proposition $\ref{prop:L^2}$ does not occur. 
\end{proposition}

Take an irreducible curve $C_0^{\prime}$ as in 
Lemma \ref{lm:onedimcomp}. Then by Lemma \ref{lm:bcneqzero}, 
we have $B^{\prime} \cap C_0^{\prime} \neq \emptyset$.
So let us take a point $x \in B^{\prime} \cap C_0^{\prime}$. 
Then the preimage ${f^{\prime}}^{-1} (x) \subset Y$ is stable 
under the action by $G$ on $Y$. 
By Proposition \ref{prop:branchdivisor}, however,
the set ${f^{\prime}}^{-1} (x)$ is either a point 
or a base space of the fundamental cycle of a rational double point. 
This implies that the action by $G$ on ${f^{\prime}}^{-1} (x)$ has 
a fixed point, which contradicts the definition of $\pi : Y \to X$. 
Thus we have the assertion. \qed
  
\subsection{The case $L^2 = K_Y^2 -4$}
Next we exclude cases $3$-$1$) and $3$-$2$) in 
Proposition \ref{prop:L^2}. 
In these two cases, we have $L^2 = 2(n-1)$;   
hence the canonical image $Z$ is a non-degenerate surface 
in $\mathbb{P}^n$ of minimal degree $n-1$.  
Thus from the well-known classification,   
it follows that our $Z$ is a image of the 
Hirzebruch surface $Z^{\prime} = \varSigma_d$ by the morphism associated 
with the complete linear system 
$|D^{\prime}| = |\varDelta_0 + \frac{n-1+d}{2}\varGamma|$, 
where $0 \leq d \leq n-1$ and $d \equiv n-1$ $\mathrm{mod}$ $2$ 
(see \cite{ratsfI}or \cite[Lemma 1.2]{smallc1-1}).  
Let us denote this morphism by 
$\varPhi_{D^{\prime}}: Z^{\prime} \to Z \subset \mathbb{P}^n$. 
Then $\varPhi_{D^{\prime}}$ is an embedding if $d < n-1$, 
and is the contraction of $\varDelta_0$ if $d = n-1$. 
%Hence the morphism $\varPhi_{D^{\prime}}$ gives the minimal 
%desingularization of $Z$. 
Note that in the later case, our $Z$ is a cone over a rational curve 
embedded in $\mathbb{P}^{n-1}$ by $\mathcal{O}_{\mathbb{P}^1} (n-1)$. 

For the case $d < n-1$, the lemma below is trivial. 
For the case $d = n-1$, we can give a proof by the same 
method as in \cite[Lemma 1.5]{smallc1-1}. 

\begin{lemma}
The morphism $\varPhi_L : \tilde{Y} \to Z$
lifts to a morphism $f^{\prime} : \tilde{Y} \to Z^{\prime}$. 
\end{lemma}

By the same argument as in the exclusion of case $1$), 
we see that the action by $G$ on $Y$ induces one on $Z^{\prime}$.

Let us recall the morphism $p: \tilde{Y} \to Y$ and the 
base locus of $|K_Y|$. In case $3$-$1$) in Proposition \ref{prop:L^2}, 
the morphism $p$ is the blowing-up at (possibly infinitely near) 
four points, which we shall call $y_1, \ldots , y_4$. 
Let $E_i$ denote the total transform to $\tilde{Y}$
of the $(-1)$-curve corresponding to $y_i$.  
Then we have $E = \sum_{i=1}^4 E_i$ 
and $L E_i =1$ ($1 \leq i \leq 4$). 
Since the action by $G$ on the set of base points of $|M|$ has 
no fixed point, we have only two cases: 
i) the case where $y_1, \ldots , y_4$ are four distinct 
points on $\tilde{Y}$, and ii) the case where $y_1$ and $y_2$ are 
distinct points on $\tilde{Y}$, 
and $y_{i +2}$ is infinitely near to $y_i$ for $i =1$, $2$.   
In the later case, the divisor $E_i^{\prime} = E_i - E_{i +2}$ is 
a $(-2)$-curve satisfying $L E_i^{\prime} =0$.

Meanwhile in case $3$-$2$), the morphism $p : \tilde{Y} \to Y$ 
is an isomorphism. Hence we may assume $\tilde{Y} = Y$. 
We have $|M|=|L|$ and $F = \sum_{i=1,2} F_i$, 
where $F_i$ is a fundamental cycle of a rational double point. 
Since the action on $Y$ has no fixed point, 
we have $F_1 \cap F_2 = \emptyset$; hence the generator of $G$ 
maps $F_1$ onto $F_2$. 
It follows $L F_1 = L F_2 =2$. 

In what follows, we put $T = 2E$ for case $3$-$1$), 
and $T =F$ for case $3$-$2$). Then we have 
\[
 K_{\tilde{Y}} \sim L +T.
\]   

\begin{lemma} \label{lm:Tmod4}
Let $T$ be the divisor above. 
Then $\varGamma (f^{\prime}_* T) \equiv 2$ $\mathrm{mod}$ $4$ holds. 
\end{lemma}

Proof. 
Since $d \equiv n-1 \equiv 1$ $\mathrm{mod}$ $2$, we have 
$d \neq 0$. Thus the action by $G$ on $Z^{\prime} =\varSigma_d$ 
induces one on $\mathbb{P}^1$ 
via the natural fibration $\varSigma_d \to \mathbb{P}^1$ 
of the Hirzebruch surface. 
It follows there exists a member $\varGamma_0 \in |\varGamma|$ 
stable under the action by $G$. 
Let us take a blowing-up $\tilde{X} \to X$ 
such that $\tilde{Y} = \tilde{X} \times_X Y$ holds. 
The base change $\tilde{\pi} : \tilde{Y} \to \tilde{X}$ 
is an unramified double cover satisfying 
$\mathrm{Gal} (\tilde{Y}/\tilde{X}) \simeq \mathrm{Gal} (Y/X)$. 
Then since ${f^{\prime}}^* \varGamma_0$ is stable under 
the action by $G$ on $\tilde{Y}$,   
the divisor ${f^{\prime}}^* \varGamma_0$ is a pull-back by 
$\tilde{\pi}$ of a certain divisor on $\tilde{X}$.   
Thus from $\tilde{\pi}^* K_{\tilde{X}} \sim K_{\tilde{Y}}$ and 
the Riemann--Roch theorem, we infer 
\[
  ({f^{\prime}}^* \varGamma_0)^2 
            + ({f^{\prime}}^* \varGamma_0) K_{\tilde{Y}}
   = ({f^{\prime}}^* \varGamma_0) (L + T)
   = 2 + ({f^{\prime}}^* \varGamma_0) T
   \equiv 0 \quad \mathrm{mod} \quad 4. 
\]
Hence we have the assertion. \qed

\begin{lemma} \label{lm:-1and-2curves}
The morphism $f^{\prime} : \tilde{Y} \to Z^{\prime}$ contracts 
no $(-1)$-curve on $\tilde{Y}$. 
Further, the following hold$:$

$\mathrm{i}$%
$)$ if $C$ is a $(-1)$-curve on $\tilde{Y}$ satisfying $LC=1$, 
       then $f^{\prime}_* C \sim \varGamma$%
$;$ 

$\mathrm{ii}$%
$)$ if $C$ is a $(-2)$-curve on $\tilde{Y}$ satisfying $LC=0$, 
       then $f^{\prime}$ contracts $C$. 
\end{lemma}

Proof. 
The first assertion trivially follows from 
the definition of $p : \tilde{Y} \to Y$. 
In order to prove i) and ii),  
we put $f^{\prime}_* C \sim a \varDelta_0 + b \varGamma$. 
We denote by $\theta$ the involution of $\tilde{Y}$ over 
$Z^{\prime}$. This involution exists, since $f^{\prime}$ contracts 
no $(-1)$-curve. 

First, let us prove the assertion i). 
Assume that $C$ is a $(-1)$-curve on $\tilde{Y}$ satisfying $LC=1$. 
Then since   
$L \sim {f^{\prime}}^* D^{\prime}$, 
we have  
\begin{equation} \label{eq:fprimec}
  (\varDelta_0 + d \varGamma) f^{\prime}_* C
         + \frac{n-1-d}{2} \varGamma f^{\prime}_* C =1, 
\end{equation}
where each term of the left hand is a non-negative integer. 
Thus we obtain $(\varDelta_0 + d \varGamma) f^{\prime}_* C =0$ or $1$. 
Assume that $(\varDelta_0 + d \varGamma) f^{\prime}_* C =1$. 
Then we have $f^{\prime}_* C \sim a \varDelta_0 + \varGamma$ 
and $(\frac{n-1-d}{2})a =0$. Thus, in this case,  
we only have to show $a=0$, 
which is trivial if $n-1-d \neq 0$. 
If $n-1-d=0$, then by the irreducibility of $C$, we have 
$\varDelta_0 f^{\prime}_* C = 1 -a (n-1) \geq 0$, hence $a=0$.     
Assume next that $(\varDelta_0 + d \varGamma) f^{\prime}_* C =0$.
Then by (\ref{eq:fprimec}), we obtain 
$f^{\prime}_* C = \varDelta_0 $ and $d = n-3$.  
We exclude this case as follows.    
We have ${f^{\prime}}^* \varDelta_0 = C + \theta (C) + \xi$ for a certain 
effective divisor $\xi$ exceptional with respect to $f^{\prime}$. 
It follows 
\[
({f^{\prime}}^* \varDelta_0)^2 
= (C + \theta (C) + \xi)(C + \theta (C)) 
\geq C^2 + {\theta (C)}^2 + 2 C \theta (C)
\geq -4,
\]    
hence $-2(n-3) \geq -4$.  
This contradicts $\lambda \geq 4$.
Thus we have $(\varDelta_0 + d \varGamma) f^{\prime}_* C \neq 0$, 
which completes the proof of the assertion i). 

Next, let us prove the assertion ii). 
Assume that $f^{\prime}(C)$ is a curve. Then since 
$\varPhi_{D^{\prime}}$ contracts $f^{\prime}(C)$, 
we have $d = n-1$ and $f^{\prime}(C) = \varDelta_0$. 
Note that we have $f^{\prime}_*C = \varDelta_0$
or $2\varDelta_0$, since $\deg f^{\prime} =2$. 
Assume that $f^{\prime}_*C = \varDelta_0$. Then we have 
${f^{\prime}}^* \varDelta_0 = C + \theta (C) + \xi$ for a certain 
effective divisor $\xi$ exceptional with respect to $f^{\prime}$. 
Then by the same method as in the proof of i), we obtain 
$-2(n-1) = ({f^{\prime}}^* \varDelta_0)^2 \geq -8$, 
which contradicts $\lambda \geq 4$. 
%Thus we have $f^{\prime}_*C \neq \varDelta_0$. 
Assume next that $f^{\prime}_*C = 2\varDelta_0$. 
Then we have ${f^{\prime}}^* \varDelta_0 = C + \xi$ for a certain 
effective divisor $\xi$ exceptional with respect to $f^{\prime}$. 
Then again by the same method, we obtain 
$-2(n-1) \geq -2$, which contradicts $\lambda \geq 4$.  
Thus we have the assertion ii). \qed

If our $Y$ is of case $3$-$1$) in Proposition \ref{prop:L^2}, 
then by the lemma above we have 
$f^{\prime}_* T = 2f^{\prime}_* E \sim 8 \varGamma$, 
which contradicts Lemma \ref{lm:Tmod4}. 
Thus we have the following: 
\begin{proposition}  \label{prop:case3-1exclusion}
Case $3$-$1$%
$)$ in Proposition $\ref{prop:L^2}$ does not occur. 
\end{proposition}

So in what follows, we assume that our $Y$ is of case 
$3$-$2$) in Proposition \ref{prop:L^2}.  

\begin{lemma} \label{lm:-2curves}
Let $C$ be an irreducible component of $F_1$ satisfying 
$D^{\prime} f^{\prime}_* C > 0$. Then one of the following holds$:$

$\mathrm{i}$%
$)$ $D^{\prime} f^{\prime}_* C =1$ and $f^{\prime}_* C \sim \varGamma$%
$;$  

$\mathrm{ii}$%
$)$ $D^{\prime} f^{\prime}_* C =2$ and $f^{\prime}_* C \sim 2 \varGamma$%
$;$  

$\mathrm{iii}$%
$)$ $D^{\prime} f^{\prime}_* C =2$, $f^{\prime}_* C = \varDelta_0$, 
        and $d = n - 5 =1$. 
\end{lemma}

Proof. 
First, note that if $f^{\prime}(C) = \varDelta_0$, then 
we have $C \neq \theta (C)$, where $\theta$ is the involution 
of $\tilde{Y}=Y$ over $Z^{\prime}$. 
We can verify this as follows. Let $\iota$ be 
the generator of the Galois group $G$, and 
$\iota |_{Z^{\prime}}$, the corresponding automorphism of $Z^{\prime}$. 
Then since $d \neq 0$, we have 
$f^{\prime} (\iota (C)) 
= \iota |_{Z^{\prime}} (f^{\prime}(C)) 
= \varDelta_0 = f^{\prime} (C)$. 
This means $C \neq \theta(C) = \iota (C)$, 
since we have 
$\iota (C) \subset F_2$ and $F_1 \cap F_2 = \emptyset$. 
Next, note that $C$ is a $(-2)$ curve satisfying 
$0  < D^{\prime} f^{\prime}_* C \leq D^{\prime} f^{\prime}_* F_1 = 2$. 
Then we can prove the assertion by the same method as in the 
proof of Lemma \ref{lm:-1and-2curves}.  \qed

By $D^{\prime} f^{\prime}_* F_1 =2$ together with   
Lemmas \ref{lm:-1and-2curves} and \ref{lm:-2curves}, 
we see that either of the following holds:  

a) $f^{\prime}_* F_1 = f^{\prime}_* (\iota (F_2)) \sim 2\varGamma$;     

b) $f^{\prime}_* F_1 = f^{\prime}_* (\iota (F_2)) = \varDelta_0$, 
    and $d = n-5 =1$,

\noindent
where $\iota$ is the generator of the Galois group of $G$. 
Case a) above, however, contradicts the assertion in 
Lemma \ref{lm:Tmod4}. Thus we have the following:  
 
\begin{lemma} \label{lm:F1F2image}
$f^{\prime}_* F_1 = f^{\prime}_* F_2 = \varDelta_0$
and $d =n-5 =1$. 
\end{lemma}

Now let us study the morphism 
$f^{\prime} : \tilde{Y} =Y \to Z^{\prime} = \varSigma_1$. 
Let $R^{\prime}$ be the ramification divisor of $f^{\prime}$, 
and $B^{\prime} = f^{\prime}_* R^{\prime}$, the branch divisor. 
Then by the lemma above we obtain 
\begin{equation} \label{eq:ramificationdivisor}
R^{\prime} %&\sim K_{\tilde{Y}} - {f^{\prime}}^*K_{Z^{\prime}}
           \sim {f^{\prime}}^*(3\varDelta_0 + 6\varGamma)
                 + \sum_{i =1, 2} F_i    \qquad \text{and} \qquad 
B^{\prime} \sim 2(4\varDelta_0 + 6\varGamma) . 
           %\sim {f^{\prime}}^* (-4K_{Z^{\prime}})
\end{equation}
We take the double cover of $Z^{\prime}$ with branch divisor $B^{\prime}$, 
and denote by $Y^{\sharp}$ its canonical resolution. 
Let us recall how to obtain the canonical resolution. 
Set $Z^{\prime}_0 = Z^{\prime}$ and $B^{\prime}_0 = B^{\prime}$.  
We define $Z^{\prime}_i$ and $B^{\prime}_i$ inductively as follows.  
Choose a singularity $z_i$, if any, of $B^{\prime}_{i-1}$, 
and take the blowing-up  
$q^{\prime}_i : Z^{\prime}_i \to Z^{\prime}_{i-1}$ at this point. 
We denote by $\varepsilon_i$ the $(-1)$-curve corresponding to $z_i$. 
Let $m_i$ be the multiplicity of $B^{\prime}_{i-1}$ at $z_i$, and 
$[\frac{m_i}{2}]$, the largest integer not exceeding $\frac{m_i}{2}$.  
Then we define $B^{\prime}_i$ by 
$B^{\prime}_i 
= {q^{\prime}_i}^* B^{\prime}_{i-1} - 2[\frac{m_i}{2}] \varepsilon_i$. 
For a certain $s \geq 0$,
the divisor $B^{\prime}_s$ is non-singular. 
So take the double cover 
$f^{\sharp} : \tilde{Y_s} \to Z^{\sharp} = Z^{\prime}_s$ 
with branch divisor $B^{\sharp} = B^{\prime}_s$.  
Then this $\tilde{Y_s}$ is our canonical resolution $Y^{\sharp}$. 
Put $q^{\prime} = 
(q^{\prime}_1 \circ q^{\prime}_2 \circ \cdots \circ q^{\prime}_s): 
Z^{\sharp} \to Z^{\prime}$. 
There exists a natural birational  
morphism $p^{\sharp} : Y^{\sharp} \to \tilde{Y}$ 
satisfying $q^{\prime} \circ f^{\sharp} = f^{\prime} \circ p^{\sharp}$. 
We use the same symbol $\varepsilon_i$ for the total transform to 
$Z^{\sharp}$ of the $(-1)$-curve $\varepsilon_i \subset Z^{\prime}_i$. 
Note, for our case, the action by 
the Galois group $G = \mathrm{Gal} (Y/X)$ on $\tilde{Y}$ induces 
one on $Z^{\sharp}$ and one on $Y^{\sharp}$. 
This action on $Y^{\sharp}$ is free. 

By the same method as in \cite[Section 2]{quintic}, 
we obtain the following: 

\begin{proposition} \label{prop:i1i2}
There exist $i_1$ and $i_2$ 
$(i_1 < i_2)$  
satisfying $[\frac{m_{i_1}}{2}] = [\frac{m_{i_2}}{2}] = 2$. 
For any $i \neq i_1$, $i_2$,  
the equality $[\frac{m_i}{2}] = 1$ holds. 
The morphism $p^{\sharp} : Y^{\sharp} \to \tilde{Y} = Y$ is 
a composite of two quadric transformations. 
\end{proposition}

Thus the branch divisor $B^{\prime}$ has an 
essential singularity. By the proposition above, we obtain 
\begin{equation} \label{eq:KYsharp}
 K_{Y^{\sharp}} \sim {f^{\sharp}}^* 
    ({q^{\prime}}^*(2\varDelta_0 + 3\varGamma) 
       - \varepsilon_{i_1} - \varepsilon_{i_2}) . 
\end{equation}

\begin{lemma} \label{lm:bprimeessentialsing}
Every essential singularity of 
$B^{\prime}$ lies on $\varDelta_0$. 
\end{lemma}

Proof. 
Since $f^{\prime}$ contracts no $(-1)$-curve, 
${f^{\prime}}^* B^{\prime} - 2 R^{\prime} = 2 \zeta^{\prime}$ 
holds for a certain effective divisor $\zeta^{\prime}$ on $\tilde{Y}$. 
This $\zeta^{\prime}$ satisfies  
\begin{equation} \label{eq:zetaprime}
 2 \zeta^{\prime} \sim 
2 ({f^{\prime}}^*(\varDelta_0) - \sum_{i =1, 2} F_i),  
\end{equation}
since we have (\ref{eq:ramificationdivisor}). 
Let $\zeta^{\prime} = \sum \zeta^{\prime}_i$ be the decomposition into 
connected components. 
Note that $f^{\prime}$ maps each $\zeta^{\prime}_i$ to a point 
on $Z^{\prime}$.  
Then, for any $i$ satisfying 
$f^{\prime} (\zeta^{\prime}_i) \notin \varDelta_0$, 
we infer from (\ref{eq:zetaprime}) that
${\zeta^{\prime}_i}^2 = \zeta^{\prime}_i \zeta^{\prime} =0$, 
hence $\zeta^{\prime}_i =0$, which implies the assertion. \qed  

\begin{lemma} \label{lm:fixedpartKYsharp} 
Let $\eta^{\sharp} \sim K_{Y^{\sharp}} - {p^{\sharp}}^* K_{\tilde{Y}}$ 
be the exceptional divisor corresponding to 
$p^{\sharp} : Y^{\sharp} \to \tilde{Y}$. 
Then the fixed part of $|K_{Y^{\sharp}}|$ is given by 
$\sum_{i=1,2} {p^{\sharp}}^* F_i + \eta^{\sharp}$. 
Further, the linear equivalence 
$\sum_{i=1,2} {p^{\sharp}}^* F_i + \eta^{\sharp}
\sim 
{f^{\sharp}}^*({q^{\prime}}^*\varDelta_0 
-\varepsilon_{i_1} -\varepsilon_{i_2})$
holds, where 
$i_1$ and $i_2$ are integers given 
in Proposition $\ref{prop:i1i2}$.  
\end{lemma}

Proof. 
The first assertion follows from 
$|K_{Y^{\sharp}}|= |K_{\tilde{Y}}| + \eta^{\sharp}$, 
since $|L|$ has no base point. 
The second assertion follows 
from (\ref{eq:KYsharp}) and 
$\sum {p^{\sharp}}^* F_i + \eta^{\sharp}
\sim K_{Y^{\sharp}} - {p^{\sharp}}^* L
\sim K_{Y^{\sharp}} - {p^{\sharp}}^* {f^{\prime}}^* D^{\prime}$. 
\qed

\begin{lemma} \label{lm:gamma1}
There exists a member $\varGamma_1 \in |\varGamma|$ 
contained in the 
fixed locus of the action by $G$ on 
$Z^{\prime} = \varSigma_1$. 
\end{lemma}

Proof. 
The action by $G$ on $Z^{\prime} = \varSigma_1$ 
induces one on $\mathbb{P}^1$ via the natural 
fibration $Z^{\prime} = \varSigma_1 \to \mathbb{P}^1$ 
of the Hirzebruch surface. 
Let us show that this induced action on $\mathbb{P}^1$ is 
non-trivial. 
There exists a member $\varDelta_1 \in |\varDelta_0 + \varGamma|$ 
stable under the action by G satisfying 
$\varDelta_1 \cap \varDelta_0 = \emptyset$. 
Assume that the induced action on $\mathbb{P}^1$ is trivial.    
Then this $\varDelta_1$ is contained in the fixed locus 
of the action by $G$ on $Z^{\prime}$. 
From this together with $B^{\prime} \varDelta_1 = 12$ and 
Lemma \ref{lm:bprimeessentialsing}, 
it follows that $B^{\prime}$ has   
a smooth point or a negligible singularity   
that is stable under the action by $G$. 
This, however, leads us to a contradiction by the same 
argument as in the proof of Proposition \ref{prop:case1exclusion}. 
Thus the induced action on $\mathbb{P}^1$ is non-trivial. 
Now take two fibers %, say $\varGamma_1$ and $\varGamma_2$, 
of $Z^{\prime} \to \mathbb{P}^1$ that  
lie over the fixed points of the action on $\mathbb{P}^1$. 
Since $Z^{\prime} = \varSigma_1$, one of these two fibers 
are contained in the fixed locus of the action by $G$. \qed

Let us exclude case $3$-$2$) in Proposition \ref{prop:L^2}. 
\begin{proposition}   \label{prop:case3-2exclusion}
Case $3$-$2$%
$)$ in Proposition $\ref{prop:L^2}$ does not occur. 
\end{proposition}

Proof. 
Let $\varGamma_1 \in |\varGamma|$ be the member as in 
Lemma \ref{lm:gamma1}. 
By (\ref{eq:ramificationdivisor}), we have 
$B^{\prime} \varGamma_1 =8$, hence  
$B^{\prime} \cap \varGamma_1 \neq \emptyset$. 
If a smooth point or a negligible singularity of $B^{\prime}$ 
lies on $B^{\prime} \cap \varGamma_1$, we can derive a contradiction 
by the same argument as in the proof of 
Proposition \ref{prop:case1exclusion}. 
Thus by Lemma \ref{lm:bprimeessentialsing}, we see that  
$B^{\prime} \cap \varGamma_1 = \varDelta_0 \cap \varGamma_1$ 
and that this point is an essential singularity of $B^{\prime}$. 
So we put $\varDelta_0 \cap \varGamma_1 = \{ z_1\}$, where 
the point $z_1$ is the center of the first blowing-up 
$q^{\prime}_1: Z^{\prime}_1 \to Z^{\prime}_0 = Z^{\prime}$ in 
the procedure to obtain the canonical resolution $Y^{\sharp}$. 
Then, by Proposition \ref{prop:i1i2}, 
we have $3 \leq m_1 \leq 5$. 
If $m_1$ is odd, then the strict transform 
$\varepsilon_1^{\sharp} \simeq \mathbb{P}^1 \subset Z^{\sharp}$ of 
the exceptional curve $\varepsilon_1 \subset Z_1^{\prime}$ 
is a component of $B^{\sharp}$ stable under the action by $G$. 
This, however, leads us to a 
contradiction, since the action by $G$ on $Y^{\sharp}$ is free. 
It follows $m_1 \equiv 0$ $\mathrm{mod}$ $2$, hence $m_1 =4$. 
Thus we have 
$B_1^{\prime} = {q_1^{\prime}}^* B^{\prime} - 4 \varepsilon_1$
and $B_1^{\prime}  {q_1^{\prime}}^{-1}_* (\varGamma_1) =4$, 
where the divisor ${q_1^{\prime}}^{-1}_* (\varGamma_1)$ is the 
strict transform of $\varGamma_1$ by 
$q_1^{\prime} : Z_1^{\prime} \to Z^{\prime}$. 
Note that the action by $G$ on $Z^{\prime}$ induces one on $Z_1^{\prime}$, 
and that the strict transform ${q_1^{\prime}}^{-1}_* (\varGamma_1)$ is 
contained in the fixed locus of this induced action. 
By the same argument as that on $\varGamma_1$ above, we see that 
the point $B_1^{\prime} \cap {q_1^{\prime}}^{-1}_* (\varGamma_1) 
= \varepsilon_1 \cap {q_1^{\prime}}^{-1}_* (\varGamma_1)$ 
is an essential singularity of $B_1^{\prime}$,  
that we can set $\varepsilon_1 \cap {q_1^{\prime}}^{-1}_* (\varGamma_1) = 
\{ z_2 \}$, where the point $z_2$ is the center of the second 
blowing-up $q_2^{\prime} : Z_2^{\prime} \to Z_1^{\prime}$, 
and that $m_2= 4$, where $m_2$ is the multiplicity of $B_1^{\prime}$ 
at $z_2$. 
Thus we have $i_1 =1$ and $i_2 =2$, where $i_1$ and $i_2$ are the 
integers given in Proposition \ref{prop:i1i2}. 

Now we derive a contradiction. 
Let $\varGamma_1^{\sharp}$ be the strict transform to $Z^{\sharp}$ of 
of the divisor $\varGamma_1$. 
Note that we have $z_1 \in \varGamma_1$ 
and $z_2 \in {q_1^{\prime}}^{-1}_* (\varGamma_1)$.   
Thus by Lemma \ref{lm:fixedpartKYsharp}, we obtain 
\[
 f^{\sharp}_* (\sum {p^{\sharp}}^* F_i + \eta^{\sharp}) 
      \varGamma_1^{\sharp}
 = 2( \varDelta_0 \varGamma + {\varepsilon_1}^2 + {\varepsilon_2}^2 ) 
 = -2 < 0. 
\] 
From this together with Lemma \ref{lm:F1F2image}, 
we infer that the divisor $\varGamma_1^{\sharp}$ is the  
image by $f^{\sharp}$ of an irreducible component of $\eta^{\sharp}$, 
which contradicts the equality 
$\dim (q^{\prime} \circ f^{\sharp}) (\eta^{\sharp})=
\dim (f^{\prime} \circ p^{\sharp}) (\eta^{\sharp}) = 0$. 
Hence we have the assertion. \qed 

\subsection{The case $L^2 = K_Y^2 -2$} 

Finally, we study case $2$) in Proposition \ref{prop:L^2}. 
It will turn out that $\lambda =4$ in this case, and that 
the surfaces of this case have the structure as in the statement 
of Theorem \ref{thm:completedescription}. 
In what follows, we assume that our $Y$ is of case $2$) in 
Proposition \ref{prop:L^2}, hence $ \deg Z = L^2 /2 = n$. 
Note that in this case, the morphism $p : \tilde{Y} \to Y$ is 
a blowing-up at two distinct points on $Y$. 
Let $E_1$ and $E_2$ denote the $(-1)$-curves corresponding to the
centers of this blowing-up. 
Then we have 
$p^*|K_Y| = |L| + \sum_{i=1, 2} E_i$ 
and $ L E_1 = L E_2 = 1 $. 
The Galois group $G = \mathrm{Gal} (Y/X)$ acts transitively on  
the set $\{ E_1, E_2 \}$.  
We denote by $Z^{\prime}$ the minimal desingularization of $Z$. 

\begin{lemma}  \label{lm:L^2-2resoution}
There exists a blowing-up $r : Z^{\prime} \to \mathbb{P}^2$ at 
$($possibly infinitely near$)$ $9-n$ points such 
that the anticanonical morphism 
$Z^{\prime} \to Z \subset \mathbb{P}^n$ of $Z^{\prime}$ gives 
the minimal desingularization of $Z$. 
\end{lemma}

Proof. 
Note that our $Z = \varPhi_{K_Y} (Y)$ is a non-degenerate surface in 
$\mathbb{P}^n$ of degree $n$. 
Hence our $Z$ is  one of the following 
(see \cite{ratsfI} or \cite[Section 3]{smallc1-4}): 

%\noindent
i) a projection of a surface of degree $n$ in $\mathbb{P}^{n + 1}$ 
from a point outside the surface; 

%\noindent
ii) the Veronese embedding into $\mathbb{P}^8$ of a quadric 
in $\mathbb{P}^3$ ($n=8$); 

%\noindent
iii) the anticanonical image of $\mathbb{P}^2$ 
blown up at $9-n$ points;

%\noindent 
iv) a cone over an elliptic curve in $\mathbb{P}^{n-1}$ of degree $n$. 

Since $Z^{\prime} \to Z$ is given by a complete linear system, 
case i) above is impossible for our case. 
Since $q (Y) = 0$, case iv) also is impossible. 
Thus it suffices to exclude case ii).  
In case ii), however, the divisor $L$ is linearly equivalent 
to twice a divisor on $\tilde{Y}$, which contradicts the equality 
$LE_i =1$. Hence we have the assertion. \qed

In what follows, we put $D^{\prime} = - K_{Z^{\prime}}$ and 
denote by $\varPhi_{D^{\prime}} : Z^{\prime} \to Z \subset \mathbb{P}^n$ 
the anticanonical map of $Z^{\prime}$. 
Note that the action by $G = \mathrm{Gal} (Y/X)$ 
on $\tilde{Y}$ induces one on $Z^{\prime}$.   

\begin{lemma} \label{lm:lifalilitytcriterion}
If the surface $Z^{\prime}$ has no $(-2)$-curve, 
or if every $(-2)$-curve on $Z^{\prime}$ is stable 
under the action by $G$ on $Z^{\prime}$, 
then $\varPhi_L : \tilde{Y} \to Z \subset \mathbb{P}^n$ 
lifts to a morphism $f^{\prime} : \tilde{Y} \to Z^{\prime}$. 
\end{lemma}

Proof. 
Take the shortest composite $p^{\prime} : Y^{\prime} \to \tilde{Y}$ 
of quadric transformations such that  
$Y^{\prime}$ admits 
a morphism $f^{\prime} : Y^{\prime} \to Z^{\prime}$ satisfying 
$\varPhi_L \circ p^{\prime} = \varPhi_{D^{\prime}} \circ f^{\prime}$. 
Then the action by $G$ on $\tilde{Y}$ induces one on $Y^{\prime}$. 
%Indeed, we obtain our $Y$ as follows: first, take the fiber product 
%$\tilde{Y} \times_Z Z^{\prime}$; second, take the reduction of this 
%fiber product; third, take the normalization of this reduction;  
%finally; take the minimal desingularization of this normalization.   
Note that $f^{\prime}$ contracts no $(-1)$-curve. 
This follows from $LE_i =1$ and the definition of $p^{\prime}$, 
since the surface $Y$ is of general type. 
To obtain the assertion, 
we only need to show that $p^{\prime} : Y^{\prime} \to \tilde{Y}$ 
is an isomorphism. 
Assume that $p^{\prime} : Y^{\prime} \to \tilde{Y}$ 
is not an isomorphism. 
Then there exists a $(-1)$-curve $C$ on $Y^{\prime}$ exceptional 
with respect to $p^{\prime}$. 
Since the anticanonical map 
$\varPhi_{D^{\prime}} : Z^{\prime} \to Z \subset \mathbb{P}^n$ 
contracts $f^{\prime} (C)$ to a point,  
the curve $f^{\prime} (C)$ is a $(-2)$-curve on $Z^{\prime}$,   
hence, by the assumption in the statement, 
stable under the action by $G$ on $Z^{\prime}$. 
Meanwhile by the same method as in Lemma \ref{lm:-1and-2curves}, 
we see that $f^{\prime}_* C = f^{\prime} (C)$ or $2 f^{\prime} (C)$, 
and that if $f^{\prime}_* C = f^{\prime} (C)$, then $C$ is a 
component of the ramification divisor of $f^{\prime}$.  
It follows that $C \simeq \mathbb{P}^1$ is stable under 
the action by $G$ on $Y^{\prime}$, which implies the existence of  
fixed points of this action.   
%This, however, contradicts the definition of $\pi : Y \to X$,  
%hence $p^{\prime} : Y^{\prime} \to \tilde{Y}$ is an isomorphism. 
%Thus we have the assertion. \qed
This, however, contradicts the definition of $\pi : Y \to X$. 
Thus we have the assertion.  \qed

\begin{lemma} \label{lm:fprimeE_i}
Assume that $\varPhi_L : \tilde{Y} \to Z$ lifts to a morphism 
$f^{\prime} : \tilde{Y} \to Z^{\prime}$. 
Then $f^{\prime} (E_1)$ and $f^{\prime} (E_2)$ are $(-1)$-curves 
on $Z^{\prime}$.
Further, the following hold$:$

%\noindent
$\mathrm{i}$%
$)$ $f^{\prime}_* E_i = f^{\prime} (E_i)$ for $i =1$, $2$%
$;$

%\noindent 
$\mathrm{ii}$%
$)$ the ramification divisor $R^{\prime}$ of $f^{\prime}$ satisfies 
$R^{\prime} \sim 
{f^{\prime}}^*(-2K_{Z^{\prime}}) + 2 \sum_{i =1, 2} E_i $%
$;$ 

%\noindent
$\mathrm{iii}$%
$)$ the branch divisor $B^{\prime}$ of $f^{\prime}$ satisfies 
$B^{\prime} \sim 
-4K_{Z^{\prime}} + 2 \sum_{i=1, 2} f^{\prime} (E_i)$%
$;$
 
%\noindent
$\mathrm{iv}$%
$)$ $f^{\prime} (E_1)$ and $f^{\prime} (E_2)$ are  
distinct components of the branch divisor $B^{\prime}$. 
\end{lemma}

Proof. 
The first assertion and the assertion i) follow 
from $E_i L = E_i {f^{\prime}}^* D^{\prime} =1$, which 
implies $\varPhi_L (E_i)$ is a line in $\mathbb{P}^n$. 
The assertions ii) and iii) follow from 
$D^{\prime} \sim -K_{Z^{\prime}}$ and the assertion i). 
So it suffices to prove the assertion iv). 
Let us prove the assertion iv). 
Let $\theta$ be the involution of $\tilde{Y}$ over $Z^{\prime}$. 
Since $Y$ is of general type, 
the divisors $E_1$ and $E_2$ are the only $(-1)$-curves on $\tilde{Y}$. 
It follows that if $f^{\prime}(E_1) \neq f^{\prime}(E_2)$, 
then $\theta (E_i) =  E_i$ holds for $i=1$, $2$.  
Thus we only need to show $f^{\prime}(E_1) \neq f^{\prime}(E_2)$. 
Assume that $f^{\prime}(E_1) = f^{\prime}(E_2)$. 
Then 
${f^{\prime}}^*(f^{\prime}(E_1)) = 
 {f^{\prime}}^*(f^{\prime}(E_2)) = E_1 + E_2 + \xi$ holds 
for a certain effective divisor $\xi$ 
exceptional with respect to $f^{\prime}$.  
Since we have $E_1 \cap E_2 = \emptyset$, 
we see, by the same method as in the proof of Lemma \ref{lm:-1and-2curves}, 
that $\xi ^2 = -(E_1 + E_2) \xi = 0$, hence $\xi =0$. 
It follows ${f^{\prime}}^*(f^{\prime}(E_1)) = 
 {f^{\prime}}^*(f^{\prime}(E_2)) = E_1 + E_2$.  
From this together with the assertions ii) and iii), we infer 
${f^{\prime}}^* B^{\prime} - 2R^{\prime} =0$, which implies that 
the branch divisor $B^{\prime}$ has at most negligible singularities. 
Thus by \cite[Lemma 6]{quintic}, we obtain 
%\begin{align*} 
% \chi (\mathcal{O}_{\tilde{Y}}) &= 
%  2\chi (\mathcal{O}_Z^{\prime}) 
%  + \frac{1}{2} (-2K_{Z^\prime} + \sum f^{\prime} (E_i))
%                ( -K_{Z^\prime} + \sum f^{\prime} (E_i)) \\
%  &=
%  2 + \frac{1}{2} (2 K_{Z^{\prime}}^2 
%                     - 6 K_{Z^{\prime}} f^{\prime}(E_1)  
%                     + 4 f^{\prime}(E_1)^2) \\
%  &= n + 3,
%\end{align*}
\[
  \chi (\mathcal{O}_{\tilde{Y}}) = 
  2\chi (\mathcal{O}_Z^{\prime}) 
  + \frac{1}{2} (-2K_{Z^\prime} + \sum f^{\prime} (E_i))
                ( -K_{Z^\prime} + \sum f^{\prime} (E_i)) 
  = n + 3, 
\]
which contradicts $\chi (\mathcal{O}_Y) = n + 2$. 
%This contradiction comes from the assumption 
%$f^{\prime} (E_1) = f^{\prime} (E_2)$. 
Thus we have $f^{\prime} (E_1) \neq f^{\prime} (E_2)$, 
which completes the proof of the assertion iv).  \qed

\begin{lemma}
If the surface $Y$ is of case $2$%
$)$ in Proposition $\ref{prop:L^2}$, 
then $\lambda = 4$. 
\end{lemma}

Proof. 
By Lemma \ref{lm:L^2-2resoution}, we have $n = 2 \lambda - 2 \leq 9$, 
hence $\lambda \leq 5$. Thus we only need to exclude 
the case $\lambda =5$. 
Assume $\lambda =5$. Then $r : Z^{\prime} \to \mathbb{P}^2$ 
is a blowing-up at one pint, hence  
$Z^{\prime} = \varSigma_1$. 
Thus by Lemmas \ref{lm:lifalilitytcriterion} 
and \ref{lm:fprimeE_i}, we see that 
$\varPhi_L : \tilde{Y} \to Z$ lifts to a morphism 
$f^{\prime} : \tilde{Y} \to Z^{\prime}$, 
and that $f^{\prime} (E_i)$'s are $(-1)$-curves. 
The minimal section $\varDelta_0$, however, 
is the unique $(-1)$-curve on the Hirzebruch surface $\varSigma_1$. 
Thus we have $f^{\prime} (E_1) = f^{\prime} (E_2) =\varDelta_0$, 
which contradicts Lemma \ref{lm:fprimeE_i}. 
Hence we have the assertion.      \qed 
 
Thus we only need to study the case $\lambda =4$. 
In what follows we assume $\lambda =4$, hence $n=6$. 
In this case, the morphism $r : Z^{\prime} \to \mathbb{P}^2$ is 
a blowing-up at three points. 

\begin{lemma}  \label{lm:case2structuer}
Assume that $\varPhi_L : \tilde{Y} \to Z$ lifts to a morphism 
$f^{\prime} : \tilde{Y} \to Z^{\prime}$, and that 
$f^{\prime} (E_1) \cap f^{\prime} (E_2) = \emptyset$ holds. 
Let $r^{\prime} : Z^{\prime} \to W$ denote the blowing-down 
of the two $(-1)$-curves $f^{\prime} (E_1)$ and $f^{\prime} (E_2)$. 
Then the branch divisor $B$ of the morphism 
$r^{\prime} \circ f^{\prime} : \tilde{Y} \to W$ is a 
member of the linear system $|-4 K_W|$ having $[3, 3]$-points 
at $r^{\prime} (f^{\prime} (E_1))$ 
and $r^{\prime} (f^{\prime} (E_2))$.  
Except for these two $[3,3]$-points, the branch divisor $B$  
has at most negligible singularities. 
Further, the surface $Y$ gives the minimal desingularization 
of the double cover 
$($of the surface $W$%
$)$ 
with branch divisor $B$.  
\end{lemma}

Proof. 
Note that $f^{\prime}$ contracts no $(-1)$-curve.  
Thus the divisor ${f^{\prime}}^* B^{\prime} - 2R^{\prime}$, 
linearly equivalent to 
$2(\sum {f^{\prime}}^* (f^{\prime}(E_i)) - 2 \sum E_i)$ by 
Lemma \ref{lm:fprimeE_i}, is twice a certain effective divisor 
$\zeta$ on $\tilde{Y}$.  
We have $\zeta E_j = 
(\sum {f^{\prime}}^* (f^{\prime}(E_i)) - 2 \sum E_i) E_j =1$, 
hence $\sharp (\zeta \cap E_j) = 1$ for $j =1$, $2$. 
So we put $\{ z_j \} = f^{\prime} (\zeta \cap E_j)$. Then the point 
$z_j \in f^{\prime} (E_j)$, where $1 \leq j \leq 2$, 
is an essential singularity of the branch divisor $B^{\prime}$.  
Meanwhile, by Lemma \ref{lm:fprimeE_i}, we see   
that the divisor $B^{\prime} - \sum f^{\prime} (E_i)$ is effective, 
and that $(B^{\prime} - \sum f^{\prime} (E_i)) f^{\prime}(E_j) = 3$ 
for each $j=1$, $2$, 
from which we infer 
$\mathrm{mult}_{z_j} (B^{\prime} - \sum f^{\prime} (E_i)) \leq 3$. 
If, moreover, 
we have $\mathrm{mult}_{z_j} (B^{\prime} - \sum f^{\prime} (E_i)) \leq 2$, 
then $z_j$ is a negligible singularity of the branch divisor $B^{\prime}$; 
the singularity $z_j$ of $B^{\prime}$ decomposes into a sum of points 
of multiplicity at most $2$ by the blowing-up at $z_j$.  
Thus we obtain 
$\mathrm{mult}_{z_j} (B^{\prime} - \sum f^{\prime} (E_i)) = 3$, hence  
$(B^{\prime} - \sum f^{\prime} (E_i)) \cap f^{\prime} (E_j) = \{ z_j\}$
and $\mathrm{mult}_{z_j} B^{\prime} = 4$. 
Let $q_1 \circ q_2 : Z_2^{\prime} \to Z^{\prime}$ be the blowing-up  
at $z_1$ and $z_2$, and $\varepsilon_j = (q_1 \circ q_2)^{-1}(z_j) $, 
the $(-1)$-curve corresponding to $z_j$.  
Then by the same method as in \cite[Section 2]{quintic}, we infer 
that the divisor 
$B_2^{\prime} = (q_1 \circ q_2)^* B^{\prime} - 4 \sum_{i=1,2} \varepsilon_i$ 
has at most negligible singularities, and that 
the surface $\tilde{Y}$ gives the canonical resolution of the double cover 
with branch divisor $B^{\prime}$. 
It follows that  
$B = r^{\prime}_* B^{\prime}$ %,
%which is a member of $|-4K_W|$ by Lemma \ref{lm:fprimeE_i}, 
has $[3,3]$-points at 
$r^{\prime} (f^{\prime}(E_1))$ and $r^{\prime} (f^{\prime}(E_2))$, 
that the divisor $B$ has no essential singularity  
except for these two $[3,3]$-points, 
and that the surface $Y$ gives the minimal desingularization of the 
double cover with branch divisor $B$. 
Now all we have left is the linear equivalence 
$B \sim - 4K_W$, which, however, is 
trivial by iii) in Lemma \ref{lm:fprimeE_i}.  \qed 

\begin{lemma} \label{lm:case2action}
Let $r^{\prime} : Z^{\prime} \to W$ be the blowing-down given in 
Lemma $\ref{lm:case2structuer}$. 
Then the surface $W$ is the Hirzebruch surface $\varSigma_d$ of 
degree $d=0$ or $2$. The action by $G = \mathrm{Gal} (Y/X)$ 
on $Z^{\prime}$ induces one on $W$, of which fixed locus is a 
set of four isolated points. 
Further, none of these four fixed points 
lies on the branch divisor $B$. 
\end{lemma}

Proof. 
The action by $G$ on $Z^{\prime}$ induces one on $W$, 
since the divisor $f^{\prime} (E_1) + f^{\prime} (E_2)$ is 
stable under the action by $G$. 
Note that the anticanonical system $|-K_{Z^{\prime}}|$ has 
no fixed component. From this together with $K_W^2 = K_Z^2 +2 = 8$, 
we see that $W = \varSigma_d$ for a certain integer $0 \leq d \leq 2$. 

Let us show that the class of $\varGamma$, a fiber of the Hirzebruch 
surface $W = \varSigma_d \to \mathbb{P}^1$, is stable under the action 
by $G$ on $W$. If the class of $\varGamma$ is not stable, then we 
see that $d =0$ and that the generator of $G$ maps $\varGamma$ to a 
member of the linear system $|\varDelta_0|$. 
It follows that there exists an irreducible member 
$\varDelta \in |\varDelta_0 + \varGamma|$ contained in the fixed locus 
of the action by $G$ on $W$. We have $\varDelta \cap B \neq \emptyset$, 
since $\varDelta$ is a $2$-curve. Meanwhile since the Galois group $G$ 
acts transitively on the set  
$\{ r^{\prime} (f^{\prime} (E_1)), r^{\prime} (f^{\prime} (E_2)) \}$, 
neither $r^{\prime} (f^{\prime} (E_1))$ nor $r^{\prime} (f^{\prime} (E_2))$ 
lies on $\varDelta$. Thus by Lemma \ref{lm:case2structuer}, every point 
in $\varDelta \cap B$ is at most a negligible singularity of $B$. 
Then the same argument as in the proof of 
Proposition \ref{prop:case1exclusion} leads  us to a contradiction. 
Hence the class of $\varGamma$ is stable. 

Now let us show the assertions. 
The argument above shows that  
the action by $G$ on $W$ induces one on $\mathbb{P}^1$ 
via the natural fibration of 
the Hirzebruch surface $W = \varSigma_d \to \mathbb{P}^1$. 
Note that if this induced action on $\mathbb{P}^1$ is trivial, 
then there exists an irreducible member 
$\varDelta_1 \in |\varDelta_0 + d \varGamma|$ contained in the 
fixed locus of the action by $G$ on $W$, which, 
together with the same argument as in the case of $\varDelta$ above, 
leads us to a contradiction.
Thus the induced action on $\mathbb{P}^1$ is non-trivial. 
It follows that $|\varGamma|$ has exactly two members   
stable under the action on $W$, 
which we shall call $\varGamma_1$ and $\varGamma_2$. 
The same argument as in the case of $\varDelta$ above shows that 
the induced action on $\varGamma_i$ is non-trivial for each $i =1$, $2$. 
Thus we see that $d \neq 1$, and that if $d=0$ or $2$, then the fixed locus 
of the induced action on $W$ is a set of four isolated points.  
The absence of the fixed points lying on $B$ follows from 
the same argument as in the case of $\varDelta$ above. 
\qed

By Lemmas \ref{lm:case2structuer} and \ref{lm:case2action}, 
we see that if $\varPhi_L$ lifts to $f^{\prime} : \tilde{Y} \to Z^{\prime}$, 
and if $f^{\prime} (E_1) \cap f^{\prime} (E_2) = \emptyset$, 
then our surface $X$ has the structure as in 
Theorem \ref{thm:completedescription}.  
Let us check that these two conditions are in fact satisfied. 
To do this,  
we study the arrangement of $(-1)$-curves and $(-2)$-curves 
on $Z^{\prime}$, and use Lemmas \ref{lm:lifalilitytcriterion} 
and \ref{lm:fprimeE_i}.  

Let $r_i : Z^{\prime}_i \to Z^{\prime}_{i-1}$, where 
$-2 \leq i \leq 0$, be the blowing-up such that 
$r = (r_{-2} \circ r_{-1} \circ r_{0}) : 
Z^{\prime}_0 = Z^{\prime} \to Z^{\prime}_{-3} = \mathbb{P}^2 $ holds. 
We denote by $z_i \in Z^{\prime}_{i-1}$ 
and $\varepsilon_i = r_i^{-1} (z_i)$ the center of the blowing-up $r_i$ 
and its corresponding $(-1)$-curve, respectively. 
For each $-2 \leq i \leq 0$, we denote by $\varepsilon_i^{\prime}$   
the strict transform to $Z^{\prime}$ 
of the exceptional curve $\varepsilon_i$.
For the total transform to $Z^{\prime}$ of $\varepsilon_i$, 
we use the same symbol $\varepsilon_i$.  

\begin{lemma}
Let $m \leq 2$ be a non-negative integer, 
and $C$, a $(-m)$-curve on $Z^{\prime}$ not exceptional 
with respect to $r: Z^{\prime} \to \mathbb{P}^2$. 
Then $C$ is a strict transform to $Z^{\prime}$ of a line on $\mathbb{P}^2$ 
passing exactly $m+1$ of the tree points $z_{i}$'s   
$($%
$-2 \leq i \leq 0$%
$)$. 
\end{lemma}

Proof. 
Let $l$ be a line on $\mathbb{P}^2$. 
Then we have 
$C \sim m_0 r^* (l)  - \sum_{i = -2}^0 n_i \varepsilon_i$ for 
certain integers $m_0 \geq 1$ and $n_i \geq 0$'s.
Note that $C^2 = -m$ and $-K_{Z^{\prime}}C = 2 -m$, 
since $C$ is a $(-m)$-curve. 
Thus we have  
\begin{equation} \label{eq:CK}
m_0^2 - \sum_{i= -2}^0 n_i^2 = -m,  \quad \quad 
3 m_0 - \sum_{i= -2}^0 n_i = -m + 2. 
\end{equation}
From these equalities, we infer 
\[
 5 \sum_{i= -2}^0 n_i^2 
 + \sum_{-2 \leq i < j \leq 0} (n_i - n_j)^2 
 + \sum_{i= -2}^0 (n_i + m-2)^2 
 = 9m + 4 (m-2)^2 \leq 18,
\]
hence $\sum_{i= -2}^0 n_i^2 \leq 3$. 
Thus we have $n_i^2 = n_i$ for any $-2 \leq i \leq 0$. 
By this together with the equalities (\ref{eq:CK}), we obtain 
$m_0 = 1$ and $\sum_{i= -2}^0 n_i = m +1 $. 
Thus we have the assertion. \qed

We study the arrangement of $(-1)$-curves and $(-2)$-curves on 
$Z^{\prime}$ according to the configuration of the centers $z_i$'s of 
the blowing-up $r: Z^{\prime} \to \mathbb{P}^2$.
First, we study the case where no two of the centers $z_{-2}$, 
$z_{-1}$, and $z_{0}$ are infinitely near. This case is divided 
into two cases: case $2$-$1$-$1$) and case $2$-$1$-$2$). 

Case $2$-$1$-$1$): the case where the centers $z_{-2}$, 
$z_{-1}$, and $z_0$ are not collinear. 
In this case, the surface $Z^{\prime}$ has no $(-2)$-curve. 
Thus $\varPhi_L$ lifts to 
a morphism $f^{\prime} : \tilde{Y} \to Z^{\prime}$. 
There exist exactly six $(-1)$-curves: 
$\varepsilon_{-2}$, $\varepsilon_{-1}$, $\varepsilon_{0}$, 
$r^{-1}_*(l_{-2,-1})$, $r^{-1}_*(l_{-1,0})$, and 
$r^{-1}_*(l_{-2,0})$, where $l_{i,j}$ denotes the line on 
$\mathbb{P}^2$ passing $z_i$ and $z_j$. 
Let $(X_0 : X_1: X_2)$ be homogeneous coordinates of 
$\mathbb{P}^2$ satisfying 
$z_{-2} = (1:0:0)$, $z_{-1} = (0:1:0)$, and $z_0 = (0:0:1)$. 
For each $(a,b) \in \mathbb{C}^{\times} \times \mathbb{C}^{\times}$, 
we denote by $\varphi_{(a,b)}$ the automorphism of $Z^{\prime}$ 
corresponding to the projective transformation  
$(X_0 : X_1: X_2) \mapsto (X_0 : a X_1: b X_2)$. 

Let us study the induced action by $G$ on $Z^{\prime}$. 
Let $\mathrm{Aut} (Z^{\prime})$ be the group of analytic automorphisms 
of the surface $Z^{\prime}$, and $D_6$, the dihedral group of 
degree $6$. Then we have a short exact sequence 
\[
 0 \to \mathbb{C}^{\times} \times \mathbb{C}^{\times}
   \to \mathrm{Aut}(Z^{\prime}) 
   \to D_6 \to 0,
\]
where the morphism 
$\mathbb{C}^{\times} \times \mathbb{C}^{\times} 
\to \mathrm{Aut}(Z^{\prime}) $ 
is given by $(a,b) \mapsto \varphi_{(a,b)}$, and 
the morphism $\mathrm{Aut}(Z^{\prime}) \to D_6$ 
corresponds to the transitions of six $(-1)$-curves on $Z^{\prime}$.  
Let $\varphi_{\sigma}$ and $\varphi_{\tau}$ be 
the automorphisms of $Z^{\prime}$ 
corresponding to the Cremona transformation 
$(X_0 : X_1: X_2) \mapsto 
(X_2 X_0: X_0 X_1: X_1 X_2)$ and the morphism 
$(X_0 : X_1: X_2) \mapsto (X_0 : X_2: X_1)$, respectively. 
Then we have 
\[
 (\varphi_{\sigma})^6 = \mathrm{id}_{Z^{\prime}} \qquad
 (\varphi_{\tau})^2   = \mathrm{id}_{Z^{\prime}} \qquad 
 \varphi_{\sigma} \circ \varphi_{\tau}
 \circ \varphi_{\sigma} \circ \varphi_{\tau} = \mathrm{id}_{Z^{\prime}}.
\]
Thus the short exact sequence above splits. 
We denote by $\sigma$ and $\tau$ the image by 
$\mathrm{Aut} (Z^{\prime}) \to D_6$ of  
$\varphi_{\sigma}$ and $\varphi_{\tau}$, respectively. 
We have a group homomorphism 
$G \to \mathrm{Aut} (Z^{\prime})$ 
corresponding to the action by $G$ on $Z^{\prime}$. 
Composing this homomorphism with 
$\mathrm{Aut} (Z^{\prime}) \to D_6$, 
we obtain a group homomorphism 
$\alpha : G \to D_6$. 
Note that by Lemma \ref{lm:fprimeE_i},  
the morphism $\alpha$ is an injection of $G$ into $D_6$. 
Hence the image $\alpha (G)$ is conjugate to  
$\langle \tau \rangle$, $ \langle \sigma^3 \tau  \rangle$,  
or $\langle \sigma^3 \rangle$ in $D_6$. 

Assume that the image $\alpha (G)$ is conjugate 
to $\langle \tau \rangle$ in $D_6$. 
Replacing the morphism $r : Z^{\prime} \to \mathbb{P}^2$
if necessary, we may assume that $\alpha (G) = \langle \tau \rangle$. 
Then since the Galois group $G$ acts transitively on the set 
$\{ f^{\prime} (E_1), f^{\prime} (E_2) \}$, 
we have  
$\{ f^{\prime} (E_1), f^{\prime} (E_2) \} 
= \{ r^{-1}_* (l_{-2, -1}), r^{-1}_* (l_{-2, 0}) \}$  
or $\{ \varepsilon_{-1}, \varepsilon_0 \}$, 
hence $f^{\prime} (E_1) \cap f^{\prime} (E_2) = \emptyset$. 
It follows that the surface $W$, where $r^{\prime} : Z^{\prime} \to W$ is the 
blowing-down of the two $(-1)$-curves   
$f^{\prime} (E_1)$ and $f^{\prime} (E_2)$, 
is isomorphic to the Hirzebruch surface $\varSigma_1$, 
which contradicts Lemma \ref{lm:case2action}. 
Thus $\alpha (G)$ is not conjugate to $\langle \tau \rangle$. 

Assume that the image $\alpha (G)$ is conjugate to 
$\langle \sigma^3 \tau \rangle$ in $D_6$. 
Replacing the morphism $r : Z^{\prime} \to \mathbb{P}^2$ 
if necessary, we may assume that 
$\alpha (G) = \langle \sigma^3 \tau \rangle$.  
Then the blowing-down of the two $(-1)$-curves 
$\varepsilon_{-2}$ and $r^{-1}_* (l_{-1, 0})$ gives 
a birational morphism 
$r^{\prime \prime} : Z^{\prime} \to 
\varSigma_0 = \mathbb{P}^1 \times \mathbb{P}^1 $ 
satisfying 
$r^{\prime \prime}_* (\varepsilon_0) \sim 
 r^{\prime \prime}_* (r^{-1}_* (l_{-2, -1})) \sim \varDelta_0$ and  
$r^{\prime \prime}_* (\varepsilon_{-1}) \sim 
 r^{\prime \prime}_* (r^{-1}_* (l_{-2, 0})) \sim \varGamma$. 
Note that the action by $G$ on $Z^{\prime}$ induces 
one on $\varSigma_0 = \mathbb{P}^1 \times \mathbb{P}^1$. 
We take homogeneous coordinates 
$((\xi_0 : \xi_1), (\eta_0 : \eta_1))$ of 
$\varSigma_0 = \mathbb{P}^1 \times \mathbb{P}^1$ 
in such a way that 
$r^{\prime \prime} (\varepsilon_{-2}) 
= ((1 : 0), (1 : 0))$ and 
$r^{\prime \prime} (r^{-1}_* (l_{-1, 0})) 
= ((0 : 1), (0 : 1))$ hold, 
and that 
the automorphism of $\varSigma_0$ corresponding to 
the generator $\iota$ of $G$ is given by 
$((\xi_0 : \xi_1), (\eta_0 : \eta_1)) \mapsto 
((\eta_1 : \eta_0), (\xi_1 : \xi_0))$.   
Since we have $-K_{Z^{\prime}} \sim 
{r^{\prime \prime}}^* (-K_{\varSigma_0}) 
- \varepsilon_{-2} - r^{-1}_* (l_{-1, 0})$, 
the space $H^0(\mathcal{O}_{Z^{\prime}} (-2K_{Z^{\prime}}))$ 
corresponds to a certain subspace 
$V$ of $H^0(\mathcal{O}_{\varSigma_0} (-2K_{\varSigma_0}))$. 
Every element in $V$ is a homogeneous polynomial 
$\psi (\xi_0, \xi_1, \eta_0, \eta_1)$ of bidegree $(4, 4)$ 
vanishing with multiplicity at least $2$ at 
$((1 : 0), (1 : 0))$ and $((0 : 1), (0 : 1))$.

Note that we have a natural inclusion 
$V \hookrightarrow H^0(\mathcal{O}_Y (2K_Y))$,   
since we have 
$L \sim {f^{\prime}}^* D^{\prime}$.   
We denote by $\phi$ the restriction to $V$ 
of the natural action by $G$ on $H^0_Y (\mathcal{O}(2K_Y))$. 
Let $\phi^{\prime} (\iota)$ be the automorphism of $V$ given 
by $\psi (\xi_0, \xi_1, \eta_0, \eta_1) \mapsto
\psi (\eta_1, \eta_0, \xi_1, \xi_0)$.     
Then $\iota \mapsto \phi^{\prime} (\iota)$,  
where $\iota$ is the generator of the Galois group $G$,   
gives another action $\phi^{\prime}$ by $G$ on $V$.  
Note that for any $g \in G$ and $\psi \in V$, the two elements 
$\phi (g) \psi$ and $\phi^{\prime} (g) \psi$ defines the 
same divisor on $\varSigma_0$. 
From this we infer that $\phi = c \phi^{\prime}$ for 
a certain character $c \in \mathrm{Char} (G)$.

Now let $V^+$ be the set of all elements in $V$ stable under the 
action $\phi^{\prime}$. Then by $\phi = c \phi^{\prime}$, 
we see that $V^+ \subset H^0(\mathcal{O}_X (2K_X - T_c))$ 
for a torsion divisor $T_c \in \mathrm{Pic} (X)$ 
corresponding to the character $c$. 
Meanwhile, by the Riemann--Roch theorem, we have 
$h^0(\mathcal{O}_X (2K_X - T_c)) = \chi + K_X^2 = 11$. 
The space $V^+$, however, has a base consisting of  
twelve elements: 
\[
 \xi_0^i \xi_1^{4-i} \eta_0^j \eta_1^{4-j} 
  + \xi_0^{4-j} \xi_1^{j} \eta_0^{4-i}\eta_1^{i}
  \qquad 
  ( 0 \leq i, \quad 0 \leq j, \quad 2 \leq i+j \leq 4).   
\] 
This contradicts the inequality 
$\dim V^+ \leq  h^0(\mathcal{O}_X (2K_X - T_c))$.   
Hence, the image $\alpha (G)$ is not conjugate 
to $\langle \sigma^3 \tau \rangle$ in $D_6$. 

Thus we have $\alpha (G) = \langle \sigma^3 \rangle$. 
Hence, replacing $r : Z^{\prime} \to \mathbb{P}^2$ if necessary, 
we may assume that 
$\{ f^{\prime} (E_1), f^{\prime} (E_2) \} = 
\{ \varepsilon_{-2}, r^{-1}_* (l_{-1, 0}) \}$. 
Then the surface $W$ as in Lemma \ref{lm:case2structuer}, 
obtained by blowing down the two $(-1)$-curves $f^{\prime} (E_1)$  
and $f^{\prime} (E_2)$ of $Z^{\prime}$, 
is isomorphic to the Hirzebruch surface $\varSigma_0$.  
Thus by Lemmas \ref{lm:case2structuer} and \ref{lm:case2action}, 
our surface $X$, in case $2$-$1$-$1$), has the structure 
as in the case $d=0$ in Theorem \ref{thm:completedescription}.  
 
Case $2$-$1$-$2$): 
the case where three points $z_{-2}$, $z_{-1}$, and $z_0$ 
are collinear. Let $l_{-2, -1}$ be the line on $\mathbb{P}^2$ passing 
the tree points $z_{i}$'s above. Then the strict transform 
$r^{-1}_* (l_{-2, -1})$ is the unique $(-2)$-curve on $Z^{\prime}$. 
Hence by Lemma \ref{lm:lifalilitytcriterion}, the morphism 
$\varPhi_L : \tilde{Y} \to Z$ lifts to 
$f^{\prime} : \tilde{Y} \to Z^{\prime}$. 
Meanwhile the surface $Z^{\prime}$ has exactly three $(-1)$-curves:
$\varepsilon_{-2}$, $\varepsilon_{-1}$, and $\varepsilon_0$. 
Replacing $r : Z^{\prime} \to \mathbb{P}^2$ if necessary, 
we may assume 
$\{ f^{\prime} (E_1), f^{\prime} (E_2)\} 
= \{ \varepsilon_{-2}, \varepsilon_{-1} \}$ 
by Lemma \ref{lm:fprimeE_i}.  
Let $r^{\prime} : Z^{\prime} \to W$ be the blowing-down as 
in Lemma \ref{lm:case2structuer} of the two $(-1)$-curves 
$f^{\prime} (E_1)$ and $f^{\prime} (E_2)$. 
Then we have $W = \varSigma_1$, 
$r^{\prime}_* (\varepsilon_0) = \varDelta_0$, and 
$r^{\prime}_* (r^{-1}_*(l_{-2, -1})) \sim \varGamma$, 
which contradicts Lemma \ref{lm:case2action}. 
Thus case $2$-$1$-$2$) does not occur.   

Next, we study the case where  
$z_{-2}$ and $z_{-1}$ are distinct points on $\mathbb{P}^2$, 
and $z_0$ is infinitely near to $z_{-1}$. 
We denote by $l_{-2, -1}$ the unique line on $\mathbb{P}^2$ 
passing $z_{-2}$ and $z_{-1}$. 
This case is divided into two cases: 
case $2$-$2$-$1$) and case $2$-$2$-$2$). 

Case $2$-$2$-$1$): the case where $z_0$ does not lie on the 
strict transform $(r_{-2} \circ r_{-1})^{-1}_* (l_{-2, -1})$ 
of $l_{-2, -1}$ by $r_{-2} \circ r_{-1}$.
Let $l_{-1, 0}$ be the unique line on $\mathbb{P}^2$ whose 
strict transform $(r_{-2} \circ r_{-1})^{-1}_* (l_{-1, 0})$
by $r_{-2} \circ r_{-1}$ passes $z_0$. 
Then the surface $Z^{\prime}$ has 
a unique $(-2)$-curve $\varepsilon_{-1}^{\prime}$, 
and exactly four $(-1)$-curves 
$\varepsilon_{-2}$, $\varepsilon_0$, 
$r^{-1}_* (l_{-2, -1})$, and $r^{-1}_* (l_{-1, 0})$. 
Hence, by Lemma \ref{lm:lifalilitytcriterion}, the morphism 
$\varPhi_L : \tilde{Y} \to Z$ lifts to 
$f^{\prime} : \tilde{Y} \to Z^{\prime}$. 
Note that $\{ \varepsilon_0 , r^{-1}_* (l_{-2, -1})\}$ is the 
set of all $(-1)$-curves intersecting the 
unique $(-2)$-curve $\varepsilon_{-1}^{\prime}$. 
Thus we have 
$\{ f^{\prime} (E_1), f^{\prime} (E_2)\} 
= \{ \varepsilon_0 , r^{-1}_* (l_{-2, -1}) \}$ or 
$\{ \varepsilon_{-2} , r^{-1}_* (l_{-1, 0}) \}$, 
hence, in particular, 
$f^{\prime} (E_1) \cap f^{\prime} (E_2) = \emptyset$. 
We denote by $r^{\prime} : Z^{\prime} \to W$ the 
blowing-down as in Lemma \ref{lm:case2structuer}  
of the two $(-1)$-curves $f^{\prime} (E_1)$ and $f^{\prime} (E_2)$. 
If $\{ f^{\prime} (E_1), f^{\prime} (E_2)\} 
= \{ \varepsilon_0 , r^{-1}_* (l_{-2, -1}) \}$, then we have 
$W = \varSigma_0$, 
$r^{\prime}_* (\varepsilon_{-1}^{\prime}) \sim \varDelta_0$ 
and 
$r^{\prime}_* (\varepsilon_{-2}) \sim
 r^{\prime}_* (r^{-1}_* (l_{-1, 0}))  \sim  \varGamma$. 
If on the other hand  
$\{ f^{\prime} (E_1), f^{\prime} (E_2)\} 
= \{ \varepsilon_{-2} , r^{-1}_* (l_{-1, 0}) \}$, 
then we have 
$W = \varSigma_2$, 
$r^{\prime}_* (\varepsilon_{-1}^{\prime}) = \varDelta_0$, and   
$r^{\prime}_* (\varepsilon_0) \sim
 r^{\prime}_* (r^{-1}_* (l_{-2, -1}))  \sim  \varGamma$. 
Thus by lemmas \ref{lm:case2structuer} and \ref{lm:case2action}, 
our surface $X$, in case $2$-$2$-$1$), has 
the structure as in the case $d=0$ or the case $d=2$ 
in Theorem \ref{thm:completedescription},  
according as  
$\{ f^{\prime} (E_1), f^{\prime} (E_2)\} = 
\{ \varepsilon_0 , r^{-1}_* (l_{-2, -1}) \}$ or 
$\{ f^{\prime} (E_1), f^{\prime} (E_2)\} = 
\{ \varepsilon_{-2} , r^{-1}_* (l_{-1, 0}) \}$   
respectively.  
   
Case $2$-$2$-$2$): 
the case where $z_0$ lies on the strict transform 
$(r_{-2} \circ r_{-1})^{-1}_* (l_{-2, -1})$. 
In this case, the surface $Z^{\prime}$ has 
exactly two $(-2)$-curves 
$\varepsilon_{-1}^{\prime}$ and $r^{-1}_* (l_{-2, -1})$, 
and exactly two $(-1)$-curves 
$\varepsilon_{-2}$ and $\varepsilon_0$. 
Note that every $(-2)$-curve on $Z^{\prime}$ is 
stable under the action by $G$ on $Z^{\prime}$; 
the divisor $r^{-1}_* (l_{-2, -1})$ is the unique  
$(-2)$-curve intersecting all $(-1)$-curves on $Z^{\prime}$.  
Thus by Lemma \ref{lm:lifalilitytcriterion}, the morphism 
$\varPhi_L : \tilde{Y} \to Z$ lifts 
to $f^{\prime} : \tilde{Y} \to Z^{\prime}$ . 
Then it follows from Lemma \ref{lm:fprimeE_i} that 
$\{ f^{\prime} (E_1), f^{\prime} (E_2)\} 
= \{ \varepsilon_{-2} , \varepsilon_0 \}$. 
This, however, contradicts the transitivity of the action 
by $G$ on $\{ f^{\prime} (E_1), f^{\prime} (E_2)\}$, 
since $\varepsilon_0$ is the unique $(-1)$-curve intersecting 
all $(-2)$-curves on $Z^{\prime}$. 
Thus case $2$-$2$-$2$) does not occur. 

Finally, we study the case where  
$z_{-1}$ is infinitely near to $z_{-2}$, 
and $z_0$ is infinitely near to $z_{-1}$. 
We denote by $l_{-2, -1}$ the unique line on $\mathbb{P}^2$ 
whose strict transform $(r_{-2})^{-1}_* (l_{-2, -1})$ passes $z_{-1}$. 
Note that $Z^{\prime}$ has no $(-3)$-curve, since 
the linear system $|-K_{Z^{\prime}}|$ has no fixed component. 
Thus the point $z_0$ does not lie on the strict transform 
$(r_{-1})^{-1}_* (\varepsilon_{-2})$. 
This case is divided into two cases: 
case $2$-$3$-$1$) and case $2$-$3$-$2$). 

Case $2$-$3$-$1$): 
the case where $z_0$ does not lie on the strict transform 
$(r_{-2} \circ r_{-1})^{-1}_* (l_{-2, -1})$. 
In this case, the surface $Z^{\prime}$ has exactly 
two $(-2)$-curves 
$\varepsilon_{-2}^{\prime}$ and $\varepsilon_{-1}^{\prime}$, 
and exactly two $(-1)$-curves 
$\varepsilon_0$ and $r^{-1}_* (l_{-2, -1})$. 
Since $\varepsilon_{-2}^{\prime}$ is the unique 
$(-2)$-curve intersecting no $(-1)$-curve on $Z^{\prime}$, 
every $(-2)$-curve is stable under the action by $G$ on $Z^{\prime}$. 
Thus by Lemma \ref{lm:lifalilitytcriterion}, 
the morphism $\varPhi_L : \tilde{Y} \to Z$ lifts to 
$f^{\prime} : \tilde{Y} \to Z^{\prime}$.   
Then it follows from Lemma \ref{lm:fprimeE_i} that 
$\{ f^{\prime} (E_1), f^{\prime} (E_2)\} 
= \{ \varepsilon_0 , r^{-1}_* (l_{-2, -1})\}$, 
hence $f^{\prime} (E_1) \cap f^{\prime} (E_2) = \emptyset $. 
Let $r^{\prime} : Z^{\prime} \to W$ be the blowing-down 
as in Lemma \ref{lm:case2structuer} of the two $(-1)$-curves 
$f^{\prime} (E_1)$ and $f^{\prime} (E_2)$.  
Then we have 
$W = \varSigma_2$, 
$r^{\prime}_* (\varepsilon_{-2}^{\prime}) = \varDelta_0$, 
and $r^{\prime}_* (\varepsilon_{-1}^{\prime}) \sim \varGamma$. 
Thus by Lemmas \ref{lm:case2structuer} and \ref{lm:case2action}, 
our surface $X$, in case $2$-$3$-$1$), has the structure 
as in the case $d=2$ in Theorem \ref{thm:completedescription}. 

Case $2$-$3$-$2$): 
the case where $z_0$ lies on the strict transform 
$(r_{-2} \circ r_{-1})^{-1}_* (l_{-2, -1})$. 
In this case, the surface $Z^{\prime}$  
has exactly three $(-2)$-curves 
$\varepsilon_{-2}^{\prime}$, $\varepsilon_{-1}^{\prime}$, 
and $r^{-1}_* (l_{-2, -1})$, 
and a unique $(-1)$-curve 
$\varepsilon_0$.    
Note that $\varepsilon_{-2}^{\prime}$ is the unique $(-2)$-curve 
intersecting no $(-1)$-curve on $Z^{\prime}$, 
and that $\varepsilon_{-1}^{\prime}$ is the unique $(-2)$-curve 
intersecting $\varepsilon_{-2}^{\prime}$. 
Thus every $(-2)$-curve on $Z^{\prime}$ is stable under 
the action by $G$ on $Z^{\prime}$.  
Thus by Lemma \ref{lm:lifalilitytcriterion}, 
the morphism $\varPhi_L : \tilde{Y} \to Z$ lifts to 
$f^{\prime} : \tilde{Y} \to Z^{\prime}$.   
This, however, contradicts Lemma \ref{lm:fprimeE_i},  
since $\varepsilon_0$ is the unique $(-1)$-curve on $Z^{\prime}$. 
Hence, case $2$-$3$-$2$) does not occur. 

Thus we have the following: 

\begin{proposition}  \label{prop:case2conclusion}
Assume that the surface $Y$ is of case $2$%
$)$ in Proposition 
$\ref{prop:L^2}$. Then $\lambda = 4$. 
Further,  
the surface $X$ in this case  
has the structure 
as in Theorem $\ref{thm:completedescription}$. 
\end{proposition}

\section{The case $\deg \varPhi_{K_Y} = 1$} \label{scn:deg=1}

In this section, we exclude the case $\deg \varPhi_{K_Y} =1$ 
and give a proof for Theorems \ref{thm:maintheorem} and 
\ref{thm:completedescription}. 
In what follows, we assume that $\deg \varPhi_{K_Y} =1$. 
Note that by Proposition \ref{prop:degphiKY}, we have $\lambda =4$, 
hence $K_Y^2 = 14$, $p_g (Y) =7$, and $q (Y) =0$. 
Thus our $Y$ is a canonical surface whose invariant lies on the 
Castelnuovo line. 
By \cite[Lemma 1.1]{3pg-7}, 
the canonical system $|K_Y|$ is free from base points;  
hence the canonical map 
$\varPhi_{K_Y} : Y \to \mathbb{P}^n$ is a morphism, 
where $n = 2\lambda -2 =6$.  
In what follows, we frequently use results given in \cite{3pg-7}.  

Let $\mathcal{Q} \subset \mathbb{P}^n$ be the intersection 
of all quadrics containing the canonical image $Z = \varPhi_{K_Y} (Y)$. 
By \cite[Section $1$]{3pg-7}, we obtain the following: 

\begin{proposition} \label{prop:intersectionofquadrics}
Let $\mathcal{Q}$ be the variety defined above. Then either  
of the following holds$:$

$1$%
$)$ $\mathcal{Q}$ is the image by $\varPhi_T$ of the variety 
$\mathcal{Q}^{\prime} = 
\mathbb{P} (\mathcal{O}_{\mathbb{P}^2} 
\oplus \mathcal{O}_{\mathbb{P}^2} (2))$, 
where $\varPhi_T$ is the morphism associated with 
a tautological divisor $T$ of the $\mathbb{P}^1$-bundle 
$\mathrm{pr}_{\mathcal{Q}^{\prime}} : \mathbb{P} (\mathcal{O}_{\mathbb{P}^2} 
\oplus \mathcal{O}_{\mathbb{P}^2} (2)) \to \mathbb{P}^2$%
$;$  

$2$%
$)$
$\mathcal{Q}$ is the image by $\varPhi_T$ of the variety 
$\mathcal{Q}^{\prime} =
\mathbb{P} (\bigoplus_{i=0}^2 \mathcal{O}_{\mathbb{P}^1} (a_i))$, 
where $\varPhi_T$ is the morphism associated with 
a tautological divisor $T$ of the $\mathbb{P}^2$-bundle 
$\mathrm{pr}_{\mathcal{Q}^{\prime}} : 
\mathbb{P} (\bigoplus_{i=0}^2 \mathcal{O}_{\mathbb{P}^1} (a_i))
\to \mathbb{P}^1$,  
and 
$0 \leq a_0 \leq a_1 \leq a_2$ and $\sum_{i=0}^2 a_i = n-2$.  
\end{proposition}

First, we exclude case $1$) in the proposition above.  

\begin{proposition}   \label{prop:deg1case1exclusion}
Case $1$%
$)$ in Proposition $\ref{prop:intersectionofquadrics}$ does not occur. 
\end{proposition}

Proof. Assume that our $\mathcal{Q}$ is as in case 
$1$) in Proposition \ref{prop:intersectionofquadrics}. 
Then $\mathcal{Q}$ is a cone over the Veronese surface.  
Let $p_0$ be the vertex of $\mathcal{Q}$, and 
$\varLambda$, the linear system consisting of pull-backs 
by $\varPhi_{K_Y}$ of all hyperplanes in $\mathbb{P}^n$ 
passing $p_0$.  
We denote by $\varLambda_0$ and $G_0$ its variable part and 
fixed part respectively. 
By \cite[Proof of Claim I]{3pg-7}, the linear system $\varLambda_0$ 
is free from base points and  induces  
$\varPhi_{\varLambda_0} : Y \to \mathbb{P}^{n-1}$, 
a morpshism of mapping degree $3$ onto its image. 
The image 
$\varPhi_{\varLambda_0} (Y)$ is the Veronese surface, 
i.e., the projective plane $\mathbb{P}^2$ embedded in 
$\mathbb{P}^5$ by $\mathcal{O}_{\mathbb{P}^2} (2)$. 
Note that by the definition of $\mathcal{Q}$, 
the varity $\mathcal{Q}$ and its vertex $p_0$ 
are stable under the action by $G = \mathrm{Gal} (Y/X)$ on 
$\mathbb{P}^n$.  
This implies that the subspace of 
$H^0 (\mathcal{O}_Y (K_Y))$ corresponding to $\varLambda$ is 
stable under the action by $G$ on $H^0 (\mathcal{O}_Y (K_Y))$. 
Thus the action by $G$ on $Y$ induces one on  
$\varPhi_{\varLambda_0} (Y) = \mathbb{P}^2$. 
Now let us derrive a contradiction. 
Since $G \simeq \mathbb{Z} / 2$, the 
fixed locus of this induced action contains a line $l_0$  
on $\mathbb{P}^2$. 
Then the divsor $\varPhi_{\varLambda_0}^* (l_0)$, stable 
under the action by $G$, is a pull-back by 
$\pi : Y \to X$ of that on $X$. 
We however have $\varPhi_{\varLambda_0}^* (l_0)^2 
= \deg \varPhi_{\varLambda_0} =3$, which contradicts $\deg \pi =2$. 
Thus we have the assertion.  \qed

Next, we exclude case $2$) 
in Proposition \ref{prop:intersectionofquadrics}. 

\begin{lemma} \label{lm:a0=0}
If the variety $\mathcal{Q}$ is as in case $2$%
$)$ of Proposition \ref{prop:intersectionofquadrics},   
then $a_0 =0$. 
\end{lemma}

Proof. 
Assume that our variety $\mathcal{Q}$ is as in case $2$) 
in Proposintion \ref{prop:intersectionofquadrics}  
and that $a_0 > 0$. 
Then $\varPhi_T : \mathcal{Q}^{\prime} \to \mathbb{P}^n$ is an embedding. 
We identify $\mathcal{Q}$ and $\mathcal{Q}^{\prime}$ by $\varPhi_T$.  
By the same arguement as in the proof of 
Proposintion \ref{prop:deg1case1exclusion}, we see that 
the variety $\mathcal{Q}$ 
is stable under the action by $G$ on $\mathbb{P}^n$. %;   
%hence the action on $\mathbb{P}^n$ induces one on $\mathcal{Q}$.   
Let $P$ be a fiber of the $\mathbb{P}^2$-bundle 
$\mathrm{pr}_{\mathcal{Q}^{\prime}} 
: \mathcal{Q}=\mathcal{Q}^{\prime} \to \mathbb{P}^1$. 
Then $P$ and $T$ generate 
the Picard group of $\mathcal{Q}$. 
Using this, we see easily    
that if a divisor $P^{\prime}$ on $\mathcal{Q}$ satisfies  
${P^{\prime}}^3 = K_{\mathcal{Q}}{P^{\prime}}^2 =0$ and 
$h^0 (\mathcal{O}_{\mathcal{Q}} (P^{\prime})) =2$, then 
$P^{\prime} \sim P$.    
Thus the class of $P$ is stable 
under the action by $G$ on $\mathcal{Q}$.  
It follows that this action induces one on $\mathbb{P}^1$ 
via the projection 
$\mathrm{pr}_{\mathcal{Q}^{\prime}} 
: \mathcal{Q}=\mathcal{Q}^{\prime} \to \mathbb{P}^1$, 
and that there exsits a member $P_0 \in |P|$ stable under the 
action on $\mathcal{Q}$. Now let us derrive a contradiction.    
Since $G \simeq \mathbb{Z} / 2$, the fixed locus of the 
action by $G$ on $P_0 = \mathbb{P}^2$ contains a line $l_0$.        
Hence the action on $Z$ has a fixed points. 
By \cite[Theorem 1.5]{3pg-7}, however, the surface $Z$ has at most 
rational double points as its singularities.  
Thus, by the same arguement as 
in the proof of Proposition \ref{prop:case1exclusion}, 
we infer that the action on $Y$ has fixed points,  
which contradicts the definition of $\pi : Y \to X$. \qed 

\begin{proposition} \label{prop:prop:deg1case2exclusion}
Case $2$%
$)$ in Proposition $\ref{prop:intersectionofquadrics}$ 
does not occur. 
\end{proposition}

Proof. 
Assume that our $\mathcal{Q}$ is as in case $2$) 
in Proposition \ref{prop:intersectionofquadrics}. 
Then by Lemma \ref{lm:a0=0} and \cite[Claim II]{3pg-7}, 
we have $a_0 = 0$  and $a_1 > 0$. 
It follows that our $\mathcal{Q}$ is a cone over 
the Hirzebruch surface 
$\varSigma_{a_2 - a_1}$ 
embedded in $\mathbb{P}^{n-1}$ by 
$|\varDelta_0 + a_2 \varGamma|$. 
Let $p_0$ be the vertex of $\mathcal{Q}$, 
and $\varLambda$, the linear system 
consisting of the pull-backs by $\varPhi_{K_Y}$ of 
all hyperplanes in $\mathbb{P}^n$ passing $p_0$. 
We denote by $\varLambda_0$ and $G_0$ the variable part 
and the fixed part of $\varLambda$ respectively. 
By \cite[Proof of Claim II]{3pg-7}, the linear system $\varLambda_0$ 
is free from base points and induces 
$\varPhi_{\varLambda_0} : Y \to \mathbb{P}^{n-1}$, 
a morphism of degree $3$ onto its image. 
The image $\varPhi_{\varLambda_0} (Y)$ is the Hirzebruch surface 
$\varSigma_{a_2 - a_1}$ embedded in $\mathbb{P}^{n-1}$ by 
$|\varDelta_0 + a_2 \varGamma|$. 
By the same arguement as in the proof of 
Proposition \ref{prop:deg1case1exclusion}, 
we see that the action by $G$ on $Y$ induces 
one on $\varSigma_{a_2 - a_1}$. 
The class of $\varDelta_0 + a_2 \varGamma$ and 
that of $ -K_{\varSigma_{a_2 - a_1}}$ are stable under this induced 
action on $\varSigma_{a_2 - a_1}$;
hence so are the class of $\varDelta_0$ and 
that of $\varGamma$.   
Thus there exist members $\varDelta_1 \in |\varDelta_0|$ 
and $\varGamma_1 \in |\varGamma|$ stable under the action on
$\varSigma_{a_2 - a_1}$. 
Then from 
$\varPhi_{\varLambda_0}^* (\varDelta_1) 
 \varPhi_{\varLambda_0}^* (\varGamma_1)  
 = \deg \varPhi_{\varLambda_0} =3$,  
we derrive a contradiction
by the same arguement as in the proof 
of Proposition \ref{prop:deg1case1exclusion}. 
Thus we have the assertion. \qed 

Now we are ready to prove Theorems \ref{thm:maintheorem}
and \ref{thm:completedescription}. 
\medskip

{\sc Proof of Theorems \ref{thm:maintheorem} 
and \ref{thm:completedescription}}. 
%\medskip

By Propositions \ref{prop:degphiKY},   
            \ref{prop:intersectionofquadrics}, 
            \ref{prop:deg1case1exclusion}, 
       and  \ref{prop:prop:deg1case2exclusion}, 
we have $\deg \varPhi_{K_Y} = 2$.  
Thus Theorems \ref{thm:maintheorem} and \ref{thm:completedescription} 
follow from 
Propsitions \ref{prop:L^2}, 
            \ref{prop:case1exclusion}, 
            \ref{prop:case3-1exclusion}, 
            \ref{prop:case3-2exclusion},    
     and    \ref{prop:case2conclusion}.  \qed

\begin{remark} \label{rem:onthedescription}
Let $X_{(1)}$ and $X_{(2)}$ be two minimal complex surfaces 
as in Theorem \ref{thm:completedescription}, 
$\pi_{(i)} : Y_{(i)} \to X_{(i)}$ ($i =1$, $2$), the 
unramified double cover corresponding to the torsion group,   
$f_{(i)} : Y_{(i)} \to W_{(i)} = \varSigma_{d_{(i)}}$,  
the generically two-to-one morphism as 
in Theorem \ref{thm:completedescription},  
and $B_{(i)}$, the branch divisor of $f_{(i)}$. 
Then if $X_{(1)}$ and $X_{(2)}$ are isomorphic to each other, 
so are the triplets
$(W_{(1)}, \iota |_{W_{(1)}}, B_{(1)})$ 
and $(W_{(2)}, \iota |_{W_{(2)}}, B_{(2)})$,  
where $\iota |_{W_{(i)}}$ denotes the involution of $W_{(i)}$ 
corresponding to the generator of the Galoirs group of $\pi_{(i)}$. 
This is verified as follows. 
Let $p_{(i)} : \tilde{Y}_{(i)} \to Y_{(i)}$ be the shortest composite  
of quadric transformations such that the variable part of 
$p_{(i)}^* |K_{Y_{(i)}}|$ is free from base points, and 
$r^{\prime}_{(i)} : Z_{(i)}^{\prime} \to W_{(i)} = \varSigma_{d_{(i)}}$, 
the blowing-up at two $[3,3]$-points of the branch divisor $B_{(i)}$. 
Then $f_{(i)}$ induces a morphims 
$\tilde{f}_{(i)} : \tilde{Y}_{(i)} \to Z_{(i)}^{\prime}$. The 
projection $r_{(i)}^{\prime}$ is the blowing-down of the image 
by $\tilde{f}_{(i)}$ of the exceptional divisor of 
$p_{(i)} : \tilde{Y}_{(i)} \to Y_{(i)}$. Since $Z_{(i)}^{\prime}$ is 
the minimal desingularization of the canonical image 
of $Y_{(i)}$, we have the assertion. 
\end{remark}

\section{The moduli space for the case $\chi=4$} \label{scn:onmodulispace}

In this section, we shall study the moduli space for 
surfaces as in Theorem \ref{thm:completedescription},  
and give a proof for Theorem \ref{thm:moduli}.  
For this purpose, we shall first study 
the explicit description of our surfaces in more detail. 

Let $X$ be a minimal algebraic surface with $c_1^2= 2\chi -1$, 
$\chi = 4$, and $\mathrm{Tors} \simeq \mathbb{Z} /2$. 
We denote by $\pi : Y \to X$  the 
unramified double cover corresponding to the torsion group, 
and by $p: \tilde{Y} \to Y$, the shortest composite of 
quadric transformations such that the variable part of 
$p^*|K_Y|$ is free from base points.  
Then there exist an even integer $0 \leq d \leq 2$ and 
a generically two-to-one morphism $f: Y \to W = \varSigma_d$ 
satisfying the three conditions 
given in Theorem \ref{thm:completedescription}. 
In what follows, we denote by $\iota |_{W}$ the involution of $W$ 
corresponding to the generator of the Galois group $G = \mathrm{Gal} (Y/X)$. 

Let $r^{\prime} : Z^{\prime} \to W$ be the blowing-up at two 
$[3,3]$-points, which we shall call $w_1$ and $w_2$, 
of the branch divisor $B$ of $f$. 
Then $f_W= f \circ p : \tilde{Y} \to W$ 
lifts to a morphism $f^{\prime} : \tilde{Y} \to Z^{\prime}$.
We denote by $e_i = {r^{\prime}}^{-1} (w_i)$ the exceptional divisor 
of $r^{\prime}$ lying over $w_i$.  
Let $\tilde{r} : \tilde{Z} \to Z^{\prime}$ be the blowing-up at 
two quadraple points, 
which we shall call $w_1^{\prime} \in e_1$ and $w_2^{\prime} \in e_2$, 
of the branch divisor of $f^{\prime}$. 
Then $f^{\prime}$ 
lifts to a morphism $\tilde{f} : \tilde{Y} \to \tilde{Z}$. 
We denote by $e_i^{\prime} = \tilde{r}^{-1} (w_i^{\prime})$ the 
exceptional divisor of $\tilde{r}$ lying over $w_i^{\prime}$. 
Let us use the same symbol $e_i$ for the total transform to $\tilde{Z}$ 
of the divisor $e_i \subset Z^{\prime}$. 
Then there exists a reduced member 
$\tilde{B}_0 \in |(r^{\prime} \circ \tilde{r})^* (-4K_W) 
- 3 \sum e_i - 3 \sum e_i^{\prime}|$ satisfying 
$\tilde{B}_0 \cap \tilde{r}^{-1}_* (e_1) 
= \tilde{B}_0 \cap \tilde{r}^{-1}_* (e_2) = \emptyset$ 
such that the branch divisor of $\tilde{f}$ is given by 
$\tilde{B}_0 + \sum \tilde{r}^{-1}_* (e_i)$. 
Note that the divisor $\tilde{B}_0$ has at most negligible singularities. 
In what follows, $\varDelta_0$ and $\varGamma$ denote the 
minimal section and a fiber respectively of the Hirzebruch surface 
$W = \varSigma_d \to \mathbb{P}^1$. 

\begin{lemma} \label{lm:configw1w1prime}
Let $\iota |_{Z^{\prime}}$ be the involution of $Z^{\prime}$ induced 
by the involution $\iota |_W $ of $W$. Then the configuration of 
the four points $w_1$, $w_2 = \iota |_W (w_1)$, 
$w_1^{\prime}$, 
and $w_2^{\prime} =\iota |_{Z^{\prime}} (w_1^{\prime})$ 
satisfies the following three conditions$:$

i$)$ if $d=2$, then $w_1 \notin \varDelta_0$\
$;$

ii$)$ if the two points $w_1$ and $w_2$ lie on one and the same member 
of the linear system $|\varGamma|$, then for each $i =1$, $2$, 
the point $w_i^{\prime}$ does not lie on the strict transform 
to $Z^{\prime}$ of this member$;$ 

iii$)$ if $d$ equals $0$, and the two points $w_1$ and $w_2$ lie on 
one and the same member of the linear system $|\varDelta_0|$, 
then for each $i =1$, $2$, the point $w_i^{\prime}$ does not 
lie on the strict transform to $Z^{\prime}$ of this member. 
\end{lemma}
    
Proof. 
i). Assume that $d=2$ and $w_1 \in \varDelta_0$. Then since 
$\varDelta_0$ is stable under the action by $G$ on $W$, 
we have $w_2 \in \varDelta_0$. 
Thus ${r^{\prime}}^{-1}_* (\varDelta_0)$ is a $(-4)$-curve 
on $Z^{\prime}$, hence 
${r^{\prime}}^{-1}_* (\varDelta_0) (-K_{Z^{\prime}}) < 0$. 
It follows that ${r^{\prime}}^{-1}_* (\varDelta_0)$ is a 
fixed component of the linear systme $|-K_{Z^{\prime}}|$. 
This is impossible, since by the proof of our complete descripiton   
the pull-back ${f^{\prime}}^* |-K_{Z^{\prime}}|$ is the variable part of 
$|K_{\tilde{Y}}|$. Thus we have $w_1 \notin \varDelta_0$ for the 
case $d=2$. 

ii). Assume that $w_1$ and $w_2$ lie on one and the same member   
$\varGamma_0 \in |\varGamma|$. 
Then since $w_2 = \iota|_W (w_1)$, the member $\varGamma_0$ is 
stable under the action by $G$. It follwos that $\varGamma_0$ 
passes exactly two of the fixed points of the involution $\iota |_W$. 
Moreover if $w_1^{\prime} \in {r^{\prime}}^{-1}_* (\varGamma_0)$, then 
we obtain $w_2^{\prime} \in {r^{\prime}}^{-1}_* (\varGamma_0)$, 
$(r^{\prime} \circ \tilde{r})^{-1}_* (\varGamma_0) \sim 
(r^{\prime} \circ \tilde{r})^* (\varGamma) - \sum e_i - \sum e_i^{\prime}$, 
and $\tilde{B}_0 ((r^{\prime} \circ \tilde{r})^{-1}_* (\varGamma_0))
= -4 < 0$. The last inequality implies that $\varGamma_0$
is an irreducible component of the branch divisor $B$. 
This however is impossible, since, by the condition in 
Theorem \ref{thm:completedescription}, the branch divisor $B$ cannot pass 
any fixed points of the involution $\iota |_W$. 
Thus  we have $w_1^{\prime} \notin {r^{\prime}}^{-1}_* (\varGamma_0)$. 

iii).  By the same arguement as in the proof of ii), 
we can prove iii).  \qed

\begin{remark} \label{rem:gamma1gamma2}
As shown in the proof of Lemma above, if the two points 
$w_1$ and $w_2$ lie on one and the same member 
$\varGamma_0 \in | \varGamma |$, then this $\varGamma_0$ is stable 
under the action by $G$ on $W$.  There exist exactly two 
members of $|\varGamma|$ stable under the action by $G$.  
In what follows, we denote by $\varGamma_1$ and $\varGamma_2$ these 
two members. For each $i=1$, $2$, exactly two fixed points of the 
action by $G$ lie on $\varGamma_i$.  
\end{remark}
       
Next let us show that if conversely the configuration of 
four points $w_i$'s and $w_i^{\prime}$'s satisfies the three conditions 
in Lemma \ref{lm:configw1w1prime}, then the proceedure implied by 
our structure theorem in fact produces a minimal surface with 
the desired invariants. 
Some of the results below will be used later, in our proof of 
the uniqueness of the deformation type.  
Let $W = \varSigma_d$ be the Hirzebruch surface of degree $d=0$ or $2$, 
and $\iota |_W$, the involution (\ref{eq:involution}) given in 
Remark \ref{rem:involutiondescription}. 
Take a point $w_1 \in W$ outside the fixed locus of $\iota |_W$. 
We denote by $r^{\prime} : Z^{\prime} \to W$ the blowing-up at 
two points $w_1$ and $w_2 = \iota |_W (w_1)$, 
and by $e_i = {r^{\prime}}^{-1} (w_i)$, 
the exceptional curve lying over $w_i$.  
Let $\iota |_{Z^{\prime}}$ be the involution of $Z^{\prime}$ 
induced by $\iota |_W$. 
Take a point $w_1^{\prime} \in e_1 \subset Z^{\prime}$. 
We donote by $\tilde{r} : \tilde{Z} \to Z^{\prime}$ the 
blowing-up at two points $w_1^{\prime}$ and 
$w_2^{\prime} = \iota |_{Z^{\prime}} (w_1^{\prime})$, and 
by $e_i^{\prime} = \tilde{r}^{-1} (w_i^{\prime})$, 
the exceptional curve lying over $w_i^{\prime}$.  
We use the same symbol $e_i$ for the total transform to $\tilde{Z}$ 
of the divisor $e_i$ on $Z^{\prime}$. 
We assume that the configuration of $w_i$'s and $w_i^{\prime}$'s satisfies 
the three conditions i), ii), and iii) in Lemma \ref{lm:configw1w1prime}. 

Let $\varGamma_1$ and $\varGamma_2$ be two distinct members of 
$|\varGamma|$ stable under 
the natural action by $G = \langle \iota |_W\rangle$ 
on $W$ (see Remark \ref{rem:gamma1gamma2}). 
We take the minimal section $\varDelta_0$ 
and an irreducible member 
$\varDelta_{\infty} \in |\varDelta_0 + d \varGamma |$ such that 
both are stable under the action by $G$, 
and $\varDelta_0 \cap \varDelta_{\infty} = \emptyset$ holds.      
Let $m$ be a positive integer. 
Since the divisor $m(\varDelta_0 + \frac{d+2}{2} \varGamma_1) $ 
is stable under the action by $G$, we obtain a natural action on 
$H^0 (\mathcal{O}_W (m(\varDelta_0 + \frac{d+2}{2} \varGamma) ))$ by  
identifying this space with that of meromorphic functions 
with poles at most $m(\varDelta_0 + \frac{d+2}{2} \varGamma_1 )$.

We put $ \varLambda_m = |m(\varDelta_0 + \frac{d+2}{2} \varGamma)|$, 
and denote by $\varLambda_m^+$ and $\varLambda_m^-$ the subsystems 
of $\varLambda_m$ corresponding 
to the eigenspaces of eigenvalues $+1$ and $-1$ repectively 
with respect to ${\iota |_W}^*$.   
Moreover, for an effective divisor $C$ on $\tilde{Z}$, we  put 
\begin{align}
\varLambda_m (C) &= \{ D \in \varLambda_m ; 
(r^{\prime} \circ \tilde{r})^*(D) -C \succeq 0 \} \quad
&\tilde{\varLambda}_m (C) 
= (r^{\prime} \circ \tilde{r})^* \varLambda_m (C) - C \notag \\
\varLambda_m^+ (C) &= \{ D \in \varLambda_m^+ ; 
(r^{\prime} \circ \tilde{r})^*(D) -C \succeq 0 \} \quad
&\tilde{\varLambda}_m^+ (C) 
= (r^{\prime} \circ \tilde{r})^* \varLambda_m^+ (C) - C \notag \\
\varLambda_m^- (C) &= \{ D \in \varLambda_m^- ; 
(r^{\prime} \circ \tilde{r})^*(D) -C \succeq 0 \} \quad 
&\tilde{\varLambda}_m^- (C) 
= (r^{\prime} \circ \tilde{r})^* \varLambda_m^- (C) - C ,  \notag
\end{align}
where the symbol $\succeq 0$ means effectiveness of a divisor. 
We abbreviate $\tilde{\varLambda}_m (0)$, 
$\tilde{\varLambda}_m^+ (0)$, 
and $\tilde{\varLambda}_m^- (0)$  
to  
$\tilde{\varLambda}_m$, 
$\tilde{\varLambda}_m^+$, 
and $\tilde{\varLambda}_m^-$ respectively. 
Note that if $\tilde{f} : \tilde{Y} \to \tilde{Z}$ is 
the generically two-to-one morphism obtained 
as in the begining of this section from our structure theorem, 
then we have 
$\tilde{B}_0 \in \tilde{\varLambda}_8^+ (
3 \sum e_i + 3 \sum e_i^{\prime} )$, where  
$\tilde{B}_0 + \sum \tilde{r}^{-1}_* (e_i)$ gives  
the branch divisor of $\tilde{f} : \tilde{Y} \to \tilde{Z}$.   

\begin{lemma} \label{lm:baseptfreesystems}
%Let the configuration of $w_i$'s and $w_i^{\prime}$'s satisfy 
%the three conditions i$)$, ii$)$, and iii$)$ as in 
%Lemma \ref{lm:configw1w1prime}. 
%Then we have the following$:$  

1$)$ The linear system $\tilde{\varLambda}_2^+ $ has no base point. 

2$)$ The linear system 
$\tilde{\varLambda}_2^+ (\sum e_i + \sum e_i^{\prime} )$ 
has no base point. 

3$)$ The linear system 
$\tilde{\varLambda}_8^+ (3\sum e_i + 3\sum e_i^{\prime} )$ 
has no base point. 

4$)$ The linear systems $|-K_{Z^{\prime}}|$ and $|-K_{\tilde{Z}}|$  
have no base point.  
%
%Five linear systems 
%$\tilde{\varLambda}_2^+ $,
%$\tilde{\varLambda}_2^+ (\sum e_i + \sum e_i^{\prime} )$,
%$\tilde{\varLambda}_8^+ (3\sum e_i + 3\sum e_i^{\prime} )$,  
%$|-K_{Z^{\prime}}|$, 
%and $|-K_{\tilde{Z}}|$ have not base points. 
\end{lemma}

Proof. 
Since we have 
$
 \tilde{\varLambda}_2^+ 
+  3 \tilde{\varLambda}_2^+ (\sum e_i + \sum e_i^{\prime} )
\subset 
 \tilde{\varLambda}_8^+ (3\sum e_i + 3\sum e_i^{\prime} ) 
$, 
the assertion 3) follows from the assertions 1) and 2). 

Assume that we have the assertions 1) and 2). 
Then by  
$
\tilde{\varLambda}_2^+ (\sum e_i + \sum e_i^{\prime} )
\subset  
|- K_{\tilde{Z}}| 
$, 
we see that 
the linear system $|- K_{\tilde{Z}}|$ has no base point. 
Moreover, by this together with the Riemann--Roch theorem 
and the vanishing theorem, we obtain 
$h^0 ( \mathcal{O}_{\tilde{Z}} (-K_{\tilde{Z}})) 
= \chi ( \mathcal{O}_{\tilde{Z}} ) + K_{\tilde{Z}}^2 = 5$. 
Meanwhile,   
since  $r^{\prime} : Z^{\prime} \to W$ is the blowing-up at 
two points $w_i^{\prime}$'s, 
we have 
$h^0 ( \mathcal{O}_{Z^{\prime}} (-K_{Z^{\prime}}))
\geq h^0 ( \mathcal{O}_W (-K_{W}))  -2 = 7$.  
Thus we obtain $
h^0 ( \mathcal{O}_{Z^{\prime}} (-K_{Z^{\prime}})) 
- h^0 ( \mathcal{O}_{\tilde{Z}} (-K_{\tilde{Z}})) 
\geq 2$, which implies that neither of the two points $w_i^{\prime}$'s 
is a base point of $|-K_{Z^{\prime}}|$. 
From this, we infer that $|-K_{Z^{\prime}}|$ has no base point. 
So the assertion 4) also follows from the assertions 1) and 2). 

Thus we only need to show the assertions 1) and 2). 
First, let us show the assertion 1). 
Let $C_0$ be a general member of $| \varGamma |$. 
Then since the divisor 
$2\varDelta_0 + \frac{d+2}{2} (C_0 + \iota |_W (C_0)) \in \varLambda_2$ 
is stable under the action by $G$, and the divisor 
\[
 (2\varDelta_0 + \frac{d+2}{2} (C_0 + \iota |_W (C_0))) 
 - 2(\varDelta_0 + \frac{d+2}{2} \varGamma_1)
\] 
has no support at $\varDelta_{\infty} \cap \varGamma_2$, 
the divisor  
$2\varDelta_0 + \frac{d+2}{2} (C_0 + \iota |_W (C_0))$ 
is a member of $\varLambda_2^+$.   
Thus the base locus of $\varLambda_2^+$ is contained in  
$\varDelta_0$. Using a similar argument, we can show that   
$2\varDelta_{\infty} + \frac{2- d}{2} (C_0 + \iota |_W (C_0))
\in \varLambda_2^+$, 
so that the base locus of $\varLambda_2^+$ is contained  
in $\varDelta_{\infty}$. 
Thus since $\varDelta_0 \cap \varDelta_{\infty} = \emptyset$, 
the linear system $\tilde{\varLambda}_2^+$ has 
no base point. Hence we have the assertion 1). 

Next, let us show the assertion 2). 
We shall show it by dividing our situation into several cases. 
In what follows, for each $i = 1$, $2$, we denote by $\varGamma_{(i)}$
the unique member of $| \varGamma |$ passing $w_i$. 

Case $1$-$1$: the case where $d=0$ holds, and 
the two points $w_1$ and $w_2$ lie neither on one and the same member  
of $| \varGamma |$ nor on that of $|\varDelta_0 |$. 
In this case, for each $i=1$, $2$, we denote by $\varDelta_{(i)}$
the unique member of $|\varDelta_0|$ passing $w_i$. 
This case is divided into two subcases: case $1$-$1$-$1$ and 
case $1$-$1$-$2$. 

Case $1$-$1$-$1$; the subcase of case $1$-$1$ where 
$w_1^{\prime} \notin {r^{\prime}}^{-1}_* (\varGamma_{(1)})$ and 
$w_1^{\prime} \notin {r^{\prime}}^{-1}_* (\varDelta_{(1)})$. 
In this case, take global coordinates $(s_1^{\prime}, \xi_1^{\prime})$ of 
$W \setminus (\varGamma_{(2)} \cup \varDelta_{(2)}) \simeq \mathbb{A}^2$ 
such that 
$\varGamma_{(1)}$ is given by $s_1^{\prime}  =0$, 
$\varGamma_{(2)}$ by $s_1^{\prime}  =\infty$, 
$\varDelta_{(1)}$ by $\xi_1^{\prime}  =0$, 
and $\varDelta_{(2)}$ by $\xi_1^{\prime}  =\infty$.   
Then the involution $\iota |_W$ is given by 
$(s_1^{\prime}, \xi_1^{\prime}) 
\mapsto (1 / s_1^{\prime} , 1 / \xi_1^{\prime} )$, 
and the linear system $\varLambda_2^+$ is spanned by 
the five elements 
${s_1^{\prime}}^l {\xi_1^{\prime}}^m 
+ { s_1^{\prime}}^{2-l} {\xi_1^{\prime}}^{2-m}$ 
($0 \leq l \leq2 $, $ 0 \leq m \leq2$). 
Thus the linear system $\varLambda_2^+ (\sum e_i + \sum e_i^{\prime}) $ 
is spanned by the three elements
\[
  a_0 (s_1^{\prime} + s_1^{\prime} {\xi_1^{\prime}}^2)
+ b_0 (\xi_1^{\prime} + { s_1^{\prime}}^2 \xi_1^{\prime} ),
\qquad  \quad 
s_1^{\prime} \xi_1^{\prime},
\qquad  \quad 
{s_1^{\prime}}^2  +  {\xi_1^{\prime}}^2 , 
\]
where $a_0 \neq 0$ and $b_0 \neq 0$ are certain non-zero complex numbers. 
From this, we infer that 
the set $\{ w_1 , w_2 \}$ forms the base locus 
of $\varLambda_2^+ (\sum e_i + \sum e_i^{\prime}) $, 
and that any general member of this linear system is smooth. 
By this together with 
$\sum \varGamma_{(i)} + \sum \varDelta_{(i)} \in 
\varLambda_2^+ (\sum e_i + \sum e_i^{\prime})$, 
we see that the linear system 
$\tilde{\varLambda}_2^+ (\sum e_i + \sum e_i^{\prime})$ has no base point. 

Case $1$-$1$-$2$: the subcase of case $1$-$1$ where 
$w_1^{\prime} \in {r^{\prime}}^{-1}_* (\varGamma_{(1)})$ or 
$w_1^{\prime} \in {r^{\prime}}^{-1}_* (\varDelta_{(1)})$. 
Since the proof is the same, we only give a proof for the 
case $w_1^{\prime} \in {r^{\prime}}^{-1}_* (\varGamma_{(1)})$. 
Assume that $w_1^{\prime} \in {r^{\prime}}^{-1}_* (\varGamma_{(1)})$. 
Since we have 
$C_0 + \iota |_W  (C_0) + \sum \varGamma_{(i)} \in 
\varLambda_2^+ (\sum e_i + \sum e_i^{\prime})$ 
for any general member $C_0$ of $| \varDelta_0 |$, 
the base locus of $\tilde{\varLambda}_2^+ (\sum e_i + \sum e_i^{\prime})$ is 
contained in $\sum (r^{\prime} \circ \tilde{r})^{-1}_* (\varGamma_{(i)})$.   
Meanwhile, since we have 
$C_1 + \iota |_W (C_1) \in \varLambda_2^+ (\sum e_i + \sum e_i^{\prime})$ 
for any general member $C_1$ of $\varLambda_1 (\sum e_i)$, 
the base locus of $\tilde{\varLambda}_2^+ (\sum e_i + \sum e_i^{\prime})$ is 
contained in $\sum \tilde{r}^{-1}_* (e_i)$. 
Since we have 
$(\sum (r^{\prime} \circ \tilde{r})^{-1}_* (\varGamma_{(i)}))
\cap (\sum \tilde{r}^{-1}_* (e_i) ) = \emptyset $,
we see that the linear system 
$\tilde{\varLambda}_2^+ (\sum e_i + \sum e_i^{\prime})$ 
has no base point. 

Case $1$-$2$: the case where $d=0$ holds, and the two points 
$w_1$ and $w_2$ lie on one and the same member of 
$| \varGamma |$ or $| \varDelta_0 |$. 
In this case, for each $i =1$, $2$, we denote by 
$\varDelta_{(i)}$ the unique member of $| \varDelta_0 |$ 
passing $w_i$. 
By exchanging $\varDelta_0$ and $\varGamma$ if necessary, we 
may assume that the two points $w_1$ and $w_2$ lie on one  
and the member $\varGamma_0 \in | \varGamma |$.  
Moreover, by Remark \ref{rem:gamma1gamma2}, 
by exchanging $\varGamma_1$ and $\varGamma_2$ 
if necessary, we may assume that 
$\varGamma_0 = \varGamma_{(1)} = \varGamma_{(2)}  = \varGamma_1$. 
Then this case is divided into two subcases:  
case $1$-$2$-$1$ and case $1$-$2$-$2$. 

Case $1$-$2$-$1$: the subcase of case $1$-$2$ where 
$w_1^{\prime} \notin {r^{\prime}}^{-1}_* ( \varDelta_{(1)})$. 
Note that we have assumed the condition ii) of 
Lemma \ref{lm:configw1w1prime} for our configuration,  
so that we have $w_1^{\prime}$, 
$w_2^{\prime} \notin  {r^{\prime}}^{-1}_* (\varGamma_1)$. 
For any general member $C_0 \in | \varDelta_0 |$, we have 
$C_0 + \iota |_W (C_0) + 2 \varGamma_1 
\in \varLambda_2^+ (\sum e_i + \sum e_i^{\prime})$. 
Thus the base locus of 
$\tilde{\varLambda}_2^+ (\sum e_i + \sum e_i^{\prime})$ 
is contained in $2 (r^{\prime} \circ \tilde{r})^{-1}_* (\varGamma_1) + 
\sum \tilde{r}^{-1}_* (e_i)$.  
Meanwhile for any general $C_1$ ($ \neq \varDelta_{(1)} + \varGamma_1$) 
$\in \varLambda_1 (e_1 + e_1^{\prime})$, we have  
$C_1 + \iota |_W (C_1) \in 
\varLambda_2^+ (\sum e_i + \sum e_i^{\prime})$ and   
the irreducibility and smoothness at $w_1$ of $C_1$.   
Thus the base locus of 
$\tilde{\varLambda}_2^+ (\sum e_i + \sum e_i^{\prime})$ 
is contained in 
$(r^{\prime} \circ \tilde{r})^{-1}_* (C_1 + \iota |_W (C_1))$. 
Since 
$2 (r^{\prime} \circ \tilde{r})^{-1}_* (\varGamma_1) + 
\sum \tilde{r}^{-1}_* (e_i)$ and  
$(r^{\prime} \circ \tilde{r})^{-1}_* (C_1 + \iota |_W (C_1))$ 
do not intersect each other,  
we see that the linear system
$\tilde{\varLambda}_2^+ (\sum e_i + \sum e_i^{\prime})$ has 
no base point. 
  
Case $1$-$2$-$2$: the subcase of case $1$-$2$ where 
$w_1 \in {r^{\prime}}^{-1}_* (\varDelta_{(1)})$.   
By the same argument as one given in the proof for case $1$-$2$-$1$, 
we see that the base locus of 
$\tilde{\varLambda}_2^+ (\sum e_i + \sum e_i^{\prime})$ 
is contained in $2 (r^{\prime} \circ \tilde{r})^{-1}_* (\varGamma_1) + 
\sum \tilde{r}^{-1}_* (e_i)$.  
Meanwhile, for any general member $C_1 \in | \varGamma |$, we have 
$C_1 + \iota |_W (C_1) + \sum \varDelta_{(i)} 
\in \varLambda_2^+ ( \sum e_i + \sum e_i^{\prime})$. 
Thus the base locus of 
$\tilde{\varLambda}_2^+ (\sum e_i + \sum e_i^{\prime})$ 
is contained in 
$\sum (r^{\prime} \circ \tilde{r})^{-1}_* (\varDelta_{(i)})$. 
Since $2 (r^{\prime} \circ \tilde{r})^{-1}_* (\varGamma_1) + 
\sum \tilde{r}^{-1}_* (e_i)$ and 
$\sum (r^{\prime} \circ \tilde{r})^{-1}_* (\varDelta_{(i)})$ 
do not intersect each other, we see that 
the linear system 
$\tilde{\varLambda}_2^+ (\sum e_i + \sum e_i^{\prime})$
has no base point.  

Case $2$-$1$: the case where $d=2$ holds, and 
the two points $w_1$ and $w_2$ do not lie on one and the same 
member of $| \varGamma |$.  
Note that since we have assumed 
the condition i) in Lemma \ref{lm:configw1w1prime}, we have 
$w_1 \notin \varDelta_0$. 
Note also that for this case, or more generally for case $d=2$,  
we have $\dim \varLambda_1^- =1$, and any general member 
of this linear system is an irreducible curve stable under the action by $G$ 
that passes two points $\varDelta_{\infty} \cap \varGamma_1$ and 
$\varDelta_{\infty} \cap \varGamma_2$.    
We denote by $\varDelta_1$ the unique member of $\varLambda_1^-$ 
that passes the two points $w_1$ and $w_2$. 
Then this case is divided into three subcases:  
case $2$-$1$-$1$, case $2$-$1$-$2$, and 
case $2$-$1$-$3$. 

Case $2$-$1$-$1$: 
the subcase of case $2$-$1$ where 
$w_1^{\prime} \notin {r^{\prime}}^{-1}_* (\varGamma_{(1)})$ 
and $w_1^{\prime} \notin {r^{\prime}}^{-1}_* (\varDelta_{1})$.  
Since the divisor $\varDelta_0 + 2 \varGamma_{(1)}$ is the unique 
reducible member of $\varLambda_1 (e_1 + e_1^{\prime})$, 
and we have $h^0 (\mathcal{O}_W (\varDelta_0 + 2 \varGamma)) =4$, 
any general member of $\varLambda_1 (e_1 + e_1^{\prime})$ is 
irreducible and non-singular. 
By this together with $\varDelta_0 + 2 \varGamma_{(1)} \in 
\varLambda_1 (e_1 + e_1^{\prime})$, 
we see that $\tilde{\varLambda}_1 (e_1 + e_1^{\prime})$ has 
no base point, 
and $(r^{\prime} \circ \tilde{r})^{-1}_* (C_0 + \iota |_W (C_0))
\in \tilde{\varLambda}_2^+ (\sum e_i + \sum e_i^{\prime})$    
for any general member $C_0 \in \varLambda_1 (e_1 + e_1^{\prime})$. 
Thus we see that $\tilde{\varLambda}_2^+ (\sum e_i + \sum e_i^{\prime})$ 
has no base point. 

Case $2$-$1$-$2$: 
the subcase of case $2$-$1$ where 
$w_1^{\prime} \in {r^{\prime}}^{-1}_* (\varGamma_{(1)})$. 
Since we have $2 \varDelta_0 + \sum \varGamma_{(i)} 
+ C_0 + \iota |_W (C_0) 
\in \varLambda_2^+ (\sum e_i + \sum e_i^{\prime})$
for any general member $C_0 \in | \varGamma |$, 
the base locus of 
$\tilde{\varLambda}_2^+ (\sum e_i + \sum e_i^{\prime})$ 
is contained in 
$2 (r^{\prime} \circ \tilde{r})^* (\varDelta_0) 
+ \sum (r^{\prime} \circ \tilde{r})^{-1}_* (\varGamma_{(i)})$.   
By this together with 
$2 \varDelta_1 \in \varLambda_2^+ (\sum e_i + \sum e_i^{\prime})$, 
we see that the linear system 
$\tilde{\varLambda}_2^+ (\sum e_i + \sum e_i^{\prime})$ has 
no base point. 

Case $2$-$1$-$3$: the subcase of case $2$-$1$ where 
$w_1^{\prime} \in {r^{\prime}}^{-1}_* (\varDelta_1)$. 
Since we have 
$C_0 + \varDelta_1 \in \varLambda_2^+ (\sum e_i + \sum e_i^{\prime})$
for any general member $C_0 \in \varLambda_1^-$, 
the base locus of $\tilde{\varLambda}_2^+ (\sum e_i + \sum e_i^{\prime})$
is contained in $(r^{\prime} \circ \tilde{r})^{-1}_* (\varDelta_1)$. 
By this together with 
$2 (\varDelta_0 + \sum \varGamma_{(i)}) 
\in \varLambda_2^+ (\sum e_i + \sum e_i^{\prime})$, 
we see that the linear system 
$\tilde{\varLambda}_2^+ (\sum e_i + \sum e_i^{\prime})$ has 
no base point. 

Case $2$-$2$: 
the case where $d=2$ holds, and the two points $w_1$ and $w_2$ 
lie on one and the same member of $| \varGamma |$. 
By Remark \ref{rem:gamma1gamma2}, we may 
assume $w_1$, $w_2 \in \varGamma_1$.  
Note that we have assumed the conditions i) and ii) 
of Lemma \ref{lm:configw1w1prime}
for our configuration, so that we have 
$w_1 \notin \varDelta_0$ and 
$w_1^{\prime} \notin {r^{\prime}}^{-1}_* (\varGamma_1)$. 
Since we have 
$2 \varDelta_0 + 2 \varGamma_1 + C_0 + \iota |_W (C_0) 
\in \varLambda_2^+ (\sum e_i + \sum e_i^{\prime})$  
for any general member $C_0 \in | \varGamma |$,  
the base locus of 
$\tilde{\varLambda}_2^+ (\sum e_i + \sum e_i^{\prime})$ 
is contained in 
$2 (r^{\prime} \circ \tilde{r})^* (\varDelta_0) 
+ 2 (r^{\prime} \circ \tilde{r})^{-1}_* (\varGamma_1)
+ \sum \tilde{r}^{-1}_* (e_i)$. 
Meanwhile since we have 
$h^0 (\mathcal{O}_W (\varDelta_0 + 2 \varGamma)) = 4$, we have 
$C_1 + \iota |_{W} (C_1) \in \varLambda_2^+ (\sum e_i + \sum e_i^{\prime})$
for an irreducible member 
$C_1 \in \varLambda_1 (e_1 + e_1^{\prime})$.
Since 
$2 (r^{\prime} \circ \tilde{r})^* (\varDelta_0) 
+ 2 (r^{\prime} \circ \tilde{r})^{-1}_* (\varGamma_1)
+ \sum \tilde{r}^{-1}_* (e_i)$ 
and 
$(r^{\prime} \circ \tilde{r})^{-1}_* 
(C_1 + \iota |_{W} (C_1))$ 
do not intersect each other, we see that the linear system 
$\tilde{\varLambda}_2^+ (\sum e_i + \sum e_i^{\prime})$ has 
no base point. 

Now that we have shown the absence of base points 
of $\tilde{\varLambda}_2^+ (\sum e_i + \sum e_i^{\prime})$  
for all the eight cases $1$-$1$-$1$, \ldots, $2$-$2$, 
we have the assertion 2).  \qed  

%Lemma \ref{lm:baseptfreesystems} shows  
%smoothness of  general members of 
%$\tilde{\varLambda}_8^+ (3 \sum e_i + 3 \sum e_i^{\prime})$. 
Let $\tilde{B}_0$ be a reduced member of 
$\tilde{\varLambda}_8^+ (3 \sum e_i + 3 \sum e_i^{\prime})$ 
that has at most negligible singularities, 
satisfies $\tilde{B}_0 \cap \sum \tilde{r}^{-1}_* (e_i) = \emptyset $, 
and passes no fixed point of the action by $G$ on $\tilde{Z}$. 
Existence of such $\tilde{B}_0$ is ensured 
by Lemma \ref{lm:baseptfreesystems}.   
%Let $\tilde{f} : \tilde{Y} \to \tilde{Z}$ be the 
%canonical resolution of the double cover branched along 
%$\tilde{B}_0 + \sum \tilde{r}^{-1}_* (e_i)$. 
%
Let $\tilde{Y}$ be the 
canonical resolution of the double cover of $\tilde{Z}$
branched along 
$\tilde{B}_0 + \sum \tilde{r}^{-1}_* (e_i)$, and 
$\tilde{f} : \tilde{Y} \to \tilde{Z}$, the natural projection.  
We have $ \tilde{f}^* ({\tilde{r}}^{-1}_* (e_i)) = 2 E_i$ 
for a $(-1)$-curve $E_i$ on $\tilde{Y}$. 
Let $p : \tilde{Y} \to Y$ be the blowing-down 
of $E_1$ and $E_2$. 
Then we have $|K_{\tilde{Y}}| = 
(\tilde{r} \circ \tilde{f})^* |-K_{Z^{\prime}}| + 2 \sum E_i$. 
Since $|-K_{Z^{\prime}}|$ has no base point 
by Lemma \ref{lm:baseptfreesystems}, 
we see that $Y$ is a minimal surface with 
$c_1^ 2 = 14$ and $\chi =8$. 
By \cite[Lemma 3.1]{ngsffermat} and 
\[
 (r^{\prime} \circ \tilde{r}) (\tilde{B}_0) \varGamma_1  \equiv
 (r^{\prime} \circ \tilde{r}) (\tilde{B}_0) \varGamma_2  \equiv
 (r^{\prime} \circ \tilde{r}) (\tilde{B}_0) \varDelta_0  \equiv
  0   \mod 4, 
\]
there exists a unique free lifting to $\tilde{Y}$ of the action by $G$ on $W$. 
Let $X = Y/G$ be the quotient of $Y$ by the induced free action by $G$ on $Y$. 
Then by \cite[Theorem 1]{bound'''} or \cite[(ii) in Theorem A]{on2chi-2}, 
the surface $X$ is a minimal surface with $c_1 ^2 = 2 \chi -1$, $\chi =4$, 
and $\mathrm{Tors} (X) \simeq \mathbb{Z} /2$. 
Thus we have the following:  

\begin{proposition}  \label{prop:sufficientprocedure}
Let $W = \varSigma_d$ be the Hirzebruch surface of degree $d=0$ or $2$.
Let $r^{\prime} : Z^{\prime} \to W$ be the blowing-up at 
two points $w_1$ and $w_2 = \iota |_W (w_1)$, 
where $\iota |_W$ is the involution of $W$ 
given in Remark \ref{rem:involutiondescription}, 
and $w_1$, a point outside the fixed locus of $\iota |_W$. 
Let $ \tilde{r} : \tilde{Z} \to Z^{\prime}$ be 
the blowing-up at 
two points $w_1^{\prime}$ and $w_2^{\prime} 
= \iota |_{Z^{\prime}} (w_1^{\prime})$, 
where $\iota |_{Z^{\prime}}$ is the induced involution of $Z^{\prime}$, 
and $w_1^{\prime}$, a point infinitely near to $w_1$. 
Put $e_i = {r^{\prime}}^{-1} (w_i)$ and 
$e_i^{\prime} = \tilde{r}^{-1}(w_i^{\prime})$ for 
each $i=1$, $2$, and assume that the configuration of 
$w_i$'s and $w_i^{\prime}$'s satisfies all the three conditions 
in Lemma \ref{lm:configw1w1prime}.  
Let $\tilde{B}_0$ be a reduced member 
of $\tilde{\varLambda}_8^+ ( 3 \sum e_i + 3 \sum e_i^{\prime})$ 
that has at most negligible singularities, satisfies 
$\tilde{B}_0 \cap \sum \tilde{r}^{-1}_* (e_i) = \emptyset$, 
and passes no fixed point of the induced action on $\tilde{Z}$ by 
$G = \langle \iota |_W \rangle$. 
%Let  $\tilde{f} :  \tilde{Y} \to \tilde{Z}$ denote the 
%canonical resolution of the double cover branched along 
%$\tilde{B}_0 + \sum \tilde{r}^{-1}_* (e_i)$, 
%and $p : \tilde{Y} \to Y$, the blowing-down of 
%two $(-1)$-curves $E_1 = {\tilde{f}}^{-1} ({\tilde{r}}^{-1}_* (e_1))$
%and $E_2 = {\tilde{f}}^{-1} ({\tilde{r}}^{-1}_* (e_2))$. 
%
%
Let $\tilde{Y}$ be the canonical resolution of 
the double cover of $\tilde{Z}$ branched along 
$\tilde{B}_0 + \sum \tilde{r}^{-1}_* (e_i)$, and 
$\tilde{f} :  \tilde{Y} \to \tilde{Z}$, the natural projection. 
Let $p : \tilde{Y} \to Y$ be the blowing-down of 
two $(-1)$-curves $E_1 = {\tilde{f}}^{-1} ({\tilde{r}}^{-1}_* (e_1))$
and $E_2 = {\tilde{f}}^{-1} ({\tilde{r}}^{-1}_* (e_2))$.
Then there exists a unique free lifting to $\tilde{Y}$ of the 
action by $G$ on $\tilde{Z}$, and the quotient $Y/G$ of $Y$  
by the induced free action is a minimal surface with 
$c_1^2 = 2 \chi -1$, $\chi =4$, and $\mathrm{Tors} \simeq \mathbb{Z} / 2$. 
  
\end{proposition}

Our Theorem \ref{thm:completedescription} together with 
Remark \ref{rem:involutiondescription} 
and Lemma \ref{lm:configw1w1prime} says that all minimal surfaces with 
$c_1^2 = 2 \chi -1$, $\chi =4$, and $\mathrm{Tors} \simeq \mathbb{Z} / 2$ 
are obtained by the procedure as in the proposition above. 
We use the following lemma in order to show the uniqueness of the 
deformation type. 

\begin{lemma} \label{lm:b0tildenonsing}
Let $r^{\prime} : Z^{\prime} \to W$ and 
$\tilde{r} : \tilde{Z} \to Z^{\prime}$, 
$w_i \in W$ and $w_i^{\prime} \in Z^{\prime}$ for $i=1$, $2$, 
and $e_i = {r^{\prime}}^{-1} (w_i)$ and 
$e_i^{\prime} = \tilde{r}^{-1} (w_i^{\prime})$ for $i =1$, $2$
be the morphisms, points, and divisors respectively as in 
Proposition \ref{prop:sufficientprocedure}. 
Then any general member $\tilde{B}_0$ of 
$\tilde{\varLambda}_8^+ ( 3 \sum e_i + 3 \sum e_i^{\prime})$
is non-singular and reduced. 
Further $h^i ( \mathcal{O}_{\tilde{Z}} (\tilde{B}_0)) =0$ 
holds for any positive integer $i > 0$. 
\end{lemma}

Proof. 
The first assertion follows from 3) in Lemma \ref{lm:baseptfreesystems}. 
The second assertion follows from 3) and 4) 
in Lemma \ref{lm:baseptfreesystems} and the vanishing theorem. 
\qed

Now let us show the uniqueness of the deformation type and 
the unirationality of the moduli space. 
For this purpose, we shall give another description of our 
surface $X$.  
Let $r^{\prime} : Z^{\prime} \to W$, $\tilde{r} : \tilde{Z} \to Z^{\prime}$, 
$w_i$, $w_i^{\prime}$, $e_i$, and $e_i^{\prime}$ be as in 
Proposition \ref{prop:sufficientprocedure}. 
Let $\varGamma_1$ and $\varGamma_2$ be as 
in Remark \ref{rem:gamma1gamma2}. 
We take the minimal section $\varDelta_0$ and an irreducible member 
$\varDelta_{\infty} \in |\varDelta_0 + d \varGamma|$ 
satisfying $\varDelta_0 \cap \varDelta_{\infty} = \emptyset$ 
such that both are stable under the action by $G$. 
Note that if $d = 2$, such $\varDelta_{\infty}$'s
form a one-dimensional family. 

The fixed locus of the action by $G$ on $\tilde{Z}$ is  
a set of four isolated points: 
$\{ (r^{\prime} \circ \tilde{r})^{-1} 
(\varGamma_i \cap \varDelta_j)\}_{i=1,2, \ j= 0, \infty}$. 
Let $\bar{r} : \bar{Z} \to \tilde{Z}$ be the 
blowing-up at these four points. 
For $i =1$, $2$ and $j=0$, $\infty$, we define the 
divisors $J_{i j}$ on $\bar{Z}$ as follows: \medskip

if $d = 0$, then   
$J_{i j} =  (r^{\prime} \circ \tilde{r} \circ \bar{r})^{-1} 
(\varGamma_i \cap \varDelta_j)$ for any 
$i=1$, $2$ and $j=0$, $\infty$; 

if $d = 2$, then  
$J_{1 0} =  (r^{\prime} \circ \tilde{r} \circ \bar{r})^{-1} 
(\varGamma_1 \cap \varDelta_{\infty})$, 
$J_{1 \infty} =  (r^{\prime} \circ \tilde{r} \circ \bar{r})^{-1} 
(\varGamma_1 \cap \varDelta_0)$, 
and 
$J_{2 j} =  (r^{\prime} \circ \tilde{r} \circ \bar{r})^{-1} 
(\varGamma_2 \cap \varDelta_j)$ 
for any $j = 0$, $\infty$.
\medskip

\noindent
Moreover for $i=1$, $2$ and $j= 0$, $\infty$, we define the 
divisors $\bar{\varGamma}_i$, $\bar{\varDelta}_j$, 
$\bar{e}_i$, and ${\bar{e}}_i^{\prime}$ on $\bar{Z}$ by 
$\bar{\varGamma}_i  =
(r^{\prime} \circ \tilde{r} \circ \bar{r})^{-1}_* (\varGamma_i)$, 
$\bar{\varDelta}_j  =
(r^{\prime} \circ \tilde{r} \circ \bar{r})^{-1}_* (\varDelta_j)$, 
$\bar{e}_i = 
\bar{r}^* ( \tilde{r}^{-1}_* (e_i))$, and 
${\bar{e}}_i^{\prime} = 
\bar{r}^* (e_i^{\prime})$. 
The four divisors $J_{i j}$'s form the 
set of all irreducible exceptional curves 
of $\bar{r} : \bar{Z} \to \tilde{Z}$.  
The action by $G$ on $\tilde{Z}$ lifts to one on $\bar{Z}$. 
Note that $\sum J_{i j}$ gives the fixed locus of 
the induced action by $G$ on $\bar{Z}$.   

Now let $\bar{V} = \bar{Z} / G$ be the quotient of $\bar{Z}$ 
by the induced action by $G$, and 
$\bar{\wp} : \bar{Z} \to \bar{V}$, the natural projection.  
Then $\bar{V}$ is smooth, and 
$\sum J_{i j}$ gives the ramification divisor of 
$\bar{\wp} : \bar{Z} \to \bar{V}$.  
For $i=1$, $2$ and $j=0$, $\infty$, we define 
the divisors $\bar{I}_{i j}$, $\bar{G}_i$, and $\bar{D}_j$ 
on $\bar{V}$ by 
$\bar{I}_{i j} = \bar{\wp} (J_{i j})$, 
$\bar{G}_i = \bar{\wp} (\bar{\varGamma}_i)$, and 
$\bar{D}_j = \bar{\wp} (\bar{\varDelta}_j)$. 
Moreover we define the divisors 
$\bar{\lambda}$ and ${\bar{\lambda}}^{\prime}$ on $\bar{V}$ by 
$\bar{\lambda} = \bar{\wp} (\bar{e}_1) = \bar{\wp} (\bar{e}_2)$ and 
${\bar{\lambda}}^{\prime} = \bar{\wp} (\bar{e}_1^{\prime}) 
= \bar{\wp} (\bar{e}_2^{\prime})$. 
The divisors $\bar{I}_{i j}$'s are 
non-singular rational curves with selfintersection  
${\bar{I}_{i j}}^2 = -2$. 
Note that $\sum \bar{I}_{i j}$ gives the branch divisor 
of $\bar{\wp} : \bar{Z} \to \bar{V}$. 

Let $\bar{\nu} : \bar{V} \to \tilde{V}$ be the blowing-down of 
the $(-1)$-curve $\bar{\lambda}^{\prime}$. 
For $i =1$, $2$ and $j=0$, $\infty$, we define the divisor $\tilde{I}_{i j}$ 
on $\tilde{V}$ by $\tilde{I}_{i j} = \bar{\nu} (\bar{I}_{i j})$. 
Moreover we define the divisor $\tilde{\lambda}$ 
on $\tilde{V}$ by $\tilde{\lambda} = \bar{\nu} (\bar{\lambda})$. 
The divisors $\tilde{I}_{i j}$ and $\tilde{\lambda}$ are 
non-singular rational curves with $\tilde{I}_{i j}^2 = -2$ and  
$\tilde{\lambda}^2 = -1$ respectively.  

Let $\tilde{\nu} : \tilde{V} \to V^{\prime}$ be 
the blowing-down of the $(-1)$-curve $\tilde{\lambda}$. 
For $i =1$, $2$ and $j=0$, $\infty$, we define the 
divisors $I_{i j}^{\prime}$, $G_i^{\prime}$, and $D_j^{\prime}$ 
on $V^{\prime}$ by 
$I_{i j}^{\prime} = (\tilde{\nu} \circ \bar{\nu})_* (\bar{I}_{i j})$, 
$G_i^{\prime} = (\tilde{\nu} \circ \bar{\nu})_* (\bar{G}_i)$, and 
$D_j^{\prime} = (\tilde{\nu} \circ \bar{\nu})_* (\bar{D}_j)$. 
The divisors $I_{i j}^{\prime}$, $G_i^{\prime}$, 
$D_0^{\prime}$, and $D_{\infty}^{\prime}$ are non-singular rational 
curves with ${I_{i j}^{\prime}}^2 = -2$, ${G_i^{\prime}}^2 = -1$, 
${D_0^{\prime}}^2 = - (d+2)/ 2$, and 
${D_{\infty}^{\prime}}^2 = (d-2)/ 2$ respectively. 

Let $\nu^{\prime} : V^{\prime} \to V^{\prime \prime}$ be the 
blowing-down of the two $(-1)$-curves $G_1^{\prime}$ and $G_2^{\prime}$. 
For $i=1$, $2$ and $j =0$, $\infty$, we define the divisors 
$I_{i j}^{\prime \prime}$ and $D_j^{\prime \prime}$ on $V^{\prime \prime}$ 
by $I_{i j}^{\prime \prime} = \nu^{\prime} (I_{ij}^{\prime})$ and 
$D_j^{\prime \prime} = \nu^{\prime} (D_j^{\prime})$. 
The divisors $I_{i j}^{\prime \prime}$, $D_0^{\prime \prime}$, 
and $D_{\infty}^{\prime \prime}$ are non-singular rational curves 
with ${I_{i j}^{\prime \prime}}^2 = -1$, 
${D_0^{\prime \prime}}^2 = -(d+2)/2$ and 
${D_{\infty}^{\prime \prime}}^2 = (d - 2)/2$ respectively. 

Let $\nu^{\prime \prime} : V^{\prime \prime} \to V^{\prime \prime \prime}$ 
be the blowing-down of the two $(-1)$-curves 
$I_{1 \infty}^{\prime \prime}$ and $I_{2 \infty}^{\prime \prime}$. 
For $i=1$, $2$ and $j=0$, $\infty$, we define the divisors 
$I_{i 0}^{\prime \prime \prime}$ and 
$D_j^{\prime \prime \prime}$ on $V^{\prime \prime \prime}$ by 
$I_{i 0}^{\prime \prime \prime} 
= \nu^{\prime \prime} (I_{i 0 }^{\prime \prime})$ and 
$D_j^{\prime \prime \prime} = \nu^{\prime \prime} (D_j^{\prime \prime})$. 
The divisors $I_{i 0}^{\prime \prime \prime}$, 
$D_0^{\prime \prime \prime}$, and $D_{\infty}^{\prime \prime \prime}$ 
are non-singular rational curves with 
${I_{i 0}^{\prime \prime \prime}}^2 = 0$, 
${D_0^{\prime \prime \prime}}^2 = -1$, and 
${D_{\infty}^{\prime \prime \prime}}^2 = 1$. 
By $K_{V^{\prime \prime \prime}}^2 = 8$, we see easily that 
$V^{\prime \prime \prime}$ is isomorphic to the Hirzebruch 
surface $\varSigma_1$ of degree $1$, where
$D_0^{\prime \prime \prime}$ and 
$I_{1 0}^{\prime \prime \prime} \sim I_{2 0}^{\prime \prime \prime}$
give the minimal section and the fiber class respectively. 
We define the point $v_1^{\prime \prime \prime} \in V^{\prime \prime \prime}$
by $v_1^{\prime \prime \prime} 
= \nu^{\prime \prime} (I_{1 \infty}^{\prime \prime })$. 
Note that we have 
$v_1^{\prime \prime \prime} \notin D_0^{\prime \prime \prime}$
if $d =0$, and 
$v_1^{\prime \prime \prime} \in D_0^{\prime \prime \prime}$ 
if $d=2$. 

We put 
$\nu = ( \nu^{\prime \prime} \circ  \nu^{\prime} 
\circ \tilde{\nu}  \circ \bar{\nu}) : 
\bar{V} \to V^{\prime \prime \prime} \simeq \varSigma_1$, 
and use the same symbols  
$I_{i \infty}^{\prime \prime}$, $G_i^{\prime}$, and $\tilde{\lambda}$ 
for the total transforms to $\bar{V}$ of the divisors 
$I_{i \infty}^{\prime \prime} \subset V^{\prime \prime}$, 
$G_i^{\prime} \subset V^{\prime}$, and 
$\tilde{\lambda} \subset \tilde{V}$ respectively. 
Note that the morphism 
$\nu : \bar{V} \to V^{\prime \prime \prime} \simeq \varSigma_1$ 
is a blowing-up at six points some of which are infinitely near.   

%The following description for our surfaces of the case $\chi =4$ 
%is essentially the same as the description in  
%Ciliberto-Mendes Lopes \cite[Section 1]{nonstandardpg3} 
%for the surfaces of the non-standard case 
%(see also \cite[(b) in Theorem 3.1]{onbicanonicalmaps}).   
%In the present paper, we have not put the  
%assumption of the non-birationality of bicanonical map.   
 
\begin{proposition} \label{prop:anotherdescription}
The linear system $|-4K_{\bar{V}} + \tilde{\lambda} + \bar{\lambda}^{\prime}|
= |\nu^* (-4 K_{V^{\prime \prime \prime}}) 
- 4\sum I_{i \infty}^{\prime \prime} - 4\sum G_i^{\prime} 
- 3\tilde{\lambda} -3 \bar{\lambda}^{\prime}|$
has no base point.  
Let $\bar{A}_0$ be a reduced member of 
$|-4K_{\bar{V}} + \tilde{\lambda} + \bar{\lambda}^{\prime}|$ that has 
at most negligible singularities, and satisfies 
$\bar{A}_0 \cap \bar{\lambda} = \emptyset$ and 
$\bar{A}_0 \cap \sum \bar{I}_{i j} = \emptyset$. 
%Let $\bar{h} : \bar{X} \to \bar{V}$ be the canonical resolution 
%of the double cover branched along $\bar{A}_0 + \bar{\lambda} + 
%\sum \bar{I}_{i j}$, 
%and $\bar{X} \to X$, the blowing-down of five $(-1)$-curves 
%$\bar{h}^{-1} (\bar{\lambda})$, 
%$\bar{h}^{-1} (\bar{I}_{1 0})$, 
%$\bar{h}^{-1} (\bar{I}_{2 0})$, 
%$\bar{h}^{-1} (\bar{I}_{1 \infty})$, and 
%$\bar{h}^{-1} (\bar{I}_{2 \infty})$. 
%
%
Let $\bar{X}$ be the canonical resolution 
of the double cover of $\bar{V}$ 
branched along $\bar{A}_0 + \bar{\lambda} + 
\sum \bar{I}_{i j}$, 
and $\bar{h} : \bar{X} \to \bar{V}$, 
the natural projection.  
Let $\bar{X} \to X$ be the blowing-down of five $(-1)$-curves 
$\bar{h}^{-1} (\bar{\lambda})$, 
$\bar{h}^{-1} (\bar{I}_{1 0})$, 
$\bar{h}^{-1} (\bar{I}_{2 0})$, 
$\bar{h}^{-1} (\bar{I}_{1 \infty})$, and 
$\bar{h}^{-1} (\bar{I}_{2 \infty})$. 
Then $X$ is a minimal surface with $c_1^2 = 2 \chi -1$, $\chi =4$, 
and $\mathrm{Tors} \simeq \mathbb{Z} / 2$. 
Conversely, for any minimal surface $X_{(1)}$ with these invariants,   
there exist configuration of $w_1$ and $w_1^{\prime}$ and 
a reduced member $\bar{A}_0$ as above such that 
the surface $X$ constructed by this procedure
is isomorphic to $X_{(1)}$. 
\end{proposition}

Proof. 
Since 
$K_{\bar{Z}} \sim \bar{\wp}^* (K_{\bar{V}}) + \sum J_{i j}$, 
$2J_{i j} = \bar{\wp}^* ({\bar{I}}_{i j})$,
$\bar{\wp}^* (\bar{\lambda}) = \sum \bar{e}_i$, 
$\bar{\wp}^* (\bar{\lambda}^{\prime}) = \sum \bar{e}_i^{\prime}$,
and $\tilde{\lambda} \sim \bar{\lambda} + \bar{\lambda}^{\prime}$, 
we have 
\[
 \bar{\wp}^* (-4 K_{\bar{V}} + \tilde{\lambda} + \bar{\lambda}^{\prime})
  \sim 
 \bar{r}^* ( (r^{\prime} \circ \tilde{r})^* ( - 4K_W) 
                      - 3 \sum e_i -3 \sum e_i^{\prime}). 
\]
By this together with 
\begin{equation} \label{eq:halfofsumIij}
\sum \bar{I}_{i j} = 
\nu^* (\sum I_{i 0}^{\prime \prime \prime}) - 2 \sum G_i^{\prime} \sim 
2 (I_{1 0}^{\prime \prime \prime} - \sum G_i^{\prime}), 
\end{equation}
we obtain 
\begin{align}
&\bar{r}^* \tilde{\Lambda}_8^+ (3 \sum e_i + 3 \sum e_i^{\prime}) = 
\bar{\wp}^* |- 4 K_{\bar{V}} + \tilde{\lambda} + \bar{\lambda}^{\prime}|, 
\notag \\
&\bar{r}^* \tilde{\Lambda}_8^- (3 \sum e_i + 3 \sum e_i^{\prime}) = 
\bar{\wp}^* |- 4 K_{\bar{V}} + \tilde{\lambda} + \bar{\lambda}^{\prime} 
- (I_{1 0}^{\prime \prime \prime} - \sum G_i^{\prime})|
+ \sum J_{i j}. 
\notag
\end{align}
Thus the absence of base points of 
$|- 4 K_{\bar{V}} + \tilde{\lambda} + \bar{\lambda}^{\prime}|$ 
follows from 3) of Lemma \ref{lm:baseptfreesystems}. 
Let $\bar{A}_0 \in 
|- 4 K_{\bar{V}} + \tilde{\lambda} + \bar{\lambda}^{\prime}|$,  
$\bar{h} : \bar{ X} \to \bar{V}$, and $X$ be 
a reduced member, the induced morphism, and the obtained surface 
respectively as in the statement. 
Then $\tilde{B}_0 = \bar{r} (\bar{\wp}^*(\bar{A}_0))$ satisfies 
all the conditions given in Proposition \ref{prop:sufficientprocedure}. 
Thus for this $\tilde{B}_0$, we obtain 
morphisms $\tilde{f} : \tilde{Y} \to \tilde{Z}$ and $p : \tilde{Y} \to Y$ 
as in Proposition \ref{prop:sufficientprocedure} 
and a minimal surface $Y/G$ with $c_1^2 = 2\chi -1$, $\chi =4$, 
and $\mathrm{Tors} \simeq \mathbb{Z} /2$. 
Note that the preimage by $\tilde{f}$ of the set  
$\{ (r^{\prime} \circ \tilde{r})^{-1}
(\varGamma_i \cap \varDelta_j) \}_{i=1, 2, \ j=0, \infty}$ 
is composed  of exactly eight points. 
We denote by $\bar{Y} \to \tilde{Y}$ the blowing-up at these 
eight points. Then the morphism $\tilde{f} : \tilde{Y} \to \tilde{Z}$ 
induces a generically two-to-one morphism 
$\bar{f} : \bar{Y} \to \bar{Z}$. Moreover,  
the natural free action by $G = \langle \iota |_W \rangle$ on $\tilde{Y}$ 
lifts to one on $\bar{Y}$ that is compatible with the induced 
action by $G$ on $\bar{Z}$. 
Thus $\bar{f} : \bar{Y} \to \bar{Z}$ induces 
a natural morphism $\bar{Y} /G \to \bar{V} = \bar{Z} /G$. 
Since the branch divisor of $\bar{Y} /G \to \bar{V} = \bar{Z} /G$ 
is $\bar{A}_0 + \bar{\lambda} + \sum \bar{I}_{i j}$, 
and $\bar{V}$ has no non-trivial torsion divisor, 
the morphism $\bar{Y} / G \to \bar{V}$ coincides with 
$\bar{h}: \bar{X} \to \bar{V}$. 
Thus by Proposition \ref{prop:sufficientprocedure}, 
$X \simeq Y/G$ is 
a minimal surface with $c_1^2 = 2 \chi -1$, $\chi =4$, 
and $\mathrm{Tors} \simeq \mathbb{Z} / 2$. 
The final assertion follows from Theorem \ref{thm:completedescription} 
and Lemma \ref{lm:configw1w1prime}. \qed

\begin{remark}  \label{rem:notecmdescr}
The description above for our surfaces of the case $\chi =4$ 
is almost the same as the description in  
Ciliberto-Mendes Lopes \cite[Section 1]{nonstandardpg3} 
of the surfaces of the non-standard case for the 
non-birationality of bicanonical maps  
(see also \cite[(b) in Theorem 3.1]{onbicanonicalmaps}).   
We emphasize here that in the present paper 
we have put neither the  
assumption of the non-birationality of bicanonical maps
nor the assumption of the absence of pencils of curves 
of genus $2$. 
%In fact, however, by the following proposition, it turns out that 
%our surfaces coincide exactly with those of a class treated in 
%Ciliberto-Mendes Lopes \cite{nonstandardpg3}.     
By the description above, it is almost clear that our surfaces 
coincide with those found in \cite{nonstandardpg3}. 
To be precise, however, 
by showing the following proposition, 
we shall prove that they in fact coincide.  
\end{remark}

\begin{proposition} \label{prop:bicanonical}
Any minimal surface $X$ with $c_1^2 = 2 \chi -1$, $\chi =4$, and 
$\mathrm{Tors} \simeq \mathbb{Z} /2$ has non-birational 
bicanonical map. Moreover, it has no pencil of curves of genus $2$. 
\end{proposition}

Proof. 
Let $X$ be a minimal surface with $c_1^2 = 2 \chi -1$, 
$\chi =4$, and $\mathrm{Tors} \simeq \mathbb{Z} / 2$, 
and assume that $X$ has a pencil of curves of genus $2$. 
This pencil is rational, since $c_1^2 \geq 2$. 
Let $\vartheta$ be a non-trivial $2$-torsion divisor, and 
$C$, a general member of the pencil. Then by the same method 
as in the proof of (i), Lemma $1. 2$ of \cite{on2chi-2}, 
we see that $| K_X + \vartheta| = |(\chi - 1) C| + D$ for 
a certain effective divisor $D$ satisfying 
$K_X D = 1$, $CD = 2$, $D^2 = 3 - 2\chi$, 
and $\mathcal{O}_C (D) \not\simeq \mathcal{O}_C (K_C)$.   
Since $K_X D = 1$, the divisor $D$ contains an irreducible 
curve $D_1$ satisfying $K_X D_1 = 1$, and all other components 
of $D$ are $(-2)$-curves. Then since  
$-3 \leq  D_1^2 + D_1 (D - D_1) 
= D_1 D = K_X D_1 - (\chi -1) C D_1 = 1 - (\chi -1) C D_1 $, 
we obtain $0 \leq C D_1 \leq 1$. 
Assume that $CD_1 = 1$. Then $D - D_1$ contains a $(-2)$-curve 
$D_2$ such that $C D_2 =1$, and so 
$-2 \leq D_2^2 + D_2 (D - D_2) = 
D D_2 \leq - (\chi - 1) CD_2$, hence a contradiction. 
Thus we have $C D_1 = 0$. In this case, if we let 
$D = D^{\prime} + D^{\prime \prime}$ be the decomposition of $D$ such 
that $D^{\prime}$ is the sum of all the irreducible components 
meeting $C$ with positive intersection number, and 
$D^{\prime \prime}$ is the sum of all the irreducible components 
meeting $C$ with intersection number $0$, 
then 
we have $D_1 \subset D^{\prime \prime}$ and $C D^{\prime}$ = 2. 
Since $D^{\prime}$ has at most two irreducible components and   
any irreducible component of $D^{\prime}$ is a $(-2)$-curve, 
by using exactly the same argument 
as in the proof of Lemma $2. 1$ of \cite{on2chi-2}, 
we obtain a contradiction.  Hence the surface $X$ has no pencil of 
curves of genus $2$.

Thus in order to prove Proposition \ref{prop:bicanonical}, 
it only remains to 
prove that $X$ has non-birational bicanonical map. 
For this purpose, let us use the notation 
in Proposition \ref{prop:anotherdescription}. 
We have 
$H^0 (\bar{h}_* \mathcal{O}_{\bar{X}} (2K_{\bar{X}})) 
=
H^0(\mathcal{O}_{\bar{V}}(2(K_{\bar{V}} + \varrho)))
\oplus
H^0(\mathcal{O}_{\bar{V}}(2K_{\bar{V}} + \varrho ))
$ for a certain divisor $\varrho$ with 
$\bar{A}_0 + \bar{\lambda} + \sum \bar{I}_{i,j} \sim 2 \varrho$. 
We however have 
\[
2 K_{\bar{V}} + \varrho  \quad
\sim    \quad 
\nu^* (I_{1 0}^{\prime \prime \prime})
-\sum G_i^{\prime} + \tilde{\lambda}, 
\] 
hence $h^0 (\mathcal{O}_{\bar{V}}(2K_{\bar{V}} + \varrho )) = 0$. 
This implies that the bicanonical map of $\bar{X}$ factors through 
the rational map of $\bar{V}$ associated to the linear system 
$|2 (K_{\bar{V}} + \varrho)|$. Hence we have the assertion. \qed 
\bigskip

%%%

To give a proof for Theorem \ref{thm:moduli}, we also need the 
following lemma: 

\begin{lemma} \label{lm:a0barcoh}
Let $\bar{V}$ be the smooth surface 
as in Proposition \ref{prop:anotherdescription}. 
Then for a member $\bar{A}_0 \in 
|- 4 K_{\bar{V}} + \tilde{\lambda} + \bar{\lambda}^{\prime}|$, 
the following equalities hold$:$
%\[
% h^i (\mathcal{O}_{\bar{V}} (\bar{A}_0)) = 
%    \begin{cases}
%           29  &  \text{$($$i=0$$)$} \\          
%            0  &  \text{$($$i=1$$)$} \\
%            0  &  \text{$($$i=2$$)$}. 
%    \end{cases} 
%\]
\[
 h^0 (\mathcal{O}_{\bar{V}} (\bar{A}_0)) = 29, \qquad 
 h^1 (\mathcal{O}_{\bar{V}} (\bar{A}_0)) = 0 , \qquad
 h^2 (\mathcal{O}_{\bar{V}} (\bar{A}_0)) = 0.
\]
\end{lemma}

Proof. 
Since we have 
$\bar{\wp}^* (\bar{A}_0) \in \bar{r}^* \tilde{\varLambda}_8^+ 
(3\sum e_i + 3\sum e_i^{\prime})$, 
the equality $h^i (\mathcal{O}_{\bar{V}} (\bar{A}_0)) = 0$ for 
any $i >0$ follows from Lemma \ref{lm:b0tildenonsing}. 
From this together with the Riemann--Roch theorem, 
we infer $h^0 (\mathcal{O}_{\bar{V}} (\bar{A}_0)) = 29$. \qed

As the first part of our proof for Theorem \ref{thm:moduli}, 
we shall show the following:  

\begin{lemma}  \label{lm:irredmoduli}
Any two minimal algebraic surfaces with 
$c_1^2 = 2 \chi -1$, $\chi = 4$, and 
$\mathrm{Tors} \simeq \mathbb{Z} /2$ are equivalent under deformation of 
complex structures. The coarse moduli space $\mathcal{M}$ for 
surfaces with these invariants is irreducible. 
\end{lemma}

Proof. 
Let $\mathcal{M}$ be the coarse moduli space 
for minimal surfaces $X$'s with $c_1^2 = 2 \chi -1$, $\chi = 4$, and 
$\mathrm{Tors} \simeq \mathbb{Z} /2$. 
In what follows, for a surface $X$ with these invariants, we denote by $[X]$ 
the point in $\mathcal{M}$ corresponding to the isomorphism class of $X$. 
As a reference point of $\mathcal{M}$, let us fix a surface 
$X_{(1)}$ with these invariants for which $d =0$ holds, 
the two points $w_1$ and $w_2$ lie neither on one and the same member  
of $| \varGamma |$ nor on that of $|\varDelta_0 |$, 
and $\bar{A}_0$ is smooth.  
Below, we shall give an irreducible component 
$\mathcal{M}_{(1)}$ of $\mathcal{M}$ containing $[X_{(1)}]$, 
and show that for any $X$ with these invariants  
we have $[X] \in \mathcal{M}_{(1)}$ and $X$ has 
the same deformation type as that of the reference surface $X_{(1)}$.    
We divide our situation into several cases according to 
$d$, the configuration of $w_i$'s and $w_i^{\prime}$'s, and smoothness 
of $\bar{A}_0$ for our $X$. 
In what follows, $\epsilon$ and $\epsilon_0$ will denote positive 
real numbers small enough.  We shall replace these numbers with smaller ones 
without mentioning it explicitly.

Case $1$-$1$: 
the case where $d=0$ holds, 
the two points $w_1$ and $w_2$ lie neither on one and the same member  
of $| \varGamma |$ nor on that of $|\varDelta_0 |$. 
This case splits into two subcases: case $1$-$1$-$1$ and 
case $1$-$1$-$2$. 

Case $1$-$1$-$1$:   
the subcace of case $1$-$1$ where $\bar{A}_0$ is smooth.  
From the point of view of description 
as in Proposition \ref{prop:anotherdescription}, 
this case corresponds to the case where 
$v_1^{\prime \prime \prime} 
\notin D_0^{\prime \prime \prime}$, 
$v_0^{\prime} \notin 
\sum I_{i j}^{\prime} + \sum G_i^{\prime} + \sum D_j^{\prime}$, 
and moreover $\bar{A}_0$ is smooth, where we put 
$v_1^{\prime \prime \prime} 
= \nu^{\prime \prime} (I_{1 \infty}^{\prime \prime})$ and  
$v_0^{\prime} = \tilde{\nu} (\tilde{\lambda})$.  
Note that for all $X$'s of this case $\tilde{V}$'s have 
one and the same isomorphism class. 
Let $\mathrm{pr}_{\tilde{V} \times \tilde{\lambda}} : 
\tilde{V} \times \tilde{\lambda} \to \tilde{\lambda} \simeq \mathbb{P}^1$ 
be the trivial family.
Let $\mathrm{pr}_{\tilde{V} \times \tilde{\lambda}, \tilde{V}}: 
\tilde{V} \times \tilde{\lambda} \to \tilde{V}$ 
be the first projection. 

Then we can easily construct an analytic family 
$\mathrm{pr}_{\bar{\mathcal{V}}} : 
\bar{\mathcal{V}} \to \tilde{\lambda} \simeq \mathbb{P}^1$
together with a projection
$\mathrm{pr}_{\bar{\mathcal{V}},  \tilde{V} \times \tilde{\lambda}}
: \bar{\mathcal{V}} \to \tilde{V} \times \tilde{\lambda}$
satisfying the following condition: 
for each $t \in \tilde{\lambda}$, 
the natural projection  
$\bar{V}_t = {\mathrm{pr}_{\bar{\mathcal{V}}}}^{-1} (t)
\to \tilde{V} = {\mathrm{pr}_{\tilde{V} \times \tilde{\lambda}}}^{-1} (t)$  
is the blowing-up  at $t \in \tilde{\lambda} \subset \tilde{V}$   
with exceptional divisor $\bar{\lambda}^{\prime}_t $. 
Let us denote by $\bar{\lambda}_t$ and $\tilde{\lambda}_t$
the strict transform and the total transform by 
$\bar{V}_t = {\mathrm{pr}_{\bar{\mathcal{V}}}}^{-1} (t)
\to \tilde{V} = {\mathrm{pr}_{\tilde{V} \times \tilde{\lambda}}}^{-1} (t)$
of the divisor $\tilde{\lambda}$, respectively.
We denote by $\mathrm{pr}_{\bar{\mathcal{V}}, \tilde{V}} : 
\bar{\mathcal{V}} \to \tilde{V}$ the composite of two projections  
$\mathrm{pr}_{\bar{\mathcal{V}},  \tilde{V} \times \tilde{\lambda}}$
and  $\mathrm{pr}_{\tilde{V} \times \tilde{\lambda}, \tilde{V}}$. 

%%

%%%%%%%%%%%%%%%%%%%%%%%%%%%%%

Consider the divisor 
$-4K_{\bar{\mathcal{V}}}  
+ {\mathrm{pr}_{\bar{\mathcal{V}}, \tilde{V}}}^* (\tilde{\lambda})
+ \cup_t \bar{\lambda}_t^{\prime}$ on 
$\bar{\mathcal{V}}$. 
The restriction to $\bar{V}_t$ of this divisor 
is linearly equivalent to 
$- 4K_{\bar{V}_t} + \tilde{\lambda}_t + \bar{\lambda}_t^{\prime}$.  
Since we have 
$h^1 (\mathcal{O}_{\bar{V}_t} 
(- 4K_{\bar{V}_t} + \tilde{\lambda}_t + \bar{\lambda}_t^{\prime})) =0$ and 
$h^0 (\mathcal{O}_{\bar{V}_t} 
(- 4K_{\bar{V}_t} + \tilde{\lambda}_t + \bar{\lambda}_t^{\prime})) = 29$ 
by Lemma \ref{lm:a0barcoh}, it follows that the direct image 
\[
 \mathcal{F}_0 = 
\mathrm{{pr}_{\bar{\mathcal{V}}}}_* \mathcal{O}_{\bar{\mathcal{V}}}
(-4K_{\bar{\mathcal{V}}}  
+ {\mathrm{pr}_{\bar{\mathcal{V}}, \tilde{V}}}^* (\tilde{\lambda})
+ \cup_t \bar{\lambda}_t^{\prime})
\]
is a locally free sheaf on $\tilde{\lambda} \simeq \mathbb{P}^1$ 
of rank $29$. 
We denote by $\mathcal{F}_0^{\vee}$ 
the dual sheaf 
of $\mathcal{F}_0$ on $\tilde{\lambda}$. 

Let $\mathrm{pr}_{\mathbb{P}} : 
\mathbb{P} = \mathbb{P} (\mathcal{F}_0^{\vee}) 
\to \tilde{\lambda}$ 
be the $\mathbb{P}^{28}$-bundle over $\tilde{\lambda}$ associated with  
$\mathcal{F}_0^{\vee}$.  
Then $\mathbb{P}$ is the projectivised 
total space of vector bundle $\mathcal{F}_0$. 
We consider the Cartesian diagram  
\[ 
\begin{CD}  
\bar{\mathcal{V}} \times_{\tilde{\lambda}} \mathbb{P}   
@>\text{$ $}>>  \mathbb{P} \\
@V\text{$ $}VV  @VV\text{$\mathrm{pr}_{\mathbb{P}}$}V \\
\bar{\mathcal{V}}  @>\text{$\mathrm{pr}_{\bar{\mathcal{V}}}$}>>   
\tilde{\lambda},
\end{CD}
\]
and denote by 
$\mathrm{pr}_{\bar{\mathcal{V}} 
\times_{\tilde{\lambda}} \mathbb{P}, \, \bar{\mathcal{V}}}
: \bar{\mathcal{V}} \times_{\tilde{\lambda}} \mathbb{P}
\to \bar{\mathcal{V}}$, 
$\mathrm{pr}_{\bar{\mathcal{V}} \times_{\tilde{\lambda}} 
\mathbb{P}, \, \mathbb{P}}: 
\bar{\mathcal{V}} \times_{\tilde{\lambda}} \mathbb{P}
\to \mathbb{P}$, and 
$\mathrm{pr}_{\bar{\mathcal{V}} \times_{\tilde{\lambda}} 
\mathbb{P}}: 
\bar{\mathcal{V}} \times_{\tilde{\lambda}} \mathbb{P}
\to \tilde{\lambda}$
the first projection, 
the second projection,
and the induced natural projection respectively.   

Let $\mathcal{O}_{\mathbb{P}} (1)$ be the tautological bundle of 
$ \mathrm{pr}_{\mathbb{P}} : 
\mathbb{P} = \mathbb{P} (\mathcal{F}_0^{\vee}) \to \tilde{\lambda}$. 
Then there exists a natural non-zero global section  
\[
 \varPsi_0 \in 
 H^0 ({\mathrm{pr}_{\bar{\mathcal{V}} \times_{\tilde{\lambda}} \mathbb{P}, 
\, \mathbb{P}}}^* 
\mathcal{O}_{\mathbb{P}} (1)
\otimes 
{\mathrm{pr}_{\bar{\mathcal{V}} \times_{\tilde{\lambda}} \mathbb{P}, 
\, \bar{\mathcal{V}}}}^* 
\mathcal{O}_{\bar{\mathcal{V}}}
( -4 K_{\bar{\mathcal{V}}} 
+ {\mathrm{pr}_{\bar{\mathcal{V}}, \, \tilde{V}}}^* (\tilde{\lambda})
+  \cup_t \bar{\lambda}_t^{\prime}))
\]
on $\bar{\mathcal{V}} \times_{\tilde{\lambda}} \mathbb{P}$ 
satisfying the following condition: 
for each open set $U \subset \tilde{\lambda}$ such that 
the restriction $\mathcal{F}_0 |_U$ is trivial, 
the restriction $\varPsi_0 
|_{{\mathrm{pr}_{\bar{\mathcal{V}} \times_{\tilde{\lambda}} 
\mathbb{P}}}^{-1} (U)}$ of $\varPsi_0$ to 
${\mathrm{pr}_{\bar{\mathcal{V}} \times_{\tilde{\lambda}} 
\mathbb{P}}}^{-1} (U)$ is given by 
$\varPsi_0 
|_{{\mathrm{pr}_{\bar{\mathcal{V}} \times_{\tilde{\lambda}} 
\mathbb{P}}}^{-1} (U)} = 
\sum_{i=1}^{29} a_i \psi_i$, where 
$\{ \psi_1, \ldots , \psi_{29} \}$ and 
$\{ a_1, \ldots , a_{29} \}$ are 
a base of $H^0 (\mathcal{F}_0 |_{U})$ and its dual base respectively
(note here that we have the natural isomorphism  
${\mathrm{pr}_{\mathbb{P}}}_* 
\mathcal{O}_{\mathbb{P}} (1) \simeq \mathcal{F}_0^{\vee}$). 
We denote by $\bar{\alpha}_0 = (\varPsi_0)$ 
the divisor on $\bar{\mathcal{V}} \times_{\tilde{\lambda}} \mathbb{P}$ 
defined by the section $\varPsi_0$.  

For each $u \in \mathbb{P}$, we put 
$t(u) = \mathrm{pr}_{\mathbb{P}} (u) \in \tilde{\lambda}$. 
Then we have the natural isomorphism 
${\mathrm{pr}_{\bar{\mathcal{V}} \times_{\tilde{\lambda}} \mathbb{P}, 
\, \mathbb{P} }}^{-1} (u) \simeq \bar{V}_{t(u)}$. 
Moreover, via this identification, the restriction 
$\bar{A}_{0 \, u} = \bar{\alpha}_0 
|_{{\mathrm{pr}_{\bar{\mathcal{V}} \times_{\tilde{\lambda}} \mathbb{P}, 
\, \mathbb{P} }}^{-1} (u) } $ 
$\in 
| -4 K_{\bar{V}_{t(u)}} + \tilde{\lambda}_{t(u)} 
+ \bar{\lambda}^{\prime}_{t(u)}| $ is a divisor on $\bar{V}_{t(u)}$ 
given by the local defining function  
$\sum_{i=1}^{29} a_i (u) \psi_i |_{\bar{V}_{t(u)}}$. 
Let $\mathbb{P}_0 \subset \mathbb{P}$ be the set of all $u$'s
such that ${\bar{A}_{0 \, u}}$ is a reduced smooth 
divisor satisfying 
$\bar{A}_{0 \, u} \cap 
(\bar{\lambda}_{t(u)} 
+ \sum_{i=1, 2, \ j=0, \infty} \bar{I}_{ij \, t(u)}) = \emptyset$, 
where 
$\bar{I}_{ij \, t}$ 
($t \in \tilde{\lambda}$) 
is the restriction to $\bar{V}_t$ of the divisor 
${\mathrm{pr}_{\bar{\mathcal{V}}, \tilde{V}}}^* (\tilde{I}_{ij})$.  
Then $\mathbb{P}_0$ is a non-empty Zariski open subset of 
$\mathbb{P} = \mathbb{P} (\mathcal{F}_0^{\vee})$. 
Since $\mathbb{P} \to \tilde{\lambda}$ is a $\mathbb{P}^{28}$-bundle 
over a non-singular rational curve $\tilde{\lambda} \simeq \mathbb{P}^1$, 
there exists a covering $\{ U_{\mu}^{\vee}\}_{\mu}$ of $\mathbb{P}$  by 
a finite number of Zariski open subsets  
$U_{\mu}^{\vee}$'s satisfying the following condition:  
for any $\mu$, the restriction 
$\mathcal{O}_{\mathbb{P}} (1) |_{U_{\mu}^{\vee}}$ is 
trivial, and $U_{\mu}^{\vee}$ is isomorphic to 
the $29$-dimensional linear space $\mathbb{A}^{29}$.    
We fix one such cover $\{ U_{\mu}^{\vee}\}_{\mu}$, and 
put $U_{\mu}^0 = U_{\mu}^{\vee} \cap \mathbb{P}_0$ for 
each $\mu$. 

Let  
$\mathrm{pr}_{\bar{\mathcal{V}},  V^{\prime \prime \prime}} 
: \bar{\mathcal{V}} \to V^{\prime \prime \prime}$
and 
$\mathrm{pr}_{\bar{\mathcal{V}},  V^{\prime}} 
: \bar{\mathcal{V}} \to V^{\prime}$
be the natural projections 
$\nu^{\prime \prime} \circ \nu^{\prime} \circ \tilde{\nu} \circ 
\mathrm{pr}_{\bar{\mathcal{V}},  \tilde{V}}$
and 
$\tilde{\nu} \circ \mathrm{pr}_{\bar{\mathcal{V}},  \tilde{V}}$
respectively.   
Then since the restriction 
to  
${\mathrm{pr}_{\bar{\mathcal{V}} 
\times_{\tilde{\lambda}} \mathbb{P}, \, \mathbb{P}}}^{-1} (U_{\mu}^0)$
of  
${\mathrm{pr}_{\bar{\mathcal{V}} \times_{\tilde{\lambda}} \mathbb{P}, \,
\mathbb{P}}}^*
\mathcal{O}_{\mathbb{P}}(1) $ 
is trivial,
it follows from (\ref{eq:halfofsumIij}) that 
the restriction to 
${\mathrm{pr}_{{\bar{\mathcal{V}} \times_{\tilde{\lambda}}
\mathbb{P}, \, \mathbb{P}}}}^{-1} (U_{\mu}^0)$ 
of the divisor 
\begin{equation}    \label{eq:familyofbranchdivisors}
\bar{\alpha}_0 + {\mathrm{pr}_{\bar{\mathcal{V}} \times_{\tilde{\lambda}}
\mathbb{P}, \, \bar{\mathcal{V}}} }^* 
(\cup_t \bar{\lambda}_t 
+ {\mathrm{pr}_{{\bar{\mathcal{V}} , \tilde{V}}}}^* (\sum \tilde{I}_{i j}))  
\end{equation}
is linearly equivalent to twice the restriction to 
 ${\mathrm{pr}_{{\bar{\mathcal{V}} \times_{\tilde{\lambda}}
\mathbb{P}, \, \mathbb{P}}}}^{-1} (U_{\mu}^0)$ 
of the divisor 
\[
{\mathrm{pr}_{\bar{\mathcal{V}} \times_{\tilde{\lambda}}
\mathbb{P}, \, \bar{\mathcal{V}}} }^* 
(-2 K_{\bar{\mathcal{V}}} 
+ {\mathrm{pr}_{\bar{\mathcal{V}}, \tilde{V}}}^* (\tilde{\lambda})
+ {\mathrm{pr}_{\bar{\mathcal{V}}, V^{\prime \prime \prime}}}^* 
     (I_{1 0}^{\prime \prime \prime})
- {\mathrm{pr}_{\bar{\mathcal{V}}, V^{\prime}}}^* 
     (\sum G_i^{\prime})). 
\]
So let 
$\bar{\mathcal{X}}^{(\mu)} \to 
{\mathrm{pr}_{{\bar{\mathcal{V}} \times_{\tilde{\lambda}}
\mathbb{P}, \, \mathbb{P}}}}^{-1} (U_{\mu}^0)$
be the double cover branched along 
the restriction to 
${\mathrm{pr}_{{\bar{\mathcal{V}} \times_{\tilde{\lambda}}
\mathbb{P}, \, \mathbb{P}}}}^{-1} (U_{\mu}^0)$ 
of the divisor (\ref{eq:familyofbranchdivisors}). 
Composing this morphism with 
the projection 
$\mathrm{pr}_{\bar{\mathcal{V}} \times_{\tilde{\lambda}} \mathbb{P}, \, 
\mathbb{P}} : 
\bar{\mathcal{V}} \times_{\tilde{\lambda}} \mathbb{P} \to \mathbb{P}$, 
we obtain an analytic family 
$\mathrm{pr}_{\bar{\mathcal{X}}^{(\mu)}} : 
\bar{\mathcal{X}}^{(\mu)} \to U_{\mu}^0$. 
For each $u \in U_{\mu}^0$, we put 
$\bar{X}^{(\mu)}_u = {\mathrm{pr}_{\bar{\mathcal{X}}^{(\mu)}}}^{-1} (u)$.  
The inverse image by 
$\bar{\mathcal{X}}^{(\mu)} 
\to {\mathrm{pr}_{{\bar{\mathcal{V}} \times_{\tilde{\lambda}}
\mathbb{P}, \, \mathbb{P}}}}^{-1} (U_{\mu}^0)$ 
of   
${\mathrm{pr}_{\bar{\mathcal{V}} \times_{\tilde{\lambda}}
\mathbb{P}, \, \bar{\mathcal{V}}} }^* 
(\cup_t \bar{\lambda}_t 
+ {\mathrm{pr}_{{\bar{\mathcal{V}} , \tilde{V}}}}^* (\sum \tilde{I}_{i j}))$
gives a family over $U_{\mu}^0$ whose fiber over each $u \in U_{\mu}^0$ is 
a sum of five disjoint $(-1)$-curves on $\bar{X}^{(\mu)}_u$. 
Blowing down this family of disjoint five $(-1)$-curves relatively to 
$\mathrm{pr}_{\bar{\mathcal{X}}^{(\mu)}} : 
\bar{\mathcal{X}}^{(\mu)} \to U_{\mu}^0$, 
we obtain an analytic family 
$\mathrm{pr}_{\mathcal{X}^{(\mu)}} : 
\mathcal{X}^{(\mu)} \to U_{\mu}^0$. 
For each $u \in U_{\mu}^0$, we put 
$X^{(\mu)}_u = {\mathrm{pr}_{\mathcal{X}^{(\mu)}}}^{-1} (u)$. 
Then by Proposition \ref{prop:anotherdescription}, 
for each $u \in U_{\mu}^0$, the fiber $X^{(\mu)}_u$ is 
a minimal surface with $c_1^2 = 2\chi -1$, $\chi =4$, 
and $\mathrm{Tors} \simeq \mathbb{Z} / 2$.   

Let $U_{\mu}^0 \to \mathcal{M}$ be the natural morphism 
induced from the family $\mathrm{pr}_{\mathcal{X}^{(\mu)}} : 
\mathcal{X}^{(\mu)} \to U_{\mu}^0$, i.e., the morphism given by 
$u \mapsto [X^{(\mu)}_u]$, where $[X^{(\mu)}_u]$ is a point 
in $\mathcal{M}$ corresponding to the isomorphism class of the 
fiber $X^{(\mu)}_u$. 
The two morphisms $U_{\mu_1}^0 \to \mathcal{M}$ 
and $U_{\mu_2}^0 \to \mathcal{M}$
coincide on $U_{\mu_1}^0 \cap U_{\mu_2}^0$. 
Thus gluing $U_{\mu}^0 \to \mathcal{M}$'s, 
we obtain a morphism 
$\mathbb{P}_0 \to \mathcal{M}$ 
given locally by $u \mapsto [X^{(\mu)}_u]$. 
Since $\mathbb{P}_0$ is irreducible, 
the image of this $\mathbb{P}_0 \to \mathcal{M}$ 
lies on an irreducible component of the moduli space $\mathcal{M}$. 
We fix one such irreducible component 
and denote it by $\mathcal{M}_{(1)}$. 

Now let $X$ be a minimal surface with $c_1^2 = 2\chi -1$, 
$\chi =4$, and $\mathrm{Tors} \simeq \mathbb{Z} /2$ for which 
$d=0$ holds, the two points $w_1$ and $w_2$ lie neither on 
one and the same member of $|\varDelta_0|$ 
nor on that of $|\varGamma|$, and $\bar{A}_0$ is smooth.  
Then by Proposition \ref{prop:anotherdescription} and the 
construction of $\mathrm{pr}_{\mathcal{X}^{(\mu)}} : 
\mathcal{X}^{(\mu)} \to U_{\mu}^0$ above,
there exist a $\mu$ and a $u \in U_{\mu}^0$ such that 
$X \simeq X^{(\mu)}_u$ holds. 
Since $\mathbb{P}_0$ is connected, we infer that 
$X$ has the same deformation type as  
that of the reference surface $X_{(1)}$. 
Moreover, we infer that the corresponding point $[X]$ lies on 
the irreducible component $\mathcal{M}_{(1)}$% 
%(or more precisely, on the image of $\mathbb{P}_0 \to \mathcal{M}$)
.

Case $1$-$1$-$2$: the subcase of case $1$-$1$ where 
$\bar{A}_0$ is singular. 
Let $X$ be a minimal surface with $c_1^2 = 2\chi -1$, 
$\chi =4$, and $\mathrm{Tors} \simeq \mathbb{Z} /2$ of 
this case. 
In this case, $\bar{A}_0$ has at most negligible singularities,  
and by Proposition \ref{prop:anotherdescription}, 
the linear system $|\bar{A}_0|$ has no base point.  
Thus by the same method as in 
\cite[Proof of Theorem $4$]{quintic}, 
we obtain an analytic family 
$\mathrm{pr}_{\mathcal{X}} : \mathcal{X} \to N
= \{ u \in \mathbb{C} : |u| < \epsilon \}$ 
of minimal surfaces with $c_1^2 = 2 \chi -1$, $\chi =4$, and 
$\mathrm{Tors} \simeq \mathbb{Z} / 2$ 
such that  
$X_u = {\mathrm{pr}_{\mathcal{X}}}^{-1} (u)$ is 
of case $1$-$1$-$1$
for each $u \neq 0 \in N$, 
and 
$X_0 = {\mathrm{pr}_{\mathcal{X}}}^{-1} (0) 
\simeq X$.  
%From this together with 
%the results for case $1$-$1$-$1$, 
%we infer that our $X$ has the same deformation type as 
%that of the reference surface $X_{(1)}$. 
%Moreover, we infer that the corresponding point 
%$[X]$ lies on the irreducible component $\mathcal{M}_{(1)}$ 
%%(or more precisely, on the closure of the image of 
%%$\mathbb{P}_0 \to \mathcal{M}$) 
%given in the proof for case $1$-$1$-$1$.  
From this together with 
the results for case $1$-$1$-$1$, 
we infer that $X$ has the same deformation type as 
that of the reference surface $X_{(1)}$,  
and that the point  
$[X]$ lies on the irreducible component 
$\mathcal{M}_{(1)}$ in the proof for case $1$-$1$-$1$.

%%%

%
%

Case $1$-$2$: the case where 
$d=0$ holds, 
and the two points $w_1$ and $w_2$ lie on one and the same 
member of $| \varGamma |$ or $| \varDelta_0 |$. 
In this case, by Remark \ref{rem:gamma1gamma2}, 
we may assume that the two points $w_1$ and $w_2$ lie 
on the member $\varGamma_1 \in |\varGamma|$.  
Then by Lemma \ref{lm:configw1w1prime}, 
we have $w_1^{\prime} \notin {r^{\prime}}^{-1}_* (\varGamma_1)$. 
This case is divided into two subcases: 
case $1$-$2$-$1$ and case $1$-$2$-$2$. 

Case $1$-$2$-$1$: the subcase of case $1$-$2$ where 
$\bar{A}_0$ is smooth. 
From the point of view of description 
as in Proposition \ref{prop:anotherdescription}, 
this case corresponds to the case where 
$v_1^{\prime \prime \prime} \notin D_0^{\prime \prime \prime}$,  
$v_0^{\prime} \in G_1^{\prime} \setminus (\sum_j  I_{1 j}^{\prime})$, 
and moreover $\bar{A}_0$ is smooth, where we put 
$v_1^{\prime \prime \prime} 
= \nu^{\prime \prime} (I_{1 \infty}^{\prime \prime})$ and 
$v_0^{\prime} = \tilde{\nu} (\tilde{\lambda})$. 
Let $X$ be a minimal surface with $c_1^2 = 2\chi -1$, $\chi =4$, 
and $\mathrm{Tors} \simeq \mathbb{Z} /2$ of this case. 
Let $\epsilon$ be a positive real number small enough. 
We put $N = \{ t \in \mathbb{C} : |t| < \epsilon \}$, and 
denote by 
$\mathrm{pr}_{V^{\prime} \times N} :  V^{\prime} \times N \to N $ 
the trivial family. 
Let $\mathrm{pr}_{ V^{\prime} \times  N, \, V^{\prime}} 
:  V^{\prime} \times N \to  V^{\prime}$ be the first projection. 
%%%%%%%%%%

Then we can easily construct anatitic families 
$\mathrm{pr}_{\tilde{\mathcal{V}}} : \tilde{\mathcal{V}} \to N$, 
$\mathrm{pr}_{\bar{\mathcal{V}}} : \bar{\mathcal{V}} \to N$ together 
with projections 
$\mathrm{pr}_{\tilde{\mathcal{V}}, \, V^{\prime} \times N} : 
\tilde{\mathcal{V}} \to V^{\prime} \times N$, 
$\mathrm{pr}_{\bar{\mathcal{V}}, \tilde{\mathcal{V}}} : 
\bar{\mathcal{V}} \to \tilde{\mathcal{V}}$ 
satisfying the following conditions: 
for each $t \in N$, the projection   
$\tilde{V}_t = {\mathrm{pr}_{\tilde{\mathcal{V}}}}^{-1} (t) 
\to V^{\prime} = {\mathrm{pr}_{V^{\prime} \times N}}^{-1} (t)$ 
is the blowing-up at $v^{\prime} (t)$ with  
exceptional divisor $\tilde{\lambda}_t$, 
where $v^{\prime}: N \to V^{\prime} \times N$ is a holomorphic 
section of the analytic family 
$\mathrm{pr}_{V^{\prime} \times N} :  V^{\prime} \times N \to N$
; 
for each $t \in N$, the projection  
$\bar{V}_t = {\mathrm{pr}_{\bar{\mathcal{V}}}}^{-1} (t) \to 
\tilde{V}_t = {\mathrm{pr}_{\tilde{\mathcal{V}}}}^{-1} (t) $ 
is the blowing-up at $\tilde{v} (t)$ with exceptional divisor 
$\bar{\lambda}^{\prime}_t$, 
where $\tilde{v} : N \to \tilde{\mathcal{V}}$ 
is a holomorphic section of 
the analytic family 
$\mathrm{pr}_{\tilde{\mathcal{V}}} : \tilde{\mathcal{V}} \to N$
; $v^{\prime} (0) = v^{\prime}_0 (= \tilde{\nu} (\tilde{\lambda}))$
holds, and 
$\mathrm{pr}_{ V^{\prime} \times  N, \, V^{\prime}} (v^{\prime} (t)) \in 
\sum G_i^{\prime} + \sum D_j^{\prime} + \sum I_{ij}^{\prime}$ 
if and only if $t =0$
; 
$\tilde{v} (0) = \tilde{v}_0 (= \bar{\nu} (\bar{\lambda}^{\prime}) )$
holds,  
and $\tilde{v} (t) \in \tilde{\lambda}_t$ for any $t \in N$. 
Note that from the conditions above, we have in particular 
$\tilde{V}_0 = \tilde{V}$ and $\bar{V}_0 = \bar{V}$.
Let us denote by $\bar{\lambda}_{t}$ the strict transform 
of $\tilde{\lambda}_t$ by 
$\bar{V}_t = {\mathrm{pr}_{\bar{\mathcal{V}}}}^{-1} (t) \to 
\tilde{V}_t = {\mathrm{pr}_{\tilde{\mathcal{V}}}}^{-1} (t) $.   

Consider the divisor 
$-4 K_{\bar{\mathcal{V}}} + 
{\mathrm{pr}_{\bar{\mathcal{V}}, \tilde{\mathcal{V}}}}^* 
(\cup_t \tilde{\lambda}_t) + \cup_t \bar{\lambda}^{\prime}_t$ 
on $\bar{\mathcal{V}}$. 
The restriction to $\bar{V}_0 = \bar{V}$ 
of this divisor is linearly equivalent to 
$-4K_{\bar{V}} + \tilde{\lambda} + \bar{\lambda}^{\prime}$. 
Since we have 
$h^1 (\mathcal{O}_{\bar{V}} 
(-4K_{\bar{V}} + \tilde{\lambda} + \bar{\lambda}^{\prime})) =0$ 
by Lemma \ref{lm:a0barcoh},   
there exists a non-zero global section 
$\varPsi \in 
H^0 (\mathcal{O}_{\bar{\mathcal{V}}}
(-4 K_{\bar{\mathcal{V}}} + 
{\mathrm{pr}_{\bar{\mathcal{V}}, \tilde{\mathcal{V}}}}^* 
(\cup_t \tilde{\lambda}_t) + \cup_t \bar{\lambda}^{\prime}_t))$ 
on $\bar{\mathcal{V}}$ satisfying the following conditions: 
the restriction $(\varPsi) |_{\bar{V}_0}$ to $\bar{V}_0$ 
coincides with $\bar{A}_0$, 
where $( \varPsi )$ denotes 
the divisor on $\bar{\mathcal{V}}$ defined by 
the global section $\varPsi$; 
for any $t \in N$, the restriction  
$(\varPsi) |_{\bar{V}_t}$ to $\bar{V}_t$ is a reduced non-singular 
divisor on $\bar{V}_t$; 
the divisor $(\varPsi)$ does not intersects 
$\cup_t \bar{\lambda}_t + 
{\mathrm{pr}_{\bar{\mathcal{V}}, V^{\prime}}}^* (\sum I_{i j}^{\prime})$, 
where
$\mathrm{pr}_{\bar{\mathcal{V}}, V^{\prime}} : 
\bar{\mathcal{V}} \to V^{\prime}$ is  
the composite of three projections  
$\mathrm{pr}_{\bar{\mathcal{V}}, \tilde{\mathcal{V}}}$, 
$\mathrm{pr}_{\tilde{\mathcal{V}}, \, V^{\prime} \times N}$, 
and $\mathrm{pr}_{V^{\prime} \times N, \, V^{\prime}}$.  

Let $\bar{\mathcal{X}} \to \bar{\mathcal{V}}$ be the 
double cover branched along 
$( \varPsi ) + \cup_t \bar{\lambda}_t + 
{\mathrm{pr}_{\bar{\mathcal{V}}, V^{\prime}}}^* (\sum I_{i j}^{\prime})$. 
Then composing this morphism with 
the projection $\mathrm{pr}_{\bar{\mathcal{V}}} : 
\bar{\mathcal{V}} \to N$, we obtain an analytic family 
$\mathrm{pr}_{\bar{\mathcal{X}}} : 
\bar{\mathcal{X}} \to N$.   
For each $t \in N$, we put 
$\bar{X}_t = {\mathrm{pr}_{\bar{\mathcal{X}}}}^{-1} (t)$. 
The inverse image by $\bar{\mathcal{X}} \to \bar{\mathcal{V}}$ of 
$\cup_t \bar{\lambda}_t + 
{\mathrm{pr}_{\bar{\mathcal{V}}, V^{\prime}}}^* (\sum I_{i j}^{\prime})$ 
gives a family over $N$ whose fiber over each $t \in N$ is 
a disjoint union of five $(-1)$-curves on $\bar{X}_t$.  
Blowing down this family of five $(-1)$-curves relatively to 
$\mathrm{pr}_{\bar{\mathcal{X}}} : \bar{\mathcal{X}} \to N$, 
we obtain an analytic family 
$\mathrm{pr}_{\mathcal{X}} : \mathcal{X} \to N$. 
Then by the construction of 
$\mathrm{pr}_{\mathcal{X}} : \mathcal{X} \to N$ above, 
we have $X_0 = {\mathrm{pr}_{\mathcal{X}}}^{-1} (0) \simeq X$, 
and for each $t \neq 0 \in N$, the fiber 
$X_t = {\mathrm{pr}_{\mathcal{X}}}^{-1} (t)$ 
is a minimal surface with $c_1^2 = 2\chi -1$, $\chi =4$, 
and $\mathrm{Tors} \simeq \mathbb{Z} /2$ of case $1$-$1$-$1$. 
%From this together with the results for 
%case $1$-$1$-$1$, we infer that $X$ has the same deformation 
%type as that of the reference surface $X_{(1)}$. 
%Moreover, we infer that the corresponding point 
%$[X]$ lies on the irreducible component $\mathcal{M}_{(1)}$ 
%%(or more precisely, on the closure of the image of 
%%$\mathbb{P}_0 \to \mathcal{M}$) 
%given in the proof for case $1$-$1$-$1$.  
From this together with the results for 
case $1$-$1$-$1$, we infer that $X$ has the same deformation 
type as that of the reference surface $X_{(1)}$,  
and that the point $[X]$ lies on the irreducible component 
$\mathcal{M}_{(1)}$ in the proof for case $1$-$1$-$1$.

%%%%%%%%%%

Case $1$-$2$-$2$: the subcase of case $1$-$2$ where 
$\bar{A}_0$ is singular. 
%Note that in this case $\bar{A}_0$ has at most negligible 
%singularities. 
Let $X$ be a minimal surface with 
$c_1^2 = 2\chi -1$, $\chi =4$, and $\mathrm{Tors} \simeq \mathbb{Z} /2$ 
of this case. 
%Then using the same argument as in case $1$-$1$-$2$,  
%we infer from the result for case $1$-$2$-$1$   
%that $X$ has the same deformation type as 
%that of the reference surface $X_{(1)}$. 
%Moreover, we infer that the corresponding point $[X]$ lies on the 
%irreducible component $\mathcal{M}_{(1)}$ 
%%(or more precisely, on the closure of the image of 
%%$\mathbb{P}_0 \to \mathcal{M}$) 
%given in the proof for case $1$-$1$-$1$. 
%%
Then using the same argument as in case $1$-$1$-$2$,  
we infer from the results for case $1$-$2$-$1$   
that $X$ has the same deformation type as 
that of the reference surface $X_{(1)}$,  
and that the point $[X]$ lies on the 
irreducible component $\mathcal{M}_{(1)}$ 
in the proof for case $1$-$1$-$1$. 

Case $2$-$1$: the case where $d=2$ holds, and the two points 
$w_1$ and $w_2$ do not lie on one and the same member of $| \varGamma |$. 
Note that in this case we have $v_0^{ \prime } \notin D_0^{\prime} $ 
by Lemma \ref{lm:configw1w1prime}, where we put 
$v_0^{ \prime } = \tilde{\nu} (\tilde{\lambda})$. 
This case splits into two subcases: 
case $2$-$1$-$1$ and case $2$-$1$-$2$. 

Case $2$-$1$-$1$: the subcase of case $2$-$1$ where 
$\bar{A}_0$ is smooth.   
From the point of view of description 
as in Proposition \ref{prop:anotherdescription}, 
this case corresponds to the case where 
$v_1^{\prime \prime \prime} \in D_0^{\prime \prime \prime}$,  
$v_0^{\prime} \notin D_0^{\prime} 
+ \sum G_i^{\prime} + \sum I_{ij}^{\prime}$, 
and moreover $\bar{A}_0$ is smooth, 
where we put 
$v_1^{\prime \prime \prime} = 
\nu^{\prime \prime} ( I_{1 \infty}^{\prime \prime})$  
and 
$v_0^{\prime} = \tilde{\nu} (\tilde{\lambda})$. 
Let $X$ be a minimal surface with $c_1^2 = 2\chi -1$, $\chi =4$, 
and $\mathbb{Z} /2$ of this case. 
Let $\epsilon$ be a positive real number small enough. 
We put $N = \{ t \in \mathbb{C} : |t| < \epsilon \}$, 
and denote by 
$\mathrm{pr}_{V^{\prime \prime \prime} \times N }
: V^{\prime \prime \prime} \times N \to  N $ the trivial family. 
Let 
$\mathrm{pr}_{V^{\prime \prime \prime} \times N, \, V^{\prime \prime \prime}}
: V^{\prime \prime \prime} \times N \to  V^{\prime \prime \prime}$ 
be the first projection. 
%%%
Let us take a holomorphic section                   
$v_{(1)}^{\prime \prime \prime} : 
N \to V^{\prime \prime \prime} \times N$ 
satisfying the following conditions: 
$\mathrm{pr}_{V^{\prime \prime \prime} \times N, \, V^{\prime \prime \prime}}
(v_{(1)}^{\prime \prime \prime} (0)) = v_1^{\prime \prime \prime} $ 
($= \nu^{\prime \prime} (I_{1 \infty}^{\prime \prime})$) holds; 
$\mathrm{pr}_{V^{\prime \prime \prime} \times N, \, V^{\prime \prime \prime}}
(v_{(1)}^{\prime \prime \prime} (t)) \in I_{1 0}^{\prime \prime \prime}$ 
for any $t \in N$; 
$\mathrm{pr}_{V^{\prime \prime \prime} \times N, \, V^{\prime \prime \prime}}
(v_{(1)}^{\prime \prime \prime} (t)) = v_1^{\prime \prime \prime}$ 
if and only if $t=0$.

%%%
Recall that the configuration corresponding to Case $1$-$1$-$1$ was 
$v_1^{\prime \prime \prime} 
\notin D_0^{\prime \prime \prime}$, 
$v_0^{\prime} \notin 
\sum I_{i j}^{\prime} + \sum G_i^{\prime} + \sum D_j^{\prime}$. 
Thus using the holomorphic section 
$v_{(1)}^{\prime \prime \prime} : 
N \to V^{\prime \prime \prime} \times N$ above, 
and by the same arguement as in the proof for case $1$-$2$-$1$,  
we obtain an analytic family 
$\mathrm{pr}_{\mathcal{X}} : \mathcal{X} 
\to N= \{ t \in \mathbb{C} : |t| < \epsilon \}$ 
satisfying the following conditions:   
$X_0 = {\mathrm{pr}_{\mathcal{X}}}^{-1} (0) \simeq X$;  
and for each $t \neq 0 \in N$, the fiber 
$X_t = {\mathrm{pr}_{\mathcal{X}}}^{-1} (t)$ 
is a minimal surface with $c_1^2 = 2\chi -1$, $\chi =4$, 
and $\mathrm{Tors} \simeq \mathbb{Z} /2$ of case $1$-$1$-$1$. 
From this together with the results for 
case $1$-$1$-$1$, we infer that $X$ has the same deformation 
type as that of the reference surface $X_{(1)}$,  
and that the point $[X]$ lies on the irreducible component 
$\mathcal{M}_{(1)}$ in the proof for case $1$-$1$-$1$.  

%%%%%%%%%%%%%%%%%%%%%%%%

%
%Case $2$-$1$-$2$: the subcase of case $2$-$1$ where 
%$\bar{A}_0$ is singular. 
%Note that in this case $\bar{A}_0$ has at most negligible 
%singularities. Let $X$ be a minimal surface with 
%$c_1^2 = 2\chi -1$, $\chi =4$, and $\mathrm{Tors} \simeq \mathbb{Z} /2$ 
%of this case. 
%Then we have $h^1 (\mathcal{O}_{\bar{V}} (\bar{A}_0)) = 0$ by 
%Lemma \ref{lm:a0barcoh}. 
%Thus, using the same method as in the proof for 
%case $1$-$1$-$2$,  we obtain an analytic family 
%$\mathrm{pr}_{\mathcal{X}} : \mathcal{X} \to 
%N = \{ u \in \mathbb{C} : |u|  <  \epsilon \}$ 
%satisfying the following conditions:  
%the central fiber $X_0 = {\mathrm{pr}_{\mathcal{X}}}^{-1} (0)$ is 
%isomorphic to $X$; the fiber $X_u = {\mathrm{pr}_{\mathcal{X}}}^{-1} (u)$ 
%over each $u \neq 0 \in N$ is a minimal surface with $c_1^2 = 2\chi -1$, 
%$\chi = 4$, and $\mathrm{Tors} \simeq \mathbb{Z} /2$ of case $2$-$1$-$1$. 
%From this together with the results for case $2$-$1$-$1$, 
%we infer that $X$ has the same deformation type as 
%that of the reference surface $X_{(1)}$. 
%Moreover, we infer that the corresponding point $[X]$ lies on the 
%irreducible component $\mathcal{M}_{(1)}$ 
%(or more precisely, on the closure of the image of 
%$\mathbb{P}_0 \to \mathcal{M}$) given 
%in the proof for case $1$-$1$-$1$. 
%

%

Case $2$-$1$-$2$: the subcase of case $2$-$1$ where 
$\bar{A}_0$ is singular. 
%Note that in this case $\bar{A}_0$ has at most negligible 
%singularities. 
Let $X$ be a minimal surface with 
$c_1^2 = 2\chi -1$, $\chi =4$, and $\mathrm{Tors} \simeq \mathbb{Z} /2$ 
of this case. 
%Then using the same argument as in case $1$-$1$-$2$, 
%we infer from the result for case $2$-$1$-$1$ 
%that $X$ has the same deformation type as 
%that of the reference surface $X_{(1)}$. 
%Moreover, we infer that the corresponding point $[X]$ lies on the 
%irreducible component $\mathcal{M}_{(1)}$ 
%%(or more precisely, on the closure of the image of 
%%$\mathbb{P}_0 \to \mathcal{M}$) 
%given in the proof for case $1$-$1$-$1$. 
%
Then using the same argument as in case $1$-$1$-$2$, 
we infer from the results for case $2$-$1$-$1$ 
that $X$ has the same deformation type as 
that of the reference surface $X_{(1)}$,  
and that the point $[X]$ lies on the 
irreducible component $\mathcal{M}_{(1)}$ 
%(or more precisely, on the closure of the image of 
%$\mathbb{P}_0 \to \mathcal{M}$) 
in the proof for case $1$-$1$-$1$.

Case $2$-$2$: 
the case where $d=2$ holds, 
the two points $w_1$ and $w_2$ lie on one and the same 
member of $| \varGamma |$. 
Note that in this case we have 
$v_0^{\prime} \notin D_0^{\prime}$ by Lemma \ref{lm:configw1w1prime}, 
where we put 
$v_0^{\prime} = \tilde{\nu} (\tilde{\lambda})$. 
In this case, by Remark \ref{rem:gamma1gamma2}, 
we may assume that the two points $w_1$ and $w_2$ lie 
on the member $\varGamma_1 \in |\varGamma|$.  
Then by Lemma \ref{lm:configw1w1prime}, 
we have $w_1^{\prime} \notin {r^{\prime}}^{-1}_* (\varGamma_1)$. 
This case is divided into two subcases: 
case $2$-$2$-$1$ and case $2$-$2$-$2$. 

Case $2$-$2$-$1$: 
the subcase of case $2$-$2$ where $\bar{A}_0$ is 
smooth. 
From the point of view of description 
as in Proposition \ref{prop:anotherdescription}, 
this case corresponds to the case where 
$v_1^{\prime \prime \prime} \in D_0^{\prime \prime \prime}$,  
$v_0^{\prime} \in G_1^{\prime} \setminus (\sum_j  I_{1 j}^{\prime})$, 
and moreover $\bar{A}_0$ is smooth, where we put 
$v_1^{\prime \prime \prime} 
= \nu^{\prime \prime} (I_{1 \infty}^{\prime \prime})$ and 
$v_0^{\prime} = \tilde{\nu} (\tilde{\lambda})$. 
Note that by Lemma \ref{lm:configw1w1prime}, we have $\tilde{v}_0 \notin 
{\tilde{\nu}}^{-1}_* (G_1^{\prime}) $, 
where we put $\tilde{v}_0 = \bar{\nu} (\bar{\lambda}^{\prime})$. 
Let $X$ be a minimal surface with $c_1^2 = 2\chi -1$, $\chi =4$, 
and $\mathrm{Tors} \simeq \mathbb{Z} /2$ of this case. 
Let $\epsilon$ be a positive real number small enough. 
We put $N = \{ t \in \mathbb{C} : |t| < \epsilon \}$, and 
denote by 
$\mathrm{pr}_{V^{\prime} \times N} :  V^{\prime} \times N \to N $ 
the trivial family. 
Let $\mathrm{pr}_{ V^{\prime} \times  N, \, V^{\prime}} 
:  V^{\prime} \times N \to  V^{\prime}$ be the first projection. 
Let us take a holomorphic section 
$v^{\prime} : N \to V^{\prime} \times N$ satisfying the 
following conditions: 
$\mathrm{pr}_{ V^{\prime} \times  N, \, V^{\prime}} (v^{\prime}(0))
= v_0^{\prime} (= \tilde{\nu} (\tilde{\lambda}))$ holds;  
$\mathrm{pr}_{ V^{\prime} \times  N, \, V^{\prime}} (v^{\prime}(t)) \in 
D_0^{\prime} + \sum G_i^{\prime}  + \sum I_{ij}^{\prime}$ 
if and only if $t =0$.  
Recall that the configuration corresponding to case $2$-$1$-$1$ was
$v_1^{\prime \prime \prime} \in D_0^{\prime \prime \prime}$,  
$v_0^{\prime} \notin D_0^{\prime} 
+ \sum G_i^{\prime} + \sum I_{ij}^{\prime}$. 
Thus using the holomorphic section 
$v^{\prime} : N \to V^{\prime} \times N$ above, 
and by the same arguement as in 
the proof for case $1$-$2$-$1$,  
we obtain an analytic family 
$\mathrm{pr}_{\mathcal{X}} : \mathcal{X} 
\to N= \{ t \in \mathbb{C} : |t| < \epsilon \}$
satisfying the following conditions: 
$X_0 = {\mathrm{pr}_{\mathcal{X}}}^{-1} (0) \simeq X$; 
for each $t \neq 0 \in N$, the fiber 
$X_t = {\mathrm{pr}_{\mathcal{X}}}^{-1} (t)$ 
is a minimal surface with $c_1^2 = 2\chi -1$, $\chi =4$, 
and $\mathrm{Tors} \simeq \mathbb{Z} /2$ of case $2$-$1$-$1$.
From this together with the results for 
case $2$-$1$-$1$, we infer that $X$ has the same deformation 
type as that of the reference surface $X_{(1)}$,  
and that the point 
$[X]$ lies on the irreducible component $\mathcal{M}_{(1)}$ 
%(or more precisely, on the closure of the image of 
%$\mathbb{P}_0 \to \mathcal{M}$) 
in the proof for case $1$-$1$-$1$.  
%

%%%%%%%%%%%%%%

%
%Case $2$-$2$-$2$: the subcase of case $2$-$2$ 
%where $\bar{A}_0$ is singular. 
%Note that in this case $\bar{A}_0$ has at most negligible 
%singularities. Let $X$ be a minimal surface with 
%$c_1^2 = 2\chi -1$, $\chi =4$, and $\mathrm{Tors} \simeq \mathbb{Z} /2$ 
%of this case. 
%Then we have $h^1 (\mathcal{O}_{\bar{V}} (\bar{A}_0)) = 0$ by 
%Lemma \ref{lm:a0barcoh}. 
%Thus, using the same method as in the proof for 
%case $1$-$1$-$2$,  we obtain an analytic family 
%$\mathrm{pr}_{\mathcal{X}} : \mathcal{X} \to 
%N = \{ u \in \mathbb{C} : |u|  <  \epsilon \}$ 
%satisfying the following conditions:  
%the central fiber $X_0 = {\mathrm{pr}_{\mathcal{X}}}^{-1} (0)$ is 
%isomorphic to $X$; the fiber $X_u = {\mathrm{pr}_{\mathcal{X}}}^{-1} (u)$ 
%over each $u \neq 0 \in N$ is a minimal surface with $c_1^2 = 2\chi -1$, 
%$\chi = 4$, and $\mathrm{Tors} \simeq \mathbb{Z} /2$ of case $2$-$2$-$1$. 
%From this together with the results for case $2$-$2$-$1$, 
%we infer that $X$ has the same deformation type as 
%that of the reference surface $X_{(1)}$. 
%Moreover, we infer that the corresponding point $[X]$ lies on the 
%irreducible component $\mathcal{M}_{(1)}$ 
%(or more precisely, on the closure of the image of 
%$\mathbb{P}_0 \to \mathcal{M}$) given 
%in the proof for case $1$-$1$-$1$. 
%

%

Case $2$-$2$-$2$: the subcase of case $2$-$2$ 
where $\bar{A}_0$ is singular. 
%Note that in this case $\bar{A}_0$ has at most negligible 
%singularities. 
Let $X$ be a minimal surface with 
$c_1^2 = 2\chi -1$, $\chi =4$, and $\mathrm{Tors} \simeq \mathbb{Z} /2$ 
of this case.
%Then using the same argument as in case $1$-$1$-$2$, 
%we infer from the result for case $2$-$2$-$1$ 
%that $X$ has the same deformation type as 
%that of the reference surface $X_{(1)}$. 
%Moreover, we infer that the corresponding point $[X]$ lies on the 
%irreducible component $\mathcal{M}_{(1)}$ 
%%(or more precisely, on the closure of the image of 
%%$\mathbb{P}_0 \to \mathcal{M}$) 
%given in the proof for case $1$-$1$-$1$.
%
Then using the same argument as in case $1$-$1$-$2$, 
we infer from the results for case $2$-$2$-$1$ 
that $X$ has the same deformation type as 
that of the reference surface $X_{(1)}$,  
and that the point $[X]$ lies on the 
irreducible component $\mathcal{M}_{(1)}$ 
%(or more precisely, on the closure of the image of 
%$\mathbb{P}_0 \to \mathcal{M}$) 
in the proof for case $1$-$1$-$1$.

Now that we have the results for all the eight cases 
$1$-$1$-$1$, \ldots , $2$-$2$-$2$, 
we have the assertion.  \qed

Note that from the proof above, we see that 
the morphism $\mathbb{P}_0 \to \mathcal{M}$ 
in the proof for case $1$-$1$-$1$ is dominant.

Now let us prove Theorem \ref{thm:moduli}. 
\medskip

{\sc Proof of Theorem \ref{thm:moduli}}
%\medskip

Let $\mathbb{P}_0 \to \mathcal{M}$   
be the morphism $u \mapsto [X^{(\mu)}_u]$
given in the proof (for case $1$-$1$-$1$) 
of Lemma \ref{lm:irredmoduli}.   
Recall that we have $X^{(\mu_1)}_u \simeq X^{(\mu_2)}_u$ if 
$u \in U^0_{\mu_1} \cap U^0_{\mu_2}$.  
So in what follows, we abbreviate $X^{(\mu)}_u$ to $X_u$. 
Since $\mathbb{P}_0 \to \mathcal{M}$ is a dominant morphism 
from the $29$-dimensional variety $\mathbb{P}_0$,  
we only need to show that  for each $u_0 \in \mathbb{P}_0$, there exist 
at most eight $u \in \mathbb{P}_0$'s satisfying $[X_u] = [X_{u_0}]$. 
Recall also that for all $X$'s of case $1$-$1$-$1$ in the proof of 
Lemma \ref{lm:irredmoduli}, $\tilde{V}$'s have one and the same 
isomorphism class. 
In what follows, we assume that $W$, $Z^{\prime}$, $\tilde{V}$, and 
the configuration of $w_i$'s are those for $X$'s 
of case $1$-$1$-$1$.

Let $\mathrm{Aut} (W)$ be the group of analytic automorphisms  
of $W \simeq \varSigma_0$, and $\iota |_W$, the involution 
of $W$ as in Proposition \ref{prop:sufficientprocedure}.  
Let $\mathrm{Aut} (W, \iota |_W, \{ w_i \} )$ be the subgroup of 
$\mathrm{Aut} (W)$ consisting of  
all $\sigma \in \mathrm{Aut} (W)$'s satisfying 
$(\iota |_W) \circ \sigma = \sigma \circ (\iota |_W)$ and 
$\sigma (\{ w_i\}_{i=1,2}) = \{ w_i\}_{i=1,2}$.   
Since $\mathrm{Aut} (W, \iota |_W, \{ w_i \} )$ acts naturally 
on the sets 
$\{ w_i \}_{i=1, 2}$ and 
$\{\varDelta_0, \varDelta_{\infty}, \varGamma_1, \varGamma_2 \}$, 
we have corresponding group homomorphisms 
$\mathrm{Aut} (W, \iota |_W, \{ w_i \} ) \to \frak{S}_2$ and 
$\mathrm{Aut} (W, \iota |_W, \{ w_i \} ) \to D_4$, 
where $\frak{S}_2$ and $D_4$ denote 
the symmetric group of degree $2$ 
and the dihedral group of degree $4$ 
respectively. 
It is easy to see that the product 
$\mathrm{Aut} (W, \iota |_W, \{ w_i \} ) \to 
\frak{S}_2 \times D_4$ of these two morphisms is an isomorphism. 
 
Let $Z^{\prime} /G$ be the quotient of the surface $Z^{\prime}$
by the natural action by the group $G = \langle \iota |_W \rangle$. 
Then the quotient $Z^{\prime} /G$ has four nodes, and  
the natural morphism 
$\tilde{V} \to Z^{\prime} /G$ gives 
the minimal desingularization of $Z^{\prime} /G$.  
Thus via the diagram 
$\tilde{V} \to Z^{\prime}/G \leftarrow Z^{\prime} \to W$, 
the action by $\mathrm{Aut} (W, \iota |_W, \{ w_i \} )$ on 
the surface $W$ induces one on the pair $(\tilde{V}, \tilde{\lambda})$. 
Let 
$\frak{S}_2 = \langle \iota |_W \rangle 
\to \mathrm{Aut} (W, \iota |_W, \{ w_i \} )$ 
be the natural inclusion.  
Since $\frak{S}_2$ acts trivially on $\tilde{V}$ via this 
inclusion, we obtain a natural action 
by $\mathrm{Aut} (W, \iota |_W, \{ w_i \} ) / \frak{S}_2 \simeq D_4$ 
on the pair $(\tilde{V}, \tilde{\lambda})$. 
This action on $(\tilde{V}, \tilde{\lambda})$ induces one on $\mathbb{P}_0$. 
Remark \ref{rem:onthedescription} however implies that 
we have $X_{u_1} \simeq X_{u_2}$ if and only if 
two points $u_1 \in \mathbb{P}_0$ and $u_2 \in \mathbb{P}_0$ 
belong to the same orbit 
of the action by $\mathrm{Aut} (W, \iota |_W, \{ w_i \} ) / \frak{S}_2$ 
on $\mathbb{P}_0$. 
Thus by $\sharp D_4 =8$, we see that for any $u_0 \in \mathbb{P}_0$ there 
exist at most eight $u \in \mathbb{P}_0$'s satisfying 
$X_{u} \simeq X_{u_0}$. 
Hence we have the assertion.                 \qed 

\section{Appendix}

{\sc Proof of Proposition \ref{prop:deg=n+1}}. 
Let us prove Proposition  \ref{prop:deg=n+1}. 
The method we employ here is the same as the one used in 
\cite[Proof of Lemma 4.5]{evenI}, 
to which we refer the readers for details of the following argument. 
Let $Z \subset \mathbb{P}^n$, where $n \geq 4$, be a  
non-degenerate surface satisfying the assumptions 
in Proposition \ref{prop:deg=n+1}, and $Z^{\prime} \to Z$, 
its minimal desingularization.  
Since we have $\deg Z < 2n -2$, and  
$Z^{\prime} \to Z$ is given by a complete linear system $|D^{\prime}|$, 
the surface $Z^{\prime}$ is a rational surface 
not isomorphic to $\mathbb{P}^2$. 
Thus, for an integer $d$, the surface 
$Z^{\prime}$ admits a birational morphism 
$r : Z^{\prime} \to \varSigma_d  = Z^{\prime}_0$ 
onto the Hirzebruch surface $\varSigma_d$ of degree $d$. 
Let $D^{\prime}_0$ be a general member of the linear system 
$r_* |D^{\prime}|$, and $\varepsilon^{\prime}_i$'s, the total 
transforms to $Z^{\prime}$ of the $(-1)$-curves appearing at 
the blowings up in $Z^{\prime} \to Z$. Then we have 
$D^{\prime} \sim 
r^* D^{\prime}_0 - \sum_{i=1}^s m_i\varepsilon^{\prime}_i$, 
where $m_i$'s, $s \in \mathbb{Z}$. 

\begin{lemma}
There exists an $r: Z^{\prime} \to Z^{\prime}_0$ as above 
such that for any $i$'s, the equality $m_i =1$ holds.  
\end{lemma}

Proof of Lemma. 
Note that the general member $D^{\prime}$ is a non-singular irreducible curve 
on $Z^{\prime}$. If $h^1 (\mathcal{O}_{D^{\prime}}(D^{\prime})) > 0$, 
then by Clifford's theorem on special divisors, we have 
${D^{\prime}}^2 \geq 2 (h^0 (\mathcal{O}_{D^{\prime}}(D^{\prime}))-1)$, 
which contradicts $n \geq 4$. 
Thus we have $h^1 (\mathcal{O}_{D^{\prime}}(D^{\prime})) = 0$. 
From this together with 
the natural short exact sequence 
$0 \to \mathcal{O}_{Z^{\prime}} \to 
\mathcal{O}_{Z^{\prime}}(D^{\prime}) \to 
\mathcal{O}_{D^{\prime}}(D^{\prime}) \to 0$ 
and the Riemann--Roch theorem, 
we infer 
\[
 \chi (\mathcal{O}_{Z^{\prime}}(D^{\prime}) ) = n+1, 
\]
$D^{\prime} K_{Z^{\prime}} = {D^{\prime}}^2 + 
2 (1 - \chi (\mathcal{O}_{Z^{\prime}}(D^{\prime}) )) = 1-n$, 
and 
$(K_{Z^{\prime}} + D^{\prime}/2 )D^{\prime}= (3-n)/2 <0$. 
Thus by Cone Theorem, we find that if 
$Z^{\prime}$ is not the Hirzebruch surface, then 
there exists a $(-1)$-curve $\varepsilon^{\prime}$ on $Z^{\prime}$ 
satisfying $(K_{Z^{\prime}} + D^{\prime}/2 ) \varepsilon^{\prime} <0$. 
Since $Z^{\prime} \to Z$ contracts no $(-1)$-curve, 
we obtain $D^{\prime} \varepsilon^{\prime} =1$. 
Let $r^{\prime} : Z^{\prime} \to Z^{\prime \prime}$ be the 
blowing-down of $\varepsilon^{\prime}$. 
We put $D^{\prime \prime} = r^{\prime}_* (D^{\prime})$. 
If $Z^{\prime \prime}$ is not the Hirzebruch surface, 
then the same argument as above ensures the existence of a  
$(-1)$-curve $\varepsilon^{\prime \prime}$ on $Z^{\prime \prime}$ 
satisfying $D^{\prime \prime} \varepsilon^{\prime \prime} =1$
(for the detail, see \cite[Lemma 4.4]{evenI}). 
We can repeat the same steps until we obtain the Hirzebruch surface. 
\qed 

In what follows, we assume our $r$ satisfies the condition in the lemma above, 
hence 
$D^{\prime} \sim r^* D^{\prime}_0 - \sum_{i=1}^s \varepsilon^{\prime}_i$. 
We put 
$D^{\prime}_0 \sim a \varDelta_0 + b \varGamma$, 
where if $d=0$, we chose $\varDelta_0$ and $\varGamma$ in 
such a way that $b \geq a$. 
Then by $\chi (\mathcal{O}_{Z^{\prime}}(D^{\prime}) ) = n+1$ and 
${D^{\prime}}^2 = n+1$, 
we obtain the following three equalities:  
\begin{align}
n+1 &= D^{\prime}( D^{\prime} - K_{Z^{\prime}})/2 +1 
    = (a+1)(b- ad/2) + a -s +1,  \notag \\
n+1 &= 2a (b - ad/2) -s, \label{eql:chidprime} \\
0   &= {D^{\prime}}^2  - \chi (\mathcal{O}_{Z^{\prime}}(D^{\prime}) )
 = (a-1)(b- ad/2) - (a+1).  \label{eql:dprimesquare-chi}
\end{align}
Note that we have $b- ad/2 \geq a$ if $d \neq 1$, 
and that $b - ad/2 \geq a/2$ if $d =1$. 
Thus by (\ref{eql:chidprime}) and (\ref{eql:dprimesquare-chi}), 
we find $a =2$, $b =d +3$, and $s= 11-n$, 
hence $D^{\prime} \sim - K_{Z^{\prime}} + r^* \varGamma$.   
Since $|D^{\prime}_0|$ has  no fixed component, 
we obtain $d \leq 3$. \qed

%%\bibliographystyle{jipsj}
%\bibliographystyle{plain}
%\bibliography{../../papers,../../books}  

\bigskip

\begin{flushright}
\begin{minipage}{25em}
Masaaki Murakami \\ 
University of Bayreuth, 
Lehrstuhl Mathematik VIII \\
Universitaetsstrasse 30, 
D-95447 Bayreuth, Germany\\
\texttt{Masaaki.Murakami@uni-bayreuth.de}  
\end{minipage}
\end{flushright}

\end{document}